\documentclass[a4paper,11pt,reqno]{amsart}

\usepackage{latexsym,dsfont,amssymb,amsmath,amsthm,amscd}

\usepackage[all]{xy}

\chardef\bslash=`\\ % p. 424, TeXbook

\numberwithin{equation}{section}

\newtheorem{theorem}{Theorem}[section]
\newtheorem{corollary}[theorem]{Corollary}
\newtheorem{lemma}[theorem]{Lemma}
\newtheorem{proposition}[theorem]{Proposition}
\theoremstyle{remark}
\newtheorem{remark}[theorem]{Remark}
\newtheorem{example}[theorem]{Example}
\theoremstyle{definition}

\newcommand\bp{\begin{proof}}
\newcommand\ep{\end{proof}}

\newcommand\Dhat{{\hat\Delta}}

\newcommand\CC{{\mathcal C}}
\newcommand\D{{\mathcal D}}
\newcommand\E{{\mathcal E}}
\newcommand\F{{\mathcal F}}
\newcommand\G{{\mathcal G}}
\newcommand\RR{{\mathcal R}}
\newcommand\U{{\mathcal U}}

\newcommand\bb{{\mathfrak b}}
\newcommand\g{{\mathfrak g}}
\newcommand\h{{\mathfrak h}}

\newcommand\GL{\operatorname{GL}}
\newcommand\End{\operatorname{End}}
\newcommand\Hom{\operatorname{Hom}}

\newcommand\Mat{\operatorname{Mat}}
\newcommand\Mod{\operatorname{-Mod}_f}
\newcommand\Nat{\operatorname{Nat}}

\newcommand\Tr{\operatorname{Tr}}
\newcommand\tr{\operatorname{tr}}
\newcommand\Vect{\mathcal Vec}

\newcommand{\ad}{\operatorname{ad}}

\newcommand{\C}{{\mathbb C}}
\newcommand{\DD}{{\mathbb D}}
\newcommand{\N}{{\mathbb N}}
\newcommand{\Q}{{\mathbb Q}}
\newcommand{\R}{{\mathbb R}}

\newcommand\Z{{\mathbb Z}}

\newcommand\eps{\varepsilon}

\newcommand\enu[1]{\smallskip\newline\makebox[5mm][l]{\rm(#1)}}

\begin{document}

\title[Equivalence of categories]
{Notes on the Kazhdan-Lusztig theorem on equivalence of the Drinfeld
category and the category of $U_q\g$-modules}

%\thanks{}

\author[S. Neshveyev]{Sergey Neshveyev}
\address{Department of Mathematics, University of Oslo,
P.O. Box 1053 Blindern, NO-0316 Oslo, Norway.}
\email{sergeyn@math.uio.no}

\author[L. Tuset]{Lars Tuset}
\address{Faculty of Engineering, Oslo University College,
P.O. Box 4 St.~Olavs plass, NO-0130 Oslo, Norway.}
\email{Lars.Tuset@iu.hio.no}

%\thanks{Supported by the Research Council of Norway.}

\date{November 27, 2007}

\begin{abstract}
We discuss the proof of Kazhdan and
Lusztig of the equivalence of the Drinfeld category~$\D(\g,\hbar)$
of $\g$-modules and the category of finite dimensional
$U_q\g$-modules, $q=e^{\pi i\hbar}$, for $\hbar\in\C\setminus\Q^*$.
Aiming at operator algebraists the result is formulated as the
existence for each $\hbar\in i\R$ of a normalized unitary
$2$-cochain $\F$ on the dual $\hat G$ of a compact simple Lie group
$G$ such that the convolution algebra of $G$  with the coproduct
twisted by $\F$ is $*$-isomorphic to the convolution algebra of the
$q$-deformation $G_q$ of $G$, while the coboundary of $\F^{-1}$
coincides with Drinfeld's KZ-associator defined via monodromy of the
Knizhnik-Zamolodchikov equations.
\end{abstract}

\maketitle

\bigskip

%%%%%%%%%%%%%%%%%%%%%%%%%%%%%%%%%
\section*{Introduction}
%%%%%%%%%%%%%%%%%%%%%%%%%%%%%%%%%

One of the most beautiful and important results in quantum groups is
the theorem of Drinfeld~\cite{Dr1,Dr2} stating that the category of
$U_h\g$-modules is equivalent to a category of $\g$-modules with the
usual tensor product but with nontrivial associativity morphisms
defined by the monodromy of the Knizhnik-Zamolodchikov equations
from conformal field theory. In defining the latter category, known
as the Drinfeld category, Drinfeld was inspired by a result of Kohno
which states that the representation of the braid group defined by
the universal $R$-matrix of $U_h\g$ is equivalent to the monodromy
representation of the KZ-equations. Drinfeld proved equivalence of the categories working in the context of quasi-Hopf algebras, which are
generalizations of Hopf algebras and are algebraic counterparts of
monoidal categories with quasi-fiber functors. In this language the result
says that there exists $\F\in (U\g\otimes U\g)[[h]]$ such that the
coproduct $\Dhat_h$ on $U_h\g\cong U\g[[h]]$ is given by
$\Dhat_h=\F\Dhat(\cdot)\F^{-1}$ and that
$$
(\iota\otimes\Dhat)(\F^{-1})(1\otimes \F^{-1})
(\F\otimes 1)(\Dhat\otimes\iota )(\F)
$$
coincides with the element $\Phi_{KZ}$ defining the associativity
morphisms in the Drinfeld category. Drinfeld worked in the formal
deformation setting and gave two different proofs. Another proof of the equivalence of the categories that works for all
irrational complex parameters was given a few years later by Kazhdan
and Lusztig~\cite{KL1,KL2}. Their approach was then used by Etingof
and Kazhdan~\cite{EK1} to solve the problem of existence of
quantization of an arbitrary Lie bialgebra.

The result of Kazhdan and Lusztig can again be formulated in
algebraic terms, that is, there exists an analogue of the twist $\F$
in the analytic setting. In~\cite{NT2} we observed that such an
element can be used to construct a deformation of the Dirac operator
on quantum groups that gives rise to spectral triples. These notes
originated from a desire to understand better properties of $\F$ for
the study of these quantum Dirac operators. Another motivation is
that the result of Kazhdan and Lusztig is not usually formulated in the form we need. Even though the formulation we are using
should be obvious to a careful reader, to refer this
away to a series of papers totaling several hundred pages seems
inappropriate. What makes the situation more complicated is
that Kazhdan and Lusztig prove a more general result allowing
rational deformation parameters, in which case the Drinfeld category
has to be replaced by a category of modules over the affine Lie
algebra $\hat\g$.

The notes are organized as follows.

Section~\ref{s1} contains categorical preliminaries. The main point
is Drinfeld's notion of a quasi-Hopf algebra~\cite{Dr1}. Since the
monoidal categories we are interested in are infinite, one has to
understand the coproduct in the multiplier sense, so we talk about
discrete quasi-Hopf algebras. Modulo this nuance Section~\ref{s1}
contains the standard dictionary between categorical and algebraic
terms: monoidal categories and quasi-bialgebras, equivalence of
categories and isomorphism of quasi-bialgebras up to twisting, weak tensor
functors and comonoids, rigidity and existence of coinverse.

In Section~\ref{s2} we introduce the Drinfeld category
$\D(\g,\hbar)$, $\hbar\in\C\setminus\Q^*$. As mentioned above, it is
the category of finite dimensional $\g$-modules with the usual
tensor product but with nontrivial associativity
morphisms~$\Phi_{KZ}$ defined  via monodromy of the KZ-equations.
Alternatively one can think of the associator $\Phi_{KZ}$ as a
$3$-cocycle on the dual discrete group $\hat G$. We follow
Drinfeld's original argument~\cite{Dr1,Dr2} to prove that
$\D(\g,\hbar)$ is indeed a braided monoidal category. Remark
that by specialization and analytic continuation this can
be deduced directly from the
formal deformation case, which is a bit more convenient to deal
with. The simplifications are however not significant, so to avoid
confusion we work entirely in the analytic setting. Remark also
that there is a somewhat more conceptual proof showing that
$\D(\g,\hbar)$ is the monoidal category which corresponds to a genus zero
modular functor, see e.g.~\cite{BK}. But as everywhere in these
notes we sacrifice generality in favor of a hands-on approach.

In Section~\ref{s3} we formulate the main result, that is, equivalence of
$\D(\g,\hbar)$ and the category $\CC(\g,\hbar)$ of finite
dimensional admissible $U_q\g$-modules, $q=e^{\pi i\hbar}$.
Furthermore, the functor $\D(\g,\hbar)\to\CC(\g,\hbar)$ defining this equivalence can be chosen
such that its composition with the forgetful functor
$\CC(\g,\hbar)\to \Vect$ is naturally isomorphic to the forgetful
functor $\D(\g,\hbar)\to\Vect$. This means that the equivalence can
be expressed in algebraic terms, that is, the corresponding
quasi-bialgebras are isomorphic up to twisting. The proof of this theorem
occupies the remaining part of the paper. In fact, we prove it only
for generic $\hbar$. A simple compactness argument then shows that
the result holds at least for all $\hbar\in i\R$, which is the most
interesting case from the operator algebra point of view.

The actual proof starts in Section~\ref{s4}. Since we want a functor
isomorphic to the forgetful one, we first of all need a tensor
structure on the forgetful functor $\D(\g,\hbar)\to\Vect$. If we
have a module~$M$ representing this functor then to have a weak
tensor structure on the functor is the same thing as having a comonoid
structure on $M$. Clearly, no finite dimensional $\g$-module can
represent the forgetful functor. In Section~\ref{s4} we define a
representing object $M$ in a completion of $\D(\g,\hbar)$. It can be
thought of as an object in an ind-pro-category, but we prefer to think
of it as a topological $\g$-module.

In Section~\ref{s5} we define a comonoid structure on $M$ thus
endowing the functor $\Hom_\g(M,\cdot)$ with a weak tensor
structure. We then check that for generic $\hbar$ we in fact get a
tensor structure. This already implies that Drinfeld's KZ-associator
is a coboundary for generic $\hbar$. It is interesting to note that up to
this point the only properties of $\Phi_{KZ}$ which have been used
are analytic dependence on the parameter $\hbar$ and that the
associator acts trivially on the highest weight subspaces. We end
the section with an algorithm of how to explicitly find $\F$ such
that $\Phi_{KZ}$ is a coboundary of $\F^{-1}$. The word explicit
should however be taken with a grain of salt, as one has to
make choices depending on values of solutions of differential
equations.

In Section~\ref{sequivalence} we show that $U_q\g$ acts by natural
transformations on the functor $\Hom_\g(M,\cdot)$, allowing the latter to be regarded as a functor $\D(\g,\hbar)\to\CC(\g,\hbar)$. We
finally check that this is an equivalence of categories for generic
$\hbar$. Although the idea of the definition of this action of
$U_q\g$ is not difficult to convey, the right normalization of the
maps involved requires an ingenious choice, which is ultimately
dictated by classical identities for hypergeometric functions. This
is by far the most technical part of the proof of Kazhdan and
Lusztig, and here we omit a couple of the
most tedious computations.

\bigskip

\section{Quasi-bialgebras and monoidal categories}\label{s1}

%\subsection{Braided tensor categories}

A monoidal category $\CC$ is a category with a bifunctor
$\otimes\colon\CC\times\CC\to\CC$, $(U,V)\mapsto U\otimes V$, which
is associative up to a natural isomorphism
$$
\alpha\colon (U\otimes V)\otimes W\to U\otimes (V\otimes W)
$$
and has an object which is the unit $\bf 1$ up to natural
isomorphisms
$$
\lambda\colon {\bf 1}\otimes U\to U,\ \ \
\rho\colon U\otimes {\bf 1}\to U,
$$
such that $\lambda =\rho\colon {\bf 1}\otimes {\bf 1}\to {\bf 1}$
and such that the pentagonal diagram
$$
\xymatrix{
(U\otimes (V\otimes W))\otimes X \ar[d]_{\alpha_{1,23,4}}
& ((U\otimes V)\otimes W)\otimes X \ar[l]_{\alpha\otimes\iota}
\ar[r]^{\alpha_{12,3,4}}
& (U\otimes V)\otimes (W\otimes X) \ar[d]^{\alpha_{1,2,34}}\\
U\otimes ((V\otimes W)\otimes X) \ar[rr]^{\iota\otimes\alpha} &  &
U\otimes (V\otimes (W\otimes X))}
$$
and the triangle diagram
$$
\xymatrix{ (U\otimes{\bf 1})\otimes V \ar[rr]^\alpha
\ar[dr]_{\rho\otimes\iota} & & U\otimes({\bf 1}\otimes V)
\ar[dl]^{\iota\otimes\lambda}\\
& U\otimes V &}
$$
commute.

We say that $\CC$ has strict unit if both $\lambda$ and $\rho$ are
the identity morphisms. If also $\alpha$ is the identity, then $\CC$
is called a strict monoidal category.

A braiding in a monoidal category $\CC$ is a natural isomorphism
$\sigma\colon U\otimes V\to V\otimes U$ such that $\lambda\sigma
(U\otimes{\bf 1})=\rho (U\otimes {\bf 1})$ and such that the
hexagonal diagram
$$
\xymatrix{
(V\otimes U)\otimes W \ar[d]_\alpha &
(U\otimes V)\otimes W\ar[r]^\alpha\ar[l]_{\sigma\otimes\iota}
& U\otimes (V\otimes W)\ar[d]^{\sigma_{1,23}}\\
V\otimes (U\otimes W) \ar[r]^{\iota\otimes\sigma} & V\otimes
(W\otimes U) & (V\otimes W)\otimes U\ar[l]_\alpha }
$$
and the same diagram with $\sigma$ replaced by $\sigma^{-1}$ both
commute. 

We say that a category is $\C$-linear if it is abelian, the sets
$\Hom(U,V)$ are vector spaces over $\C$ and composition of morphisms
is bilinear. Of course, when the monoidal category is $\C$-linear
the tensor functor $\otimes$ is required to be bilinear on
morphisms.

A $\C$-linear category is called semisimple if any object is a
finite direct sum of simple objects.

\medskip

%\subsection{Braided tensor functors}

A (weak) quasi-tensor functor between monoidal categories $\CC$ and
$\CC'$ is a functor $F\colon \CC\to\CC'$ together with a (morphism)
isomorphism $F_0\colon{\bf 1}'\to F({\bf 1})$ in $\CC'$ and natural
(morphisms) isomorphisms
$$F_2\colon F(U)\otimes F(V)\to F(U\otimes V).$$
When the categories are braided then $F$ is called braided if the
diagram
$$
\xymatrix{
F(U)\otimes F(V) \ar[r]^{F_2} \ar[d]_{\sigma'}
& F(U\otimes V)\ar[d]^{F(\sigma)}\\
F(V)\otimes F(U) \ar[r]^{F_2} & F(V\otimes U)}
$$
commutes.

A (weak) quasi-tensor functor is called a (weak) tensor functor if
the diagram
\begin{equation}\label{emonodass}
\xymatrix{ (F(U)\otimes F(V))\otimes F(W) \ar[d]_{\alpha'}
\ar[r]^{\ \ F_2\otimes\iota}
 & F(U\otimes V)\otimes F(W)\ar[r]^{F_2}
& F((U\otimes V)\otimes W)\ar[d]^{F(\alpha)}\\
F(U)\otimes (F(V)\otimes F(W)) \ar[r]^{\ \ \iota\otimes F_2} &
F(U)\otimes F(V\otimes W) \ar[r]^{F_2}
& F(U\otimes (V\otimes W))
}
\end{equation}
and the diagrams
$$
\xymatrix{
F({\bf 1})\otimes F(U) \ar[r]^{F_2} &
F({\bf 1}\otimes U)\ar[d]^{F(\lambda)}\\
{\bf 1}'\otimes F(U) \ar[r]^{\lambda'}
\ar[u]^{F_0\otimes\iota} & F(U)}\ \ \ \ \ \
\xymatrix{
F(U)\otimes F({\bf 1})  \ar[r]^{F_2}  &
F(U\otimes{\bf 1})\ar[d]^{F(\rho)}\\
F(U)\otimes {\bf 1}' \ar[r]^{\rho'} \ar[u]^{\iota\otimes F_0}& F(U)}
$$
commute.

We say that a natural isomorphism $\eta\colon F\to G$ between two
(weak) (quasi-)tensor functors $\CC\to\CC'$ is monoidal if the
diagrams
$$
\xymatrix{
F(U)\otimes F(V) \ar[r]^{F_2} \ar[d]_{\eta\otimes\eta} &
F(U\otimes V)\ar[d]^{\eta}\\
G(U)\otimes G(V) \ar[r]^{G_2} & G(U\otimes V)}\ \ \ \ \ \
\xymatrix{
& {\bf 1}'\ar[dr]^{G_0} \ar[dl]_{F_0} &\\
F({\bf 1})\ar[rr]^{\eta} & & G({\bf 1})}
$$
commute.

An equivalence between two monoidal categories is called monoidal if
the functors and the natural isomorphisms defining the equivalence
are monoidal. If the functors are also ($\C$-linear) (braided) then
we speak of a ($\C$-linear) (braided) monoidal equivalence.

According to a theorem of Mac Lane any monoidal category can be
strictified, i.e. it is monoidally equivalent to a strict monoidal
category, and if the category is ($\C$-linear) (braided) then the
equivalence can be chosen to be ($\C$-linear) (braided).  This is
useful for obtaining new identities for morphisms from known ones:
it implies that an identity holds  if it can be proved assuming that
the associativity morphisms are trivial. As is customary we regard
the $\C$-linear monoidal category~$\Vect$ of finite dimensional
vector spaces as strict.

\medskip

%\subsection{Discrete quasi-bialgebras}

Consider now a direct sum
$A=\oplus_{\lambda\in\Lambda}\End(V_\lambda )$ of full matrix
algebras. Define $M(A)$ as the algebraic product
$\prod_{\lambda\in\Lambda}\End(V_\lambda )$. If $B$ is another such
algebra, we say that a homomorphism $\varphi\colon A\to M(B)$ is
nondegenerate if $\varphi(A)B=B$.

Let $A\Mod$ denote the $\C$-linear category of nondegenerate finite
dimensional $A$-modules, so $A\Mod$ is semisimple with simple
objects $\{V_\lambda\}_\lambda$. We would like $A\Mod$ to be
monoidal with tensor product and strict unit  $\C$ defined in the
usual way via nondegenerate homomorphisms
$$\Delta\colon A\to M(A\otimes A)
=\prod_{\lambda,\mu }\End(V_\lambda\otimes V_\mu ),\ \ \
\varepsilon\colon A\to\C,$$
and with associativity morphisms
$(U\otimes V )\otimes W\to U\otimes (V\otimes W )$
given by acting with an element $\Phi\in M(A\otimes A\otimes A)$.
This is indeed the case if and only if $\Phi$ is invertible and
$$
(\varepsilon\otimes\iota)\Delta=\iota
=(\iota\otimes\varepsilon)\Delta, \ \
(\iota\otimes\varepsilon\otimes\iota)\Phi =1\otimes 1,\ \
(\iota\otimes\Delta)\Delta
=\Phi(\Delta\otimes\iota)\Delta(\cdot)\Phi^{-1},
$$
\begin{equation} \label{equasi3}
(\iota\otimes\iota\otimes\Delta)(\Phi )
(\Delta\otimes\iota\otimes\iota)(\Phi )
=(1\otimes\Phi)(\iota\otimes\Delta\otimes\iota)(\Phi)(\Phi\otimes
1).
\end{equation}
We then call $A$ a {\em discrete quasi-bialgebra} with coproduct
$\Delta$, counit $\varepsilon$ and associator $\Phi$. Remark that
equation~\eqref{equasi3} corresponds to the pentagonal diagram.
Notice also that by definition $A\Mod$ is strict if and only if
$\Phi =1\otimes1\otimes 1$.

If we also have an element $\RR\in M(A\otimes A)$ and let
$\Sigma\colon U\otimes V\to V\otimes U$ denote the flip, then
$\Sigma\RR\colon U\otimes V\to V\otimes U$ is a braiding if and only
if $\Delta^{op} =\RR\Delta(\cdot )\RR^{-1}$ and
\begin{equation}
\label{equasi1}
(\Delta\otimes\iota)(\RR)=\Phi_{312}\RR_{13}\Phi^{-1}_{132}
\RR_{23}\Phi, \ \
(\iota\otimes\Delta)(\RR)=\Phi^{-1}_{231}\RR_{13}\Phi_{213}
\RR_{12}\Phi^{-1}.
\end{equation}
In this case we speak of a {\em quasitriangular discrete
quasi-bialgebra} with $R$-matrix $\RR$. Equations~(\ref{equasi1})
correspond to the hexagonal diagrams.

Note that the forgetful functor $F\colon A\Mod\to \Vect$ is a
quasi-tensor functor with $F_0$ and $F_2$ the identity morphisms. It
is a tensor functor if and only if $\Phi=1\otimes1\otimes1$.

\medskip

%\subsection{Monoidal equivalence and equivalent quasi-bialgebras}

By a twist in a (quasitriangular) discrete quasi-bialgebra $A$ we
mean an invertible element $\F$ in $M(A\otimes A)$ such that
$(\eps\otimes\iota)(\F)=(\iota\otimes\eps)(\F)=1$. The twisting
$A_\F$ of $A$ by $\F$ is then the (quasitriangular) discrete
quasi-bialgebra with comultiplication
$\Delta_\F=\F\Delta(\cdot)\F^{-1}$, counit $\eps_\F=\eps$,
associator
$$
\Phi_\F
=(1\otimes \F)(\iota\otimes\Delta )(\F)
\Phi (\Delta\otimes\iota )(\F^{-1})(\F^{-1}\otimes 1)
$$
(and $R$-matrix  $\RR_\F=\F_ {21}\RR\F^{-1}$).

\begin{proposition} \label{QAvsME}
Let $A$ and $A'$ be (quasitriangular) discrete quasi-bialgebras,
$F\colon A\Mod\to\Vect$ and $F\colon A'\Mod\to \Vect$ the forgetful
quasi-tensor functors. Then \enu{i} the (quasitriangular) discrete
quasi-bialgebras $A'$ and $A$ are isomorphic if and only if there
exists a $\C$-linear (braided) monoidal equivalence $E\colon
A\Mod\to A'\Mod$ such that $F'E$ and $F$ are monoidally naturally
isomorphic; \enu{ii} the (quasitriangular) discrete quasi-bialgebra
$A'$ is isomorphic to a twisting $A_\F$ of $A$ if and only if there
exists a $\C$-linear (braided) monoidal equivalence $E\colon
A\Mod\to A'\Mod$ such that $F'E$ and $F$ are naturally isomorphic.
\end{proposition}

If $A$ and $A'$ are finite dimensional and quasi-Hopf (see below)
then one does not need a natural isomorphism of $F'E$ and $F$ in (ii),
that is, $A'$ is isomorphic to a twisting of $A$ if and only if the
categories $A\Mod$ and $A'\Mod$ are $\C$-linear (braided) monoidally
equivalent \cite{EO}. This is no longer true in the infinite dimensional
case~\cite{B}.

\bp[Proof of Proposition~\ref{QAvsME}]
Assume first that  we have an isomorphism $\varphi\colon A'\to
A_\F$. Then by restriction of scalars $\varphi$ gives a functor
$E\colon A\Mod\to A'\Mod$. We make it a tensor functor by letting
$E_0=\iota$ and $E_2=\F^{-1}$. It is easy to see that $E$ is a
$\C$-linear (braided) monoidal equivalence. Furthermore, ignoring
the quasi-tensor structure we have $F'E=F$, and if $\F=1\otimes 1$
then $F'E=F$ as quasi-tensor functors.

\smallskip

Conversely, assume we have a $\C$-linear (braided) monoidal
equivalence $E\colon A\Mod\to A'\Mod$  and a natural isomorphism
$\eta\colon F\to F'E$. The algebra $M(A)$ can be identified with the
algebra $\Nat(F)$ of natural transformations of the forgetful
functor $F$ to itself, and similarly $M(A')=\Nat(F')$. The map $\varphi\colon\Nat(F')\to\Nat(F)$ defined by
$\varphi(a')=\eta^{-1}a'\eta$ is then an isomorphism of algebras.

Identifying $M(A\otimes A)$ with $\Nat(F\otimes F)$, we define
$\F\in M(A\otimes A)$ by the diagram
$$
\xymatrix{U\otimes V \ar[r]^\F
\ar[d]_\eta & U\otimes V\ar[d]^{\eta\otimes\eta}\\
E(U\otimes V) & E(U)\otimes E(V) \ar[l]_{E_2}}
$$
In other words, we have $\F_{U,V}
=(\eta^{-1}_U\otimes\eta^{-1}_V)E_2^{-1}\eta_{U\otimes V}$. The
element $\F$ is clearly invertible. It is easy to see that it has
the property $(\eps\otimes\iota)(\F) =(\iota\otimes\eps)(\F)=1$ if
and only if the maps $E_0,\eta\colon \C\to E(\C)$ coincide. This is
the case if $\eta$ is a monoidal natural isomorphism, and can be
achieved in general by rescaling $\eta$. Furthermore, if $\eta$ is
monoidal then $\F$ is the identity map.

The element $\Delta(a)$ considered as an element of $\Nat(F\otimes
F)$ is given by $\Delta(a)_{U,V}=a_{U\otimes V}$. For $a'\in M(A')$
we then have
\begin{align*}
(\F^{-1}(\varphi\otimes\varphi)\Delta'(a')\F)_{U,V} &=\F^{-1}_{U,V}
(\eta^{-1}_U\otimes\eta^{-1}_V)
a'_{E(U)\otimes E(V)}(\eta_U\otimes\eta_V) \F_{U,V}\\
&= \eta^{-1}_{U\otimes V}E_2 a'_{E(U)\otimes E(V)}
E_2^{-1}\eta_{U\otimes V} =
\eta^{-1}_{U\otimes V}a'_{E(U\otimes V)} \eta_{U\otimes V}\\
&=\varphi(a')_{U\otimes V} =(\Delta\varphi(a'))_{U,V},
\end{align*}
so $(\varphi\otimes\varphi)\Delta'\varphi^{-1}=\Delta_\F$.

The diagram \eqref{emonodass} for the tensor functor $E$ reads as
$$
\Phi'=(\iota\otimes E_2^{-1})E_2^{-1}E(\Phi)E_2(E_2\otimes\iota).
$$
Using that $E_2\otimes\iota\colon (E(U)\otimes E(V))\otimes E(W)\to
E(U\otimes V)\otimes E(W)$ is
$(\eta\otimes\iota)(\F^{-1}\otimes\iota)
(\eta^{-1}\otimes\eta^{-1}\otimes\iota)$, and that $E_2\colon
E(U\otimes V)\otimes E(W)\to E((U\otimes V)\otimes W)$ is
$$
\eta_{(U\otimes V)\otimes W}\F^{-1}_{U\otimes
V,W}(\eta^{-1}_{U\otimes V}\otimes\eta^{-1}_W)
=\eta(\Delta\otimes\iota)(\F^{-1})(\eta^{-1}\otimes\eta^{-1}),
$$
we see that $E_2(E_2\otimes\iota)$ in the expression above equals
$\eta(\Delta\otimes\iota)(\F^{-1})
(\F^{-1}\otimes\iota)(\eta^{-1}\otimes\eta^{-1} \otimes\eta^{-1})$.
Using a similar expression for $(\iota\otimes E_2^{-1})E_2^{-1}$ we
get
$$
\Phi'=(\eta\otimes\eta\otimes\eta) (\iota\otimes
\F)(\iota\otimes\Delta )(\F)
\eta^{-1}E(\Phi)\eta(\Delta\otimes\iota)(\F^{-1})
(\F^{-1}\otimes\iota)(\eta^{-1}\otimes\eta^{-1} \otimes\eta^{-1}).
$$
Since $\eta^{-1}E(\Phi)\eta=\Phi$, this is exactly the equality
$\Phi'=(\varphi^{-1}\otimes \varphi^{-1}\otimes
\varphi^{-1})(\Phi_\F)$.

Finally, if our quasi-bialgebras are quasitriangular and the functor
$E$ is braided, we have a commutative diagram
$$
\xymatrix{U\otimes V\ar[d]_{\Sigma\RR}\ar[r]^{\eta\ \ } & E(U\otimes
V) \ar[d]_{E(\Sigma\RR)} \ar[r]^{E_2^{-1}\ } & E(U)\otimes
E(V)\ar[d]^{\Sigma\RR'} \ar[rr]^{\ \ \ \eta^{-1}\otimes\eta^{-1}}&
& U\otimes V\ar[d]^{\Sigma(\varphi\otimes\varphi)(\RR')}\\
V\otimes U\ar[r]^{\eta\ \ } & E(V\otimes U) \ar[r]^{E_2^{-1}\ }&
E(V)\otimes E(U) \ar[rr]^{\ \ \ \eta^{-1}\otimes\eta^{-1}}& &
V\otimes U}
$$
Therefore $\Sigma(\varphi\otimes\varphi)(\RR')=\F\Sigma\RR\F^{-1}$,
that is, $(\varphi\otimes\varphi)(\RR')=\RR_\F$. \ep

We will be interested in the case when $A'$ is a bialgebra, so
$\Phi'=1\otimes1\otimes 1$.  In this case $F'$ is a tensor functor,
so if $E\colon A\Mod\to A'\Mod$ is a monoidal equivalence then
$F'E\colon A\Mod\to\Vect$ is a tensor functor. Therefore to show
that $A'$ is isomorphic to a twisting of $A$, by part (ii) of the
above proposition, we at least need a tensor functor $A\Mod\to\Vect$
which is naturally isomorphic to the forgetful functor.

We remark the following consequence of the proof of the above
proposition: if $E\colon A\Mod\to\Vect$ is a $\C$-linear functor and $\eta\colon
F\to E$ is a natural isomorphism then there is a one-to-one
correspondence between weak tensor structures on $E$ and elements
$\G\in M(A\otimes A)$ such that
$(\eps\otimes\iota)(\G)=1=(\iota\otimes\eps)(\G)$ and
$$
\Phi(\Delta\otimes\iota)(\G)(\G\otimes1)
=(\iota\otimes\Delta)(\G)(1\otimes\G).
$$
Furthermore, $E$ is a tensor functor if and only if $\G$ is invertible,
and then $\Phi_\F=1\otimes1\otimes1$ with $\F=\G^{-1}$.

To define a tensor structure on a functor isomorphic to the
forgetful one, it is convenient to use the following notion. An
object $M$ in a monoidal category $\CC$ with strict unit is called a
comonoid if it comes with two morphisms
$$
\eps\colon M\to \mathbf1,\ \ \delta\colon M\to M\otimes M
$$
such that $(\eps\otimes\iota)\delta=\iota=(\iota\otimes\eps)\delta$
and $(\iota\otimes\delta)\delta=\alpha(\delta\otimes\iota)\delta$.

\begin{lemma} \label{lcomonoid}
Let $A$ be a discrete quasi-bialgebra, $M$ an object in $A\Mod$.
Then there is a one-to-one correspondence between \enu{i} weak
tensor structures on the functor $\Hom(M,\cdot) \colon A\Mod\to
\Vect$; \enu{ii} comonoid structures on $M$.
\end{lemma}

\bp If $M$ is a comonoid then we define $E_2\colon
\Hom(M,U)\otimes\Hom(M,V)\to\Hom(M,U\otimes V)$ by $f\otimes
g\mapsto (f\otimes g)\delta$ and $E_0\colon \mathbf
1=\C\to\Hom(M,\C)$ by $E_0(1)=\eps$.

\smallskip

Conversely, if the functor $E=\Hom(M,\cdot)$ is endowed with a weak
tensor structure, we define $\delta\colon M\to M\otimes M$ as the
image of $\iota\otimes\iota$ under the map
$$
E_2\colon \Hom(M,M)\otimes\Hom(M,M)\to\Hom(M,M\otimes M),
$$
and $\eps\colon M\to\C$ as the image of $1\in\C$ under the map
$E_0\colon\C\to\Hom(M,\C)$. Using naturality of $E_2$ one checks
that the image of $f\otimes g$ under the map $E_2\colon
\Hom(M,U)\otimes\Hom(M,V)\to\Hom(M,U\otimes V)$ is $(f\otimes
g)\delta$. It is then straightforward to check that the axioms of a
weak tensor functor translate into the defining properties of a comonoid.
\ep

We are of course interested in the case when the functor
$\Hom(M,\cdot)$ is naturally isomorphic to the forgetful one.
Clearly, no such object $M$ exists in $A\Mod$ unless $A$ is finite
dimensional. So one needs to extend the category $A\Mod$ to make the
lemma useful. We do not try to do this in general, as depending on
the situation different extensions might be useful.

Remark that in the finite dimensional case the unique object up to
isomorphism, representing the forgetful functor, is the module
$A$; namely, $\Hom(A,U)\to U$, $f\mapsto f(1)$, is a natural
isomorphism. In this case the lemma and the discussion before it
show that there exists a one-to-one correspondence between comonoid
structures on $A$ and elements $\G\in A\otimes A$ such that
$(\eps\otimes\iota)(\G)=1=(\iota\otimes\eps)(\G)$ and
$\Phi(\Delta\otimes\iota)(\G)(\G\otimes1)
=(\iota\otimes\Delta)(\G)(1\otimes\G)$. Explicitly, given such an
element $\G$ one defines $\delta\colon A\to A\otimes A$ by
$\delta(a)=\Delta(a)\G$.

\medskip

Let $A$ be a (quasitriangular) discrete quasi-bialgebra. By a
$*$-operation on $A$ we mean an antilinear involutive
antihomomorphism $x\mapsto x^*$ on $A$ such that
$\Delta(x^*)=\Delta(x)^*$, $\eps(x^*)=\overline{\eps(x)}$, $\Phi$ is
unitary (and $\RR^*=\RR_{21}$). We also require any element of the
form  $1+x^*x$ to be invertible in $M(A)$, so that $A$ can be completed
to a C$^*$-algebra.

\begin{proposition} \label{Starreduction}
Let $A$ and $A'$ be (quasitriangular) discrete $*$-quasi-bialgebras.
Suppose $A'$ is isomorphic to $A_\E$ for a twist $\E$. Then there
exists a unitary twist $\F$ such that $A'$ and $A_\F$ are
$*$-isomorphic.
\end{proposition}

\bp Let $\varphi\colon A'\to A_\E$ be an isomorphism. Since every
homomorphism of full matrix algebras (with the standard
$*$-operation) is equivalent to a $*$-homomorphism, there exists an
invertible element $u\in M(A)$ such that the homomorphism
$\varphi_u:=u\varphi(\cdot)u^{-1}$ is $*$-preserving. We normalize
$u$ such that $\eps(u)=1$. Then $\E_u=(u\otimes u)\E\Delta(u^{-1})$
is a twist and it is easy to check that $\varphi_u\colon A'\to
A_{\E_u}$ is an isomorphism.

Therefore we may assume that $\varphi$ is $*$-preserving. Consider
the polar decomposition $\E=\F|\E|$. Then $\F$ is a unitary twist
and we claim that $\varphi$ is an isomorphism of discrete $*$-quasi-bialgebras $A'$ and $A_\F$.
As $\varphi\colon A'\to A_\F$ is $*$-preserving, we just have to
check that $A_\E=A_\F$.

Applying the $*$-operation to the identity
$(\varphi\otimes\varphi)\Delta'=\E\Delta\varphi(\cdot)\E^{-1}$, we
get
$$
(\varphi\otimes\varphi)\Delta'=(\E^{-1})^*\Delta\varphi(\cdot)\E^*.
$$
It follows that $\E^*\E$ commutes with the image of $\Delta$, hence
so does $|\E|$. In particular, $\Delta_\E=\Delta_\F$.

Now apply the map $T(x)=(x^*)^{-1}$ to the identity
$(\varphi\otimes\varphi\otimes\varphi)(\Phi')=\Phi_\E$.
As $T$ preserves $\Phi'$ and $\Phi$ by unitarity, we get
$(\varphi\otimes\varphi\otimes\varphi)(\Phi')=\Phi_{T(\E)}$.
Therefore
$$
(\Phi_{|\E|})_\F=\Phi_\E=\Phi_{T(\E)}=(\Phi_{|\E|^{-1}})_\F,
$$
whence $\Phi_{|\E|}=\Phi_{|\E|^{-1}}$ as
$\Delta_{|\E|}=\Delta=\Delta_{|\E|^{-1}}$. Thus
$$
(1\otimes|\E|)(\iota\otimes\Dhat)(|\E|)
\Phi(\Dhat\otimes\iota)(|\E|^{-1})(|\E|^{-1}\otimes1)
=(1\otimes|\E|^{-1})(\iota\otimes\Dhat)(|\E|^{-1})
\Phi(\Dhat\otimes\iota)(|\E|)(|\E|\otimes1).
$$
Since $(\iota\otimes\Dhat)(|\E|)$ and $1\otimes|\E|$, as well as
$|\E|\otimes1$ and $(\Dhat\otimes\iota)(|\E|)$, commute, we can
write
$$
((1\otimes|\E|)(\iota\otimes\Dhat)(|\E|))^2\Phi=
\Phi((|\E|\otimes1)(\Dhat\otimes\iota)(|\E|))^2.
$$
Consequently
$$
(1\otimes|\E|)(\iota\otimes\Dhat)(|\E|)\Phi
=\Phi(|\E|\otimes1)(\Dhat\otimes\iota)(|\E|).
$$
Thus $\Phi_{|\E|}=\Phi$, and using again
$\Delta_{|\E|}=\Delta$ we therefore get $\Phi_\E=(\Phi_{|\E|})_\F=\Phi_\F$.

Finally, assume our quasi-bialgebras are quasitriangular. Applying
the $*$-operation and then the flip to the equality
$(\varphi\otimes\varphi)(\RR')=\E_{21}\RR\E^{-1}$ we get
$$
\E_{21}\RR\E^{-1}=(\E_{21}^*)^{-1}\RR\E^*,
$$
so that $(\E^*\E)_{21}\RR=\RR\E^*\E$, whence $|\E|_{21}\RR=\RR|\E|$,
or in other words, $\RR_{|\E|}=\RR$. It follows that
$\RR_\E=(\RR_{|\E|})_\F=\RR_\F$. \ep

\medskip

%\subsection{Quasi-bialgebras from fiber functors}

We next discuss how the notion of a quasi-bialgebra arises naturally
from the Tannakian formalism. This will essentially not be used
later.

\smallskip

Let $\CC$ be a $\C$-linear monoidal category. A (quasi-)fiber
functor is a (quasi-)tensor exact faithful  $\C$-linear functor
$\CC\to\Vect$.

First one has the following reconstruction result~\cite{M}.

\begin{proposition} \label{reconstruction}
Let $\CC$ be a small $\C$-linear semisimple (braided) monoidal
category with simple strict unit. Suppose we have a quasi-fiber
functor $F\colon\CC\to\Vect$. Then there exists a (quasitriangular)
discrete quasi-bialgebra $A$ and a $\C$-linear (braided) monoidal
equivalence $E\colon \CC\to A\Mod$ such that its composition with
the forgetful functor $A\Mod\to\Vect$ is naturally isomorphic to
$F$.
\end{proposition}

Remark that by Proposition~\ref{QAvsME} such a quasi-bialgebra $A$
is unique up to isomorphism and twisting. We also remark that, as
will be clear from the proof, if $F$ is a fiber functor then $A$ can
be chosen to be a discrete bialgebra.

\begin{proof}[Proof of Proposition~\ref{reconstruction}.]
Let $\{V_\lambda\}_{\lambda\in\Lambda}$ be representatives of
isomorphism classes of the simple objects in $\CC$. Put
$A=\oplus_\lambda\End(F(V_\lambda))$. Then
$M(A)=\prod_\lambda\End(F(V_\lambda))$ can be identified with the
algebra $\Nat(F)$ of natural transformations of $F$. Regarding $F$
as a functor $E\colon\CC\to A\Mod$, we get an equivalence of $\CC$
and $A\Mod$ as $\C$-linear categories, since $E$ is exact and maps
the objects $V_\lambda$ onto all simple objects of $A\Mod$ up to
isomorphism.

Identifying $M(A\otimes A)$ with $\Nat(F\otimes F)$ and considering
$F_2$ as a natural transformation from $F\otimes F$ to
$F(\cdot\otimes \cdot)$,we define $\Delta\colon M(A)\to M(A\otimes
A)$ by $\Delta(a)=F_2^{-1}aF_2$. Define also $\eps\colon M(A)\to\C$
by $\eps(a)=a_{\mathbf 1}\in\End(F(\mathbf 1))=\C$. Finally, define
$\Phi\in M(A\otimes A\otimes A)=\Nat(F\otimes F\otimes F)$ by
$$
\Phi=(\iota\otimes F_2^{-1})F_2^{-1}F(\alpha )F_2(F_2\otimes\iota ).
$$
Then by construction $A$ becomes a discrete quasi-bialgebra and $E$
a monoidal functor.

If $\CC$ has braiding $\sigma$ then define $\RR\in M(A\otimes A)
=\Nat(F\otimes F)$ by $\RR=\Sigma F_2^{-1}F(\sigma)F_2$.
Then $A$ is quasitriangular and $E$ is braided.
\ep

\medskip

%\subsection{Rigidity and quasi-Hopf algebras}

A right (resp. left) dual to an object $U$ in a monoidal category
$\CC$ with strict unit consists of an object $U^\vee$ (resp. $^\vee
U)$ and two morphisms
$$e\colon U^\vee\otimes U\to{\mathbf 1},\ \ \
i\colon {\mathbf 1}\to U\otimes U^\vee ,\ \ \
(\hbox{resp.}\ \ e'\colon U\otimes {^\vee U}\to{\mathbf 1},\ \ \
i'\colon {\bf 1}\to {^\vee U}\otimes U )$$
such that the compositions
$$
U\xrightarrow{i\otimes\iota} (U\otimes U^\vee)\otimes U
\xrightarrow{\alpha}U\otimes (U^\vee\otimes U)
\xrightarrow{\iota\otimes e} U,
$$
$$
U^\vee\xrightarrow{\iota\otimes i} U^\vee\otimes (U\otimes U^\vee)
\xrightarrow{\alpha^{-1}}(U^\vee\otimes U)\otimes U^\vee
\xrightarrow{e\otimes\iota} U^\vee
$$
(respectively,
$$
U\xrightarrow{\iota\otimes i'} U\otimes ({}^\vee U\otimes U)
\xrightarrow{\alpha^{-1}}(U\otimes {}^\vee U)\otimes U
\xrightarrow{e'\otimes\iota} U,
$$
$$
{}^\vee U\xrightarrow{i'\otimes\iota} ({}^\vee U\otimes U)\otimes {}^\vee U
\xrightarrow{\alpha}{}^\vee U\otimes (U\otimes {}^\vee U)
\xrightarrow{\iota\otimes e'} {}^\vee U)
$$
are the identity morphisms. The category $\CC$ is called rigid if
every object has left and right duals.

If $t\in\Hom(U,V)$ then the transpose $t^\vee\colon V^\vee\to
U^\vee$ is defined as the composition
$$
V^\vee\xrightarrow{\iota\otimes i} V^\vee\otimes (U\otimes
U^\vee)\xrightarrow{\alpha^{-1}}(V^\vee\otimes U )\otimes
U^\vee\xrightarrow{(\iota\otimes t)\otimes\iota} (V^\vee\otimes
V)\otimes U^\vee\xrightarrow{e\otimes\iota} U^\vee.
$$
We then have the following identities:
%\begin{equation} \label{ttdual}
$$
(t\otimes\iota )i=(\iota\otimes t^\vee)i\colon {\mathbf 1}\to
V\otimes U^\vee,\ \ \ e(\iota\otimes t)=e(t^\vee\otimes\iota )\colon
V^\vee\otimes U\to{\mathbf 1}.
$$
%\end{equation}
This is not difficult to check directly, but is immediate if the
category is strict, which we may assume by Mac Lane's theorem. Now
if $s\in\Hom(V,W)$, assuming that $\CC$ is strict to simplify
computations, the morphism $t^\vee s^\vee$ is by definition given by
the composition
$$
W^\vee\xrightarrow{\iota\otimes i} W^\vee\otimes V\otimes
V^\vee\xrightarrow{\iota\otimes \iota\otimes t^\vee} W^\vee\otimes
V\otimes U^\vee\xrightarrow{\iota\otimes s\otimes \iota}
W^\vee\otimes W\otimes U^\vee\xrightarrow{e\otimes\iota} U^\vee.
$$
But as $(t\otimes\iota )i=(\iota\otimes t^\vee)i$, this is exactly
the definition of $(st)^\vee$. Therefore $V\mapsto V^\vee$ is a
contravariant functor of $\CC$ onto itself.

Similar arguments show that if
$\tilde{U}^\vee$ is another right dual of $U$ with corresponding
morphisms $\tilde{i}$ and $\tilde{e}$, then $\gamma
=(\tilde{e}\otimes\iota )\alpha^{-1}(\iota\otimes i)\colon
\tilde{U}^\vee\to U^\vee$ has inverse $(e\otimes\iota
)\alpha^{-1}(\iota\otimes\tilde{i})$. Also $\tilde{e}
=e(\gamma\otimes\iota )$ and $\tilde{i}=(\iota\otimes\gamma^{-1})i$.
Therefore right duals are unique up to isomorphism.
Similar statements hold for left duals. Finally note that
$$
(\iota\otimes e)\alpha(i'\otimes\iota )\colon U\to {^\vee(U^\vee)}
$$
is an isomorphism with inverse $(\iota\otimes
e')\alpha(i\otimes\iota )$, and similarly that $({}^\vee U )^\vee$
is isomorphic to $U$.

\smallskip

The category $\Vect$ is rigid with $U^\vee={}^\vee U=U^*$ and the
morphisms $e=e'$ and $i=i'$ (identifying $U^{**}$ with $U$), which
we shall denote by $e_v$ and $i_v$, are given by
$$
e_v\colon U^*\otimes U\to\C,\ f\otimes x\mapsto f(x),\ \ \hbox{and}\
\ i_v\colon \C\to U\otimes U^*,\ 1\mapsto\sum_ix_i\otimes x^i,
$$
where $\{x_i\}_i$ is a basis in $U$ and $\{x^i\}_i$ is the dual
basis in $U^*$. Then $t^\vee$ is the usual dual operator $t^*$.

\smallskip

Suppose we are given a nondegenerate anti-homomorphism $S$ of a
discrete quasi-bialgebra $A$. Then for any $A$-module $U$ we can
define an $A$-module structure on the dual space $U^*$ by
$af=f(S(a)\,\cdot)$. To make $U^*$ a right dual object we look for
morphisms
$$
e\colon U^*\otimes U\to\C,\ \ \ i\colon \C\to U\otimes U^*
$$
in the form $e=e_v(1\otimes \alpha)$ and $i=(\beta\otimes1)i_v$ for
some elements $\alpha,\beta\in M(A)$ (note that if $U$ is simple
then any linear maps $U^*\otimes U\to\C$ and $\C\to U\otimes U^*$
must be of this form). Then the maps $e$ and $i$ are morphisms if
and only if
\begin{equation} \label{edual1}
S(a_{(1)})\alpha a_{(2)} =\varepsilon (a)\alpha, \ \ \ a_{(1)}\beta
S(a_{(2)} )=\varepsilon (a)\beta
\end{equation}
as endomorphisms of $U$, and then $U^*$ is a right dual of $U$
in $A\Mod$ if and only if
\begin{equation} \label{edual2}
S(\Phi^{-1}_1)\alpha\Phi^{-1}_2 \beta S(\Phi^{-1}_3 )=1,\ \ \
\Phi_1\beta S(\Phi_2 )\alpha\Phi_3 =1,
\end{equation}
again as endomorphisms of $U$. If there exists an invertible
anti-homomorphism $S$ and elements $\alpha,\beta\in M(A)$ such that
\eqref{edual1}-\eqref{edual2} are satisfied, then we say that $A$ is
a {\em discrete quasi-Hopf algebra} with coinverse $S$. Then
$U^\vee=U^*$ with action $af=f(S(a)\,\cdot)$ is a right dual of $U$,
and $^\vee U=U^*$ with action $af=f(S^{-1}(a)\,\cdot)$ is a left
dual of $U$ with $e'=e_v(S^{-1}(\alpha)\otimes1)$ and $i'=(1\otimes
S^{-1}(\beta))i_v$.

If $\tilde{S}$ is another coinverse with corresponding elements
$\tilde{\alpha},\tilde{\beta}$, then there exists a unique
invertible $u\in A$ such that $\tilde{S}=uS(\cdot )u^{-1}$ and
$\tilde{\alpha}=u\alpha$, $\tilde{\beta}=\beta u^{-1}$. Conversely,
any $\tilde{S}$ and $\tilde{\alpha},\tilde{\beta}$ defined this way
for an invertible $u$ satisfy the same axioms as $S$ and $\alpha
,\beta$. When $\Phi =1\otimes 1\otimes 1$, then $\alpha$ and $\beta$
are inverses to each other, and setting $u=\beta$ thus gives
$\tilde\alpha=\tilde\beta =1$, so $A$ is a discrete multiplier Hopf
algebra with coinverse
 $\tilde S$ in the sense of~\cite{VD}.

\smallskip

We have explained that if a discrete quasi-bialgebra $A$  has
coinverse then $A\Mod$ is rigid. One has the following
converse~\cite{Wor,U,HO}.

\begin{proposition} \label{RigVSCo}
Let $A$ be a discrete quasi-bialgebra with $A\Mod$ rigid and such that
for every simple module $U$ the dimensions of $U$ and $U^\vee$ as
vector spaces coincide. Then $A$ has coinverse.
\end{proposition}

\bp Recall that by definition
$A=\oplus_{\lambda\in\Lambda}\End(V_\lambda)$. For each $\lambda$
the module $V_\lambda^\vee$ is simple, so there exists a unique
$\bar\lambda\in\Lambda$ such that $V_\lambda^\vee\cong
V_{\bar\lambda}$. Fix a linear isomorphism $\eta_\lambda\colon
V_\lambda^*\to V^\vee_\lambda$, which exists as the spaces
$V_\lambda$ and $V_\lambda^\vee$ by assumption have the same vector space
dimension. Then there exists a unique anti-isomorphism
$S_\lambda\colon\End(V_{\bar\lambda})\to\End(V_\lambda)$ such that
if we define an action of $\End(V_{\bar\lambda})$ on $V_\lambda^*$
by $af=f(S_\lambda(a)\cdot)$, then $\eta_\lambda$ is an
$\End(V_{\bar\lambda})$-module map. Since $V_\lambda\cong ({}^\vee
V_\lambda)^\vee$, the set $\{\bar\lambda\}_{\lambda\in\Lambda}$
coincides with $\Lambda$. Thus our anti-isomorphisms $S_\lambda$
define an anti-isomorphism $S$ of $A$ onto itself such that for each
$\lambda$ the dual module $V_\lambda^\vee$ is isomorphic to
$V_\lambda^*$ with action $af=f(S(a)\cdot)$. As explained above, the
morphisms $e\colon V_\lambda^*\otimes V_\lambda\to\C$ and
$i\colon\C\to V_\lambda\otimes V_\lambda^*$ uniquely determine
$\alpha$ and $\beta$, making $S$ a coinverse.

\smallskip

In more categorical terms the above proof goes as follows. Identify
$M(A)$ with the algebra $\Nat(F)$ of natural transformations of the
forgetful functor $F$. Extend isomorphisms $V_\lambda^*\cong
V_\lambda^\vee$ to a natural isomorphism $\eta$ from the functor
$U\mapsto F(U)^*$ to the functor $U\mapsto F(U^\vee)$. Then $S$,
$\alpha$ and $\beta$ are defined by
$$
S(a)_U=\eta^*(a_{U^\vee})^*(\eta^*)^{-1},\ \ \alpha_U
=(\iota\otimes e)(\iota\otimes\eta\otimes\iota)(i_v \otimes\iota ),\
\ \beta_U =(\iota\otimes e_v
)(\iota\otimes\eta^{-1}\otimes\iota)(i\otimes\iota).
$$
\ep

In the case when $A$ is finite dimensional the assumption on the
dimensions of $U$ and $U^\vee$ is automatically
satisfied~\cite{Scha2}. The following example from~\cite{Scha1} (see
also \cite{Zhu}) shows that this is not the case in general.

\begin{example}
Let $G$ be a discrete group, $B\subset G$ a subgroup such that each
double coset $BgB$ contains finitely many right and left cosets of
$B$. Consider the category $\CC$ of $G$-graded $B$-bimodules
$M=\oplus_{g\in G}M_g$ such that $M_g$ is a finite dimensional
complex vector space for each $g$, and $M_g\ne0$ only for $g$ in
finitely many double cosets of $B$. Define a tensor structure on
$\CC$ by $M\otimes N=M\otimes_{\C[B]} N$. Note that $\C[B]$ is a
unit object. The category is rigid with right and left dual $M^\vee$
given by $M_g^\vee=(M_{g^{-1}})^*$ and $B$-bimodule structure given
by $(b_1fb_2)(x)=f(b_2xb_1)$ for $f\in(M_{g^{-1}})^*$ and $x\in
M_{(b_1gb_2)^{-1}}$. The morphism $i\colon \C[B]\to M\otimes M^\vee$
is defined by $i(e)=\sum_{g\in G/B}\sum_{j\in I_g}x_{g,j}\otimes
x_{g,j}^*$, where $\{x_{g,j}\}_{j\in I_g}$ is a basis in $M_g$ with
dual basis $\{x_{g,j}^*\}_{j\in I_g}$, and $e\colon M^\vee\otimes
M\to\C[B]$ is defined by $e(f\otimes x)=f(x(gh)^{-1})gh$ for $f\in
M^\vee_g$ and $x\in M_h$ when $gh\in B$, and $e(f\otimes x)=0$ when
$gh\notin B$.

The category $\CC$ is in general not semisimple. To define a
semisimple subcategory consider a functor $E$ from $\CC$ to the
category of $B\backslash G$-graded finite dimensional right
$B$-modules defined by $E(M)=\C\otimes_{\C[B]}M$. It is not
difficult to see that $E$ is an equivalence of categories.
Furthermore, using $E$ the simple objects of $\CC$ can be described
as follows: the modules that are supported on a single double coset
$BgB$ (so that $M_h=0$ for $h\notin BgB$), and such that the right
action of $B\cap g^{-1}Bg$ on $E(M)_{Bg}$ is irreducible. Consider
now only those modules in $\CC$ which decompose into simple ones
such that the corresponding action of $B\cap g^{-1}Bg$ factors
through a finite group. Equivalently, we define a semisimple
subcategory $\CC_0$ of $\CC$ consisting of modules $M$ such that the
right action of $B$ on $E(M)$ factors through a finite group. Yet
another equivalent condition is that $xb=(gbg^{-1})x$ for all $g\in
G$, $x\in M_g$ and $b$ in a finite index subgroup of $B$ (where we
use the convention that $(gbg^{-1})x=0$ if $gbg^{-1}\notin B$).
Using the latter characterization we see that $\CC_0$ is closed
under tensor product, and if $M$ is in $\CC_0$ then $M^\vee$ is also in $\CC_0$.

Consider the functor $F\colon\CC_0\to\Vect$ defined by
$F(M)=\C\otimes_{\C[B]}M$. To make it a quasi-fiber functor fix a
set of representatives $R$ for $B\backslash G$. Then
$F(M)\cong\oplus_{g\in R}M_g$. For $g\in G$ denote by $[g]\in R$ the
representative of $Bg$. Then define $F_2$ as the composition of the
canonical isomorphisms
$$
F(M)\otimes F(N)\cong \bigoplus_{g,h\in R} M_{[gh^{-1}]} \otimes
N_h\cong\bigoplus_{g,h\in R} M_{gh^{-1}}\otimes
N_h\cong\bigoplus_{g\in R}(M\otimes N)_g\cong F(M\otimes N),
$$
where in the second step we used the isomorphisms $M_{[gh^{-1}]}\to
M_{gh^{-1}}$ given by $x\mapsto (gh^{-1}[gh^{-1}]^{-1})x$. Thus by
Proposition~\ref{reconstruction} the functor $F\colon\CC_0\to\Vect$
defines a discrete quasi-bialgebra $A$ such that $A\Mod$ is rigid.
Notice now that the dimensions of $F(M)$ and $F(M^\vee)$ can be
different. Indeed, let $D=BgB$ be a double coset, $M=\C[D]$. Then
$M^\vee=\C[D^{-1}]$. We have $\dim F(M)=|B\backslash D|$ and $\dim
F(M^\vee)=|B\backslash D^{-1}|=|D/B|$. A simple example where these
dimensions can be different is the $ax+b$ groups
$G=\begin{pmatrix}\Q^* & \Q\\ 0 & 1\end{pmatrix}$,
$B=\begin{pmatrix}1 & \Z\\ 0 & 1\end{pmatrix}$. So in this case the
discrete quasi-bialgebra $A$ with rigid monoidal category $A\Mod$
fails to be quasi-Hopf.
\end{example}

\bigskip

%%%%%%%%%%%%%%%%%%%%%%%%%%%%%%%%%
\section{The Drinfeld category} \label{s2}
%%%%%%%%%%%%%%%%%%%%%%%%%%%%%%%%%

Let $G$ be a simply connected simple compact Lie group, $\g$ its
complexified Lie algebra. Consider the tensor category $\CC(\g)$ of
finite dimensional $\g$-modules. For each $\hbar\in\C\setminus\Q$ we
shall introduce new associativity morphisms in $\CC(\g)$ via
monodromy of the Knizhnik-Zamolodchikov equations.

Consider the $\ad$-invariant symmetric form on $\g$ normalized such
that if we choose a maximal torus in $G$ and denote by $\h\subset\g$
be the corresponding Cartan subalgebra, then for the dual form on
$\h^*$ we have $(\alpha,\alpha)=2$ for short roots. In other words,
if $(a_{ij})_{1\le i,j\le r}$ is the Cartan matrix of $\g$, and
$d_1,\dots,d_r$ the coprime positive integers such that
$(d_ia_{ij})_{i,j}$ is symmetric, then $(\alpha_i,\alpha_j)
=d_ia_{ij}$ for a chosen system $\{\alpha_1,\dots,\alpha_r\}$ of
simple roots.  Let $t=\sum_i x_i\otimes x^i\in\g\otimes\g$ be the
element defined by this form, so $\{x_i\}_i$ is a basis in $\g$ and
$\{x^i\}_i$ is the dual basis. Since $t$ is defined by an invariant
form, it is $\g$-invariant, that is, $[t,\Dhat(x)]=0$ for all $x\in
U\g$, where $\Dhat\colon U\g\to U\g\otimes U\g$ is the
comultiplication. Remark also that by definition of $\Dhat$ we have
\begin{equation} \label{etcocycle}
(\Dhat\otimes\iota)(t)=t_{13}+t_{23},\ \
(\iota\otimes\Dhat)(t)=t_{12}+t_{13}.
\end{equation}

Let $V_1,\dots,V_n$ be finite dimensional $\g$-modules. Denote by
$Y_n$ the set of points $(z_1,\dots,z_n)\in\C^n$ such that $z_i\ne
z_j$ for $i\ne j$. The KZ$_n$ equations is the system of
differential equations
$$
\frac{\partial v}{\partial z_i}=\hbar\sum_{j\ne i}
\frac{t_{ij}}{z_i-z_j}v,\ \ \ i=1,\dots,n,
$$
where $v\colon Y_n\to V_1\otimes\dots\otimes V_n$. This system is
consistent in the sense that the differential operators
$\nabla_i=\frac{\partial}{\partial z_i}-\hbar\sum_{j\ne
i}\frac{t_{ij}}{z_i-z_j}$ commute with each other, or equivalently,
they define a flat holomorphic connection on the trivial vector
bundle over $Y_n$ with fiber $V_1\otimes\dots\otimes V_n$. This can
be checked using that $t$ is symmetric and that
$[t_{ij}+t_{jk},t_{ik}]=0$, which follows from \eqref{etcocycle} and
$\g$-invariance of $t$.

The consistency of the KZ$_n$ equations implies that locally for
each $z^0\in Y_n$ and $v_0\in V_1\otimes\dots\otimes V_n$ there
exists a unique holomorphic solution $v$ with $v(z^0)=v_0$. If
$\gamma\colon[0,1]\to Y_n$ is a path starting at $\gamma(0)=z^0$,
then this solution can be analytically continued along $\gamma$. The
map $v_0\mapsto v(\gamma(1))$ defines a linear isomorphism
$M_\gamma$ of $V_1\otimes\dots\otimes V_n$ onto itself. The
monodromy operator $M_\gamma$ depends only on the homotopy class of
$\gamma$. In particular, for each base point $z^0\in Y_n$ we get a
representation of the fundamental group $\pi_1(Y_n;z^0)$ on
$V_1\otimes\dots\otimes V_n$ by monodromy operators. Recall that
$\pi_1(Y_n;z^0)$ is isomorphic to the pure braid group $PB_n$,
which is the kernel of the homomorphism $B_n\to S_n$. If
$V_1=\dots=V_n$ then the monodromy representation extends to the
whole braid group $B_n$; we shall briefly return to this a bit
later.

The new associativity morphism $(V_1\otimes V_2)\otimes V_3\to
V_1\otimes (V_2\otimes V_3)$ will be a certain operator which
appears naturally in computing the monodromy representations for
KZ$_3$, it can be thought of as the monodromy operator from the
asymptotic zone $|z_2-z_1|\ll |z_3- z_1|$ to the zone $|z_3-z_2|\ll
|z_3-z_1|$. To proceed rigorously we need to recall a few facts
about differential equations with regular singularities. Observe
first that if
$$
v(z_1,z_2,z_3)=(z_3-z_1)^{\hbar(t_{12}+t_{23}+t_{13})}
w\left(\frac{z_2-z_1}{z_3-z_1}\right),
$$
then $v$ is a solution of KZ$_3$ if and only if $w$ is a solution of
the equation
\begin{equation} \label{eKZ3'}
w'(z)=\hbar\left(\frac{t_{12}}{z}+\frac{t_{23}}{z-1}\right)w(z),
\end{equation}
which we call the modified KZ$_3$ equation.

\begin{proposition} \label{standasymp}
Let $V$ be a finite dimensional vector space, $z\mapsto
A(z)\in\End(V)$ a holomorphic function on the unit disc $\DD$.
Assume $A(0)$ has no eigenvalues that differ by a nonzero integer.
Then the equation
$$
xG'(x)=A(x)G(x)
$$
for $G\colon(0,1)\to\GL(V)$ has a unique solution such that the
function $H(x)=G(x)x^{-A(0)}$ extends to a holomorphic function on
$\DD$ with value $1$ at $0$.

Furthermore, if $G(\cdot\,;\hbar)$ is an analogous solution of
$xG'(x;\hbar)=\hbar A(x)G(x;\hbar)$, which is well-defined for all
$\hbar$ outside the discrete set $\Lambda=\{n(\lambda-\mu)^{-1}\mid
n\in\N,\ \lambda\ne\mu, \ \lambda\ \hbox{and}\ \mu \ \hbox{are
eigenvalues of}\ A(0)\}$, then $H(x;\hbar)=G(x;\hbar)x^{-\hbar
A(0)}$ is analytic on $\DD\times(\C\setminus\Lambda)$.
\end{proposition}

\bp We shall give a proof of this standard result (see
e.g.~\cite{Wa}), mainly to remind how the assumption on $A(0)$
is used.

Write $A(z)=\sum^\infty_{n=0}A_nz^n$. We look for $G(x)$ in the
form $H(x)x^{A_0}$, where $H(x)=\sum^\infty_{n=0}H_nx^n$ with
$H_0=1$. Then $H$ must satisfy the equation
\begin{equation} \label{eGmod}
xH'(x)=A(x)H(x)-H(x)A_0,
\end{equation}
or equivalently, $[A_0,H_n]-nH_n=-\sum_{i=0}^{n-1}A_{n-i}H_i$ for
all $n\ge1$. The operator $\ad_{A_0}-n$ on $\End(V)$ has zero kernel
exactly when $A_0$ has no eigenvalues that differ by~$n$. So by our
assumptions there exist unique $H_n$ satisfying the above
conditions. We then have to check that the series $\sum_nH_nx^n$ is
convergent in the unit disc. Choose $c>0$ such that
$\|(\ad_{A_0}-n)^{-1}\|\le c$ for all $n\ge1$. Define numbers $h_n$
recursively by $h_0=1$, $h_n=c\sum_{i=0}^{n-1}\|A_{n-i}\|h_i$ for
$n\ge1$. We clearly have $\|H_n\|\le h_n$. On the other hand, by
construction the formal power series $h(x)=\sum^\infty_{n=0}h_nx^n$
satisfies the equation $h(x)-1=\varphi(x)h(x)$, where
$\varphi(x)=c\sum_{n\ge1}\|A_n\|x^n$. Since $\varphi$ is analytic on
$\DD$ and $\varphi(0)=0$, we see that $h(x)=(1-\varphi(x))^{-1}$ is
convergent in a neighbourhood of zero. Hence $\sum_nH_nx^n$ is also
convergent in the same neighbourhood. Since a solution of
\eqref{eGmod} can be continued analytically along any path in
$\DD\setminus\{0\}$, we conclude that the convergence must hold on
the whole disc. Furthermore, as $G(x)$ is invertible for small $x$,
it must be invertible everywhere.

Finally, if $A(z;\hbar)$ is analytic in two variables, then the
above argument implies that for any bounded open set $U$
such that the assumption on $A(0;\hbar)$ is satisfied for all
$\hbar\in\overline{U}$, there exists a neighbourhood $W$ of zero such
that the corresponding solution $H(x;\hbar)$ of \eqref{eGmod} with
$A$ replaced by $A(\cdot\,;\hbar)$, is analytic on $W\times U$.
Fixing $x_0\in W\setminus\{0\}$, we can consider $H(\cdot\,;\hbar)$
as a solution of a differential equation depending analytically on a
parameter and with the analytic initial value $H(x_0;\hbar)$ at
$x=x_0$. Hence $H(\cdot\,;\cdot)$ is analytic on $\DD\times U$. \ep

\begin{remark} \label{lSpec}
Uniqueness of $G$ is equivalent to the following
statement: if $A$ is an operator with no eigenvalues that differ by
a nonzero integer, and the function $x\mapsto x^ATx^{-A}$ defined for
positive~$x$ extends to an analytic function in a neighbourhood of
zero with value $1$ at $x=0$, then $T=1$. This is easy to see
directly. More generally, if $x^ATx^{-A}$ extends to an analytic
function in a neighbourhood of zero then $A$ and $T$
commute.\footnote{In the formal deformation setting a similar result
holds without any assumption on the spectrum of~$A$. Namely, if
$x\mapsto x^{hA}Tx^{-hA}\in\Mat_n(\C)[[h]]$ extends analytically,
meaning that every coefficient in the power series extends
analytically, then $A$ and $T$ commute. Indeed, we have
$x^{hA}Tx^{-hA}=T+h[A,T]\log x+\ldots$, which forces $[A,T]=0$.
Moreover, we see that already existence of the limit of
$x^{hA}Tx^{-hA}$ as $x\to0^+$ implies that $A$ and $T$ commute. As a
result replacing analytic functions by formal power series would
simplify some of the subsequent arguments.}
\end{remark}

We will also need a multivariable version of
Proposition~\ref{standasymp}.

\begin{proposition}  \label{standasymp2}
Let $A_1,\dots,A_m\colon \DD^m\to\End(V)$ be analytic functions.
Assume the differential operators
$\nabla_i=z_i\frac{\partial}{\partial z_i}-A_i(z)$, $1\le i\le m$,
pairwise commute. Assume also that none of the operators~$A_i(0)$
has eigenvalues which differ by a nonzero integer. Then the system
of equations
$$
x_i\frac{\partial G}{\partial x_i}(x)=A_i(x)G(x), \ \ 1\le i\le m,
$$
has a unique $\GL(V)$-valued solution on $(0,1)^m$ such that the
function $G(x)x_1^{-A_1(0)}\dots x_m^{-A_m(0)}$ extends to an
analytic function on $\DD^m$ with value $1$ at $x=0$.
\end{proposition}

Remark that the flatness condition $[\nabla_i,\nabla_j]=0$ reads as
$\displaystyle z_i\frac{\partial A_j}{\partial
z_i}-z_j\frac{\partial A_i}{\partial z_j}=[A_i,A_j]$. In particular,
it implies that $[A_i(0),A_j(0)]=0$.

\bp The proposition can be proved by induction on $m$. To simplify
the notation we shall only sketch a proof for $m=2$, which is
actually the only case we shall need later.

The unknown function
$H(x_1,x_2)=G(x_1,x_2)x_1^{-A_1(0)}x_2^{-A_2(0)}$ must satisfy the
system of equations
\begin{align}
x_1\frac{\partial H}{\partial x_1}&=A_1H-HA_1(0),\label{ea1}\\
x_2\frac{\partial H}{\partial x_2}&=A_2H-HA_2(0).\label{ea2}
\end{align}
By the proof of Proposition~\ref{standasymp} equation \eqref{ea1}
for $x_2=0$ has a unique holomorphic solution $H_0$ with $H_0(0)=1$.
Using that $[\nabla_1,\nabla_2]=0$ it is easy to check that
$A_2(\cdot,0)H_0$ is a holomorphic solution of~\eqref{ea1} (for
$x_2=0$) with initial value $A_2(0)$ at $x_1=0$, hence
$A_2(x_1,0)H_0(x_1)=H_0(x_1)A_2(0)$ for all~$x_1$ by uniqueness.
Then an argument similar to that in the proof of
Proposition~\ref{standasymp} shows that in a neighbourhood of zero
there exists a unique holomorphic solution of \eqref{ea2} of the
form $H(x_1,x_2)=\sum_{n=0}^\infty H_n(x_1)x_2^n$, so that
$H(x_1,0)=H_0(x_1)$ for small $x_1$. It remains to show that $H$
also satisfies~\eqref{ea1}. For this one checks, using
$[\nabla_1,\nabla_2]=0$, that
$$
x_1\frac{\partial H}{\partial x_1}-A_1H+HA_1(0)
$$
is again a solution of \eqref{ea2}. Since it is zero at $x_2=0$, we
conclude that it is zero everywhere. \ep

Turning to the modified KZ$_3$ equation \eqref{eKZ3'}, consider more
generally the equation
\begin{equation} \label{eKZ3'01}
w'(z)=\left(\frac{A}{z}+\frac{B}{z-1}\right)w(z),
\end{equation}
where $A$ and $B$ are operators on a finite dimensional vector space
$V$ such that neither $A$ nor $B$ has eigenvalues that differ by a
nonzero integer. By Proposition~\ref{standasymp} there is a unique
$\GL(V)$-valued solution $G_0(x)$ on $(0,1)$ such that
$G_0(x)x^{-A}$ extends to a holomorphic function on $\DD$ with value
$1$ at $0$. Fix $x^0\in(0,1)$. If $w_0\in V$ then $G_0(x)w_0$ is a
solution of \eqref{eKZ3'01} with initial value~$G_0(x^0)w_0$. If we
continue it analytically along a loop $\gamma_0$ starting at $x^0$
and turning around $0$ counterclockwise then at the end point we get
$G_0(x^0)e^{2\pi i A}w_0$. Thus the monodromy operator defined by
$\gamma_0$ is $G_0(x^0)e^{2\pi i A}G_0(x^0)^{-1}$. Using the change of variables $z\mapsto 1-z$ we similarly conclude that there
is a unique $\GL(V)$-valued solution $G_1(x)$ of \eqref{eKZ3'01} such that
$G_1(1-x)x^{-B}$ extends to a holomorphic function on $\DD$ with
value $1$ at $0$. Then the monodromy operator defined by a loop
$\gamma_1$ starting at $x^0$ and turning around $1$ counterclockwise
is $G_1(x^0)e^{2\pi i B}G_1(x^0)^{-1}$. The fundamental group of
$\C\setminus\{0,1\}$ with the base point $x^0$ is freely generated
by the classes $[\gamma_0]$ and $[\gamma_1]$ of $\gamma_0$
and~$\gamma_1$. Therefore the monodromy representation defined by
equation \eqref{eKZ3'01} with the base point $x^0$ is
\begin{equation}\label{emonodr}
[\gamma_0]\mapsto G_0(x^0)e^{2\pi i A}G_0(x^0)^{-1},\ \
[\gamma_1]\mapsto G_1(x^0)e^{2\pi i B}G_1(x^0)^{-1}.
\end{equation}
The operator $\Phi(A,B)=G_1(x)^{-1}G_0(x)$ does not depend on $x$,
since a solution of \eqref{eKZ3'01} is determined by its initial
value. We then see that the above representation is equivalent to
the representation
$$
[\gamma_0]\mapsto e^{2\pi i A},\ \
[\gamma_1]\mapsto \Phi(A,B)^{-1}e^{2\pi i B}\Phi(A,B),
$$
which does not depend on the choice of the base point. In fact it
can be interpreted as the monodromy representation with the base
point $0$ as follows.

Let $\Gamma$ be the space of solutions of \eqref{eKZ3'01} on
$(0,1)$. For each $x^0\in(0,1)$ denote by $\pi_{x^0}\colon
V\to\Gamma$ the isomorphism such that $\pi_{x^0}(w_0)$ is the
solution of \eqref{eKZ3'01} with initial value $w_0$ at $x^0$. If
$\gamma$ is a curve in $(0,1)$ then the monodromy operator
$M_\gamma$ is $\pi_{\gamma(1)}^{-1}\pi_{\gamma(0)}$. Define
$\pi_{x^0}\colon V\to\Gamma$ for $x^0=0,1$ by letting
$\pi_0(w_0)=G_0(\cdot)w_0$ and $\pi_1(w_0)=G_1(\cdot)w_0$. Then
$G_0(x^0)=\pi_{x^0}^{-1}\pi_0$, $G_1(x^0)=\pi_{x^0}^{-1}\pi_1$ and
$\Phi(A,B)=\pi_1^{-1}\pi_0$ can be thought of as the monodromies
from $0$ to $x^0$, from $1$ to $x^0$, and from $0$ to $1$,
respectively. This interpretation agrees with formulas
\eqref{emonodr} since the monodromy operator defined by an
infinitesimal loop around zero should of course be $e^{2\pi i A}$.

It is sometimes convenient to define $\pi_0$ as follows. Let $w_0$
be an eigenvector of $A$ with eigenvalue~$\lambda$. Then
$G_0(x)w_0=x^\lambda G_0(x)x^{-A}w_0$. Therefore $u=\pi_0(w_0)$ is a
solution of \eqref{eKZ3'01} such that $x^{-\lambda}u(x)$
extends to a holomorphic function on $\DD$ with value $w_0$ at $0$.
This completely determines $\pi_0$ if $A$ is diagonalizable.
Similarly, if $w_0$ is an eigenvector of $B$ with eigenvalue
$\lambda$ then $u=\pi_1(w_0)$ is a solution of~\eqref{eKZ3'01} such
that $x^{-\lambda}u(1-x)$ extends to a holomorphic
function on $\DD$ with value $w_0$ at~$0$. \label{MonDef}

\smallskip

We remark the following simple properties of $\Phi(A,B)$: if an
operator $C$ commutes with $A$ and~$B$ then it also commutes with
$\Phi(A,B)$, and in addition $\Phi(A,B)$ coincides with
$\Phi(A+C,B)$ and $\Phi(A,B+C)$ if the latter operators are
well-defined. Indeed, to prove the first claim observe that
$e^{sC}G_0(\cdot)e^{-sC}$ has the defining properties of $G_0$ for
every $s\in\R$, hence it coincides with $G_0$, so $G_0$ commutes
with $C$ and similarly $G_1$ commutes with $C$. For the second claim
observe that if we replace $A$ by $A+C$ then $G_0(x)$ and $G_1(x)$
get replaced by $G_0(x)x^C$ and $G_1(x)x^C$, whence
$\Phi(A+C,B)=\Phi(A,B)$. In particular, if $A$ and $B$ commute then
$\Phi(A,B)=\Phi(0,0)=1$.

Furthermore, by the second part of Proposition~\ref{standasymp} for
any fixed $A$ and $B$ the function $\C\ni\hbar\mapsto\Phi(\hbar
A,\hbar B)$ is well-defined and analytic outside a discrete set.
This discrete set does not contain zero, more precisely, $\Phi(\hbar
A,\hbar B)$ is defined at least for
$|\hbar|<(2\max\{r(A),r(B)\})^{-1}$, where $r$ denotes the spectral
radius. It can be shown~\cite{Dr2} that the first terms of the
Taylor series look like
$$
\Phi(\hbar A,\hbar B)=1-\hbar^2\zeta(2)[A,B]-\hbar^3\zeta(3)
([A,[A,B]]+[B,[A,B]])+\dots,
$$
where $\zeta$ is the Riemann zeta function; see~\cite{Ka,LM} for more
on this expansion.

Finally, if $V$ is a Hermitian vector space and $A^*=-A$, $B^*=-B$,
then $\Phi(A,B)$ is unitary. Indeed, for any $x_0\in(0,1)$ the
function $G_0(\cdot)G_0(x_0)^{-1}$ with value $1$ at $x=x_0$ takes
values in the unitary group, being an integral curve of a
time-dependent right-invariant vector field on this group. Letting
$x_0\to0$ in the equality
$G_0(x)=(G_0(x)G_0(x_0)^{-1})(G_0(x_0)x_0^{-A})x_0^A$, we conclude
that $G_0(x)$ is unitary for any $x\in(0,1)$. We similarly see that
$G_1(x)$ is unitary, and hence $\Phi(A,B)$ is unitary as well.

\smallskip

Returning to the modified KZ$_3$ equation notice first that the
image of the element $t$ in $\End(V_1\otimes V_2)$ has rational
eigenvalues for any finite dimensional $\g$-modules $V_1$ and $V_2$.
To see this, we need to recall that
\begin{equation} \label{eCas1}
t=\frac{1}{2}(\Dhat(C)-1\otimes C-C\otimes 1),
\end{equation}
where $C=\sum_ix_ix^i$ is the Casimir, and that the spectrum of $C$
consists of rational numbers since the image of $C$ under an
irreducible representation with highest weight $\lambda$ is
$(\lambda,\lambda+2\rho)$, where $\rho$ is half the sum of
the positive roots. It follows that for any fixed
$\hbar\in\C\setminus\Q^*$ and all finite dimensional $\g$-modules
$V_1$, $V_2$ and $V_3$ we have a well-defined natural isomorphism
$\Phi(\hbar t_{12},\hbar t_{23})$ of $V=V_1\otimes V_2\otimes V_3$
onto itself. Consider the $\GL(V)$-valued solutions $G_0$ and $G_1$
of \eqref{eKZ3'} as described above. Then
$$
W_i(x_1,x_2,x_3)=(x_3-x_1)^{\hbar(t_{12}+t_{23}+t_{13})}
G_i\left(\frac{x_2-x_1}{x_3-x_1}\right),\ \ i=0,1,
$$
are $\GL(V)$-valued solutions of KZ$_3$ on $\{x_1<x_2<x_3\}$. We
have $\Phi(\hbar t_{12},\hbar t_{23})=W_1(z^0)^{-1}W_0(z^0)$ for any
$z^0=(x^0_1,x^0_2,x^0_3)$. Furthermore, our considerations imply
that $\Phi(\hbar t_{12},\hbar t_{23})$ can be thought of as the
monodromy operator of KZ$_3$ from the asymptotic zone $x_2-x_1\ll
x_3- x_1$ to the zone $x_3-x_2\ll x_3-x_1$, and by conjugating by
$W_0(z^0)^{-1}$ the monodromy operators of KZ$_3$ with the base
point $z^0$ can be written as expressions of $e^{\pi i \hbar t}$ and
$\Phi(\hbar t_{12},\hbar t_{23})$\footnote{ To be precise our
discussion of the monodromy of the modified KZ$_3$ equation is not
quite enough for this conclusion because the additional factor
$(x_3-x_1)^{\hbar(t_{12}+t_{23}+t_{13})}$ has nontrivial
monodromy. In other words, the monodromy of the KZ$_3$ equations
does not reduce completely to that of the modified KZ$_3$ equation.
This is not surprising since the map $Y_3\to\C\setminus\{0,1\}$,
$z=(x_1,x_2,x_3)\mapsto x=\frac{x_2-x_1}{x_3-x_1}$, induces a
surjective homomorphism of the fundamental groups which is however
not injective. Namely, consider the standard generators $g_1$ and
$g_2$ of $B_3$. It is known that $PB_3$ is generated by $g_1^2$,
$g_2^2$ and $g_2g_1^2g_2^{-1}$. For $z^0=(x^0_1,x^0_2,x^0_3)$ with
$x^0_1<x^0_2<x^0_3$, represent $g_i$ by a path~$\tilde\gamma_i$
interchanging $x^0_i$ with $x^0_{i+1}$ such that $x^0_i$ passes
below $x^0_{i+1}$. Then the images of $g_1^2$ and $g_2^2$ in
$\pi_1(\C\setminus\{0,1\};x^0)$ can be represented by the curves
$\gamma_0$ and $\gamma_1$ introduced earlier, so the monodromy
operators of KZ$_3$ corresponding to $g_1^2$ and $g_2^2$ with the
base point $z^0$  are $W_0(z^0)e^{2\pi i\hbar t_{12}}W_0(z^0)^{-1}$
and $W_1(z^0)e^{2\pi i\hbar t_{23}}W_1(z^0)^{-1}$. But we still have
to compute the operator corresponding to $g_2g_1^2g_2^{-1}$.
Consider a more general problem. By embedding $V_1\otimes V_2\otimes
V_3$ into $(V_1\oplus V_2\oplus V_3)^{\otimes 3}$ we may assume
$V_1=V_2=V_3=W$. Extend the representation of $PB_3$ to a
representation of $B_3$ on $V=W^{\otimes 3}$ defined by
$g_1\mapsto\Sigma_{12}M_{\tilde\gamma_1}$ and
$g_2\mapsto\Sigma_{23}M_{\tilde\gamma_2}$. If $x^0$ is the image of
$z^0=\tilde\gamma_1(0)$ in $\C\setminus\{0,1\}$ then the image of
$\tilde\gamma_1(1)$ is $\frac{x^0}{x^0-1}$. It follows that
$M_{\tilde\gamma_1}=W_0(\tilde\gamma(1))W_0(\tilde\gamma(0))^{-1}
=(1-x^0)^{\hbar(t_{12}+t_{23}+t_{13})}G_0\left(\frac{x^0}{x^0-1}\right)
G_0(x^0)^{-1}$. Here $G_0\left(\frac{x^0}{x^0-1}\right)$ is obtained
by analytic continuation of $G_0$ along the image of
$\tilde\gamma_1$, that is, by going through the upper half-plane.
It is not difficult to see that
$\Sigma_{12}(1-x)^{\hbar(t_{12}+t_{23}+t_{13})}
G_0\left(\frac{x}{x-1}\right)\Sigma_{12}=G_0(x)e^{\pi i\hbar
t_{12}}$, by checking that the left hand side is a solution of the
modified KZ$_3$ equation. It follows that
$\Sigma_{12}M_{\tilde\gamma_1}=G_0(x^0)e^{\pi i\hbar
t_{12}}\Sigma_{12}G_0(x^0)^{-1} =W_0(z^0)e^{\pi i\hbar
t_{12}}\Sigma_{12}W_0(z^0)^{-1}$. Similarly one checks that
$\Sigma_{23}M_{\tilde\gamma_2}=W_1(z^0)e^{\pi i\hbar
t_{23}}\Sigma_{23}W_1(z^0)^{-1}$. Thus by conjugating by $W_0(z^0)^{-1}$ we see that the representation of $B_3$ on $V$ is equivalent to
the one given by
$g_1\mapsto e^{\pi i\hbar t_{12}}\Sigma_{12}$, $g_2\mapsto 
\Phi(\hbar t_{12},\hbar t_{23})^{-1}e^{\pi i\hbar t_{23}}\Sigma_{23}
\Phi(\hbar t_{12},\hbar t_{23})$. }, which can be thought of as monodromy operators with the base point at infinity in the
asymptotic zone $x_2-x_1\ll x_3- x_1$.

\smallskip

\begin{theorem} \label{Drinfeldcat}
Let $\hbar\in\C\setminus\Q^*$. Denote by $\D(\g,\hbar)$ the category
of finite dimensional $\g$-modules. Then the standard tensor
product, $\alpha=\Phi(\hbar t_{12},\hbar t_{23})$ and $\sigma=\Sigma
e^{\pi i\hbar t}$ define on $\D(\g,\hbar)$ a structure of a braided
monoidal category.
\end{theorem}

By definition $\D(\g,\hbar)$ is the category of non-degenerate
finite dimensional $\widehat{\C[G]}$-modules,
where~$\widehat{\C[G]}$ is the discrete bialgebra of matrix
coefficients of finite dimensional representations of $G$ with
convolution product, coproduct $\Dhat(g)=g\otimes g$ and counit
$\hat\eps(g)=1$ for $g\in G\subset M(\widehat{\C[G]})$. We can then
reformulate Theorem~\ref{Drinfeldcat} by saying that
$(\widehat{\C[G]},\Dhat,\hat\eps,\Phi(\hbar t_{12},\hbar
t_{23}),e^{\pi i\hbar t})$ is a quasitriangular discrete
quasi-bialgebra. Remark that the algebra $M(\widehat{\C[G]})$ can be
identified with the algebra $\U(G)$ of closed densely defined
operators affiliated with the group von Neumann algebra $W^*(G)$ of
$G$.

The element $\Phi(\hbar t_{12},\hbar t_{23})\in\U(G\times G\times
G)$ is called the Drinfeld associator and is often denoted
by~$\Phi_{KZ}$. Since from now on we are not going to consider any
other associativity morphisms apart from the trivial one and
$\Phi(\hbar t_{12},\hbar t_{23})$, we write $\Phi$ instead of
$\Phi(\hbar t_{12},\hbar t_{23})$ if the value of $\hbar$ is clear
from the context.

\begin{proof}[Proof of Theorem~\ref{Drinfeldcat}.]
The only nontrivial relations that we have to check are
(\ref{equasi3}) and (\ref{equasi1}) with $\RR=e^{\pi i\hbar t}$.

To prove (\ref{equasi3}) consider the system KZ$_4$ in the real
simply connected domain $\{x_1<x_2<x_3<x_4\}$.
Put
$$
T=t_{12}+t_{13}+t_{14}+t_{23}+t_{24}+t_{34}.
$$
Note that $T$ commutes with $t_{ij}$ for all $i$ and $j$. We
consider five solutions of KZ$_4$ in our domain of the form
$(x_4-x_1)^{\hbar T}F(u,v)$, where $u$ and $v$ are certain fractions
of $x_j-x_i$ corresponding to five asymptotic zones. Each asymptotic
zone is associated to a vertex of the pentagon diagram according to
the following rule: if $V_i$ and $V_j$ are between parentheses and
$V_k$ is outside, then $|x_j-x_i|\ll |x_k-x_i|$. E.g. the zone
corresponding to $((V_1\otimes V_2)\otimes V_3)\otimes V_4$ is
$x_2-x_1\ll x_3-x_1\ll x_4-x_1$, and we claim that there exist a
unique $\GL$-valued solution $W_1$ of KZ$_4$ of the form
$$
W_1(x_1,x_2,x_3,x_4)=(x_4-x_1)^{\hbar T}
F_1\left(\frac{x_2-x_1}{x_3-x_1},\frac{x_3-x_1}{x_4-x_1}\right),
$$
and a function $H_1(\cdot,\cdot)$ analytic on $\DD^2$ such that
$H_1(0,0)=1$ and
$$
F_1(u,v)=H_1(u,v)u^{\hbar t_{12}}v^{\hbar(t_{12}+t_{13}+t_{23})}
\ \ \hbox{for}\ u,v\in(0,1).
$$
Indeed, one checks that $F_1$ must satisfy the system of equations
\begin{align}
u\frac{\partial F_1}{\partial u}
&=\hbar\left(t_{12}+\frac{u}{u-1}t_{23}
+\frac{uv}{uv-1}t_{24}\right)F_1,\label{eF1}\\
v\frac{\partial F_1}{\partial v}
&=\hbar\left(t_{12}+t_{13}+t_{23}+\frac{uv}{uv-1}t_{24}
+\frac{v}{v-1}t_{34}\right)F_1.\nonumber
\end{align}
By Proposition~\ref{standasymp2} this system has a unique solution
of the required form. Similarly there exist solutions
$W_2,W_3,W_4,W_5$ such that
\begin{align*}
W_2(x_1,x_2,x_3,x_4)&=(x_4-x_1)^{\hbar T}
F_2\left(\frac{x_3-x_2}{x_3-x_1},\frac{x_3-x_1}{x_4-x_1}\right),\\
%\hbox{ for }x_3-x_2<x_3-x_1<x_4-x_1,\\
W_3(x_1,x_2,x_3,x_4)&=(x_4-x_1)^{\hbar T}
F_3\left(\frac{x_3-x_2}{x_4-x_2},\frac{x_4-x_2}{x_4-x_1}\right),\\
%\hbox{ for }x_3-x_2<x_4-x_2<x_4-x_1,\\
W_4(x_1,x_2,x_3,x_4)&=(x_4-x_1)^{\hbar T}
F_4\left(\frac{x_4-x_3}{x_4-x_2},\frac{x_4-x_2}{x_4-x_1}\right),\\
%\hbox{ for }x_4-x_3<x_4-x_2<x_4-x_1,\\
W_5(x_1,x_2,x_3,x_4)&=(x_4-x_1)^{\hbar T}
F_5\left(\frac{x_2-x_1}{x_4-x_1},\frac{x_4-x_3}{x_4-x_1}\right)
%\hbox{ for }x_2-x_1<x_4-x_1,\ x_4-x_3<x_4-x_1,
\end{align*}
and holomorphic functions $H_i(\cdot,\cdot)$, $i=2,3,4,5$, in a
neighbourhood of zero with $H_i(0,0)=1$ and such that for
positive $u,v$ we have
\begin{align*}
F_2(u,v)&=H_2(u,v)u^{\hbar t_{23}}v^{\hbar(t_{12}+t_{13}+t_{23})},\\
F_3(u,v)&=H_3(u,v)u^{\hbar t_{23}}v^{\hbar(t_{23}+t_{24}+t_{34})},\\
F_4(u,v)&=H_4(u,v)u^{\hbar t_{34}}v^{\hbar(t_{23}+t_{24}+t_{34})},\\
F_5(u,v)&=H_5(u,v)u^{\hbar t_{12}}v^{\hbar t_{34}}.
\end{align*}
Explicitly, one checks that $F_2,F_3,F_4$ and $F_5$ satisfy
\begin{align}
&\begin{cases}u\frac{\partial F_2}{\partial u}=\hbar\left(t_{23}
+\frac{u}{u-1}t_{12}+\frac{uv}{1-v+uv}t_{24}\right)F_2,\\
v\frac{\partial F_2}{\partial v}=\hbar\left(t_{12}+t_{13}+t_{23}
+\frac{uv-v}{1-v+uv}t_{24}+\frac{v}{v-1}t_{34}\right)F_2,
\end{cases}\label{eF2}\\
&\begin{cases}u\frac{\partial F_3}{\partial u}=\hbar\left(t_{23}
+\frac{uv}{uv-v+1}t_{13}+\frac{u}{u-1}t_{34}\right)F_3,\\
v\frac{\partial F_3}{\partial v}=\hbar\left(t_{23}+t_{24}+t_{34}
+\frac{v}{v-1}t_{12}+\frac{v-uv}{v-uv-1}t_{13}\right)F_3,
\end{cases}\label{eF3}\\
&\begin{cases}u\frac{\partial F_4}{\partial u}=\hbar\left(t_{34}
+\frac{uv}{uv-1}t_{13}+\frac{u}{u-1}t_{23}\right)F_4,\nonumber\\
v\frac{\partial F_4}{\partial v}=\hbar\left(t_{23}+t_{24}+t_{34}
+\frac{v}{v-1}t_{12}+\frac{uv}{uv-1}t_{13}\right)F_4,
\end{cases}\\
&\begin{cases}u\frac{\partial F_5}{\partial u}=\hbar\left(t_{12}
+\frac{u}{u+v-1}t_{23}+\frac{u}{u-1}t_{24}\right)F_5,\\
v\frac{\partial F_5}{\partial v}=\hbar\left(t_{34}
+\frac{v}{v-1}t_{13}+\frac{v}{u+v-1}t_{23}\right)F_5.
\end{cases}\nonumber
\end{align}
It turns out that the solutions $W_i$ are related as follows:
\begin{align}
W_1&=W_2(\Phi\otimes1),\label{epent1}\\
W_2&=W_3(\iota\otimes\hat\Delta\otimes\iota)(\Phi)\label{epent2},\\
W_3&=W_4(1\otimes\Phi),\nonumber\\
W_4&=W_5(\iota\otimes\iota\otimes\hat\Delta)(\Phi^{-1}),\nonumber\\
W_5&=W_1(\hat\Delta\otimes\iota\otimes\iota)(\Phi^{-1}),\nonumber
\end{align}
which immediately implies (\ref{equasi3}). We shall only check
(\ref{epent1}) and (\ref{epent2}).

\smallskip

To prove (\ref{epent1}) denote by $\Theta$ the operator such that
$W_1=W_2\Theta$. Then
$$
F_1(u,v)=F_2(1-u,v)\Theta.
$$
For any fixed $u\in(0,1)$ the functions $v\mapsto
F_1(u,v)v^{-\hbar(t_{12}+t_{13}+t_{23})}$ and $v\mapsto
F_2(1-u,v)v^{-\hbar(t_{12}+t_{13}+t_{23})}$ extend analytically to a
neighbourhood of zero. It follows that
$v^{\hbar(t_{12}+t_{13}+t_{23})}\Theta
v^{-\hbar(t_{12}+t_{13}+t_{23})}$ extends analytically as well. By
Remark~\ref{lSpec} this is possible only when $\Theta$ commutes with
$t_{12}+t_{13}+t_{23}$. It follows that
$$
F_1(u,v)v^{-\hbar(t_{12}+t_{13}+t_{23})}
=F_2(1-u,v)v^{-\hbar(t_{12}+t_{13}+t_{23})}\Theta.
$$
Letting $v=0$ in this equality and introducing
$g_1(u)=F_1(u,v)v^{-\hbar(t_{12}+t_{13}+t_{23})}|_{v=0}
=H_1(u,0)u^{\hbar t_{12}}$ and
$g_2(u)=F_2(u,v)v^{-\hbar(t_{12}+t_{13}+t_{23})}|_{v=0}
=H_2(u,0)u^{\hbar t_{23}}$, we then get
$$
g_1(u)=g_2(1-u)\Theta.
$$
Furthermore, letting $v=0$ in (\ref{eF1}) and in the first equation
of (\ref{eF2}), we see that $g_1$ and $g_2$ satisfy
$$
u\frac{dg_1}{du}
=\hbar\left(t_{12}+\frac{u}{u-1}t_{23}\right)g_1,\ \
u\frac{dg_2}{du}=\hbar\left(t_{23}
+\frac{u}{u-1}t_{12}\right)g_2.
$$
The functions $g_1(u)u^{-\hbar t_{12}}=H_1(u,0)$ and
$g_2(u)u^{-\hbar t_{23}}=H_2(u,0)$ extend to analytic functions on
the unit disc with value $1$ at $0$. Thus by definition
$$
\Theta=\Phi(\hbar t_{12},\hbar t_{23})=\Phi\otimes1.
$$

\smallskip

To prove (\ref{epent2}) denote again by $\Theta$ the element such
that $W_2=W_3\Theta$. Then
$$
F_2(u,v)=F_3\left(\frac{uv}{1-v+uv},1-v+uv\right)\Theta.
$$
As in the argument for (\ref{epent1}), but now fixing $v$ instead of
$u$, we first conclude that $\Theta$ commutes with $t_{23}$. Thus
$$
F_2(u,v)u^{-\hbar t_{23}}v^{-\hbar t_{23}}
=F_3\left(\frac{uv}{1-v+uv},1-v+uv\right)
\left(\frac{uv}{1-v+uv}\right)^{-\hbar t_{23}}
(1-v+uv)^{-\hbar t_{23}}\Theta.
$$
So letting $u=0$ and introducing $g_i(v)=F_i(u,v)u^{-\hbar
t_{23}}v^{-\hbar t_{23}}|_{u=0}$ for $i=2,3$, we get
$$
g_2(v)=g_3(1-v)\Theta.
$$
Furthermore, from the second equations in (\ref{eF2}) and
(\ref{eF3}) we obtain
$$
v\frac{dg_2}{dv}=\hbar\left(t_{12}+t_{13}
+\frac{v}{v-1}(t_{24}+t_{34})\right)g_2,\ \
v\frac{dg_3}{dv}=\hbar\left(t_{24}+t_{34}
+\frac{v}{v-1}(t_{12}+t_{13})\right)g_3.
$$
The functions $g_2(v)v^{-\hbar(t_{12}+t_{13})}=H_2(0,v)$ and
$g_3(v)v^{-\hbar(t_{24}+t_{34})}=H_3(0,v)$ extend to analytic
functions in the unit disc with value $1$ at $0$. Therefore
$$
\Theta=\Phi(\hbar t_{12}+\hbar t_{13},\hbar t_{24}+\hbar t_{34}).
$$
As $t_{12}+t_{13}=(\iota\otimes\hat\Delta\otimes\iota)(t_{12})$
and $t_{24}+t_{34}=(\iota\otimes\hat\Delta\otimes\iota)(t_{23})$,
we get
$$
\Theta=(\iota\otimes\hat\Delta\otimes\iota) (\Phi(\hbar t_{12},\hbar
t_{23})) =(\iota\otimes\hat\Delta\otimes\iota)(\Phi).
$$

\medskip

To prove~\eqref{equasi1} observe that the second
relation in (\ref{equasi1}) follows from the first one by flipping
the first and the third factors and using that $t=t_{21}$ and
$\Phi_{321}=\Phi^{-1}$. The latter equality is easily obtained from
the change of variables $z\mapsto 1-z$ in~(\ref{eKZ3'}).

Turning to the proof of the first identity in (\ref{equasi1}),
consider the system KZ$_3$ in the simply connected space
$$
\Gamma=\{(z_1,z_2,z_3)\in Y_3 \mid \Im z_1\le\Im z_2\le\Im z_3\}.
$$
Consider the real domain $\{(x_1,x_2,x_3) \mid x_1<x_2<x_3\}$ and
two $\GL$-valued solutions of KZ$_3$ in this domain of the form
\begin{align*}
W_0(x_1,x_2,x_3)&=(x_3-x_1)^{\hbar(t_{12}+t_{23}+t_{13})}
H_0\left(\frac{x_2-x_1}{x_3-x_1}\right)
\left(\frac{x_2-x_1}{x_3-x_1}\right)^{\hbar t_{12}},\\
W_1(x_1,x_2,x_3)&=(x_3-x_1)^{\hbar(t_{12}+t_{23}+t_{13})}
H_1\left(\frac{x_3-x_2}{x_3-x_1}\right)
\left(\frac{x_3-x_2}{x_3-x_1}\right)^{\hbar t_{23}}.
\end{align*}
Similarly we have solutions of KZ$_3$ in the real domain
$\{(x_1,x_2,x_3) \mid x_1<x_3<x_2\}$ such that
\begin{align*}
W_2(x_1,x_2,x_3)&=(x_2-x_1)^{\hbar(t_{12}+t_{23}+t_{13})}
H_2\left(\frac{x_3-x_1}{x_2-x_1}\right)
\left(\frac{x_3-x_1}{x_2-x_1}\right)^{\hbar t_{13}},\\
W_3(x_1,x_2,x_3)&=(x_2-x_1)^{\hbar(t_{12}+t_{23}+t_{13})}
H_3\left(\frac{x_2-x_3}{x_2-x_1}\right)
\left(\frac{x_2-x_3}{x_2-x_1}\right)^{\hbar t_{23}},
\end{align*}
and solutions in the real domain $\{(x_1,x_2,x_3) \mid x_3<x_1<x_2\}$
such that
\begin{align*}
W_4(x_1,x_2,x_3)&=(x_2-x_3)^{\hbar(t_{12}+t_{23}+t_{13})}
H_4\left(\frac{x_1-x_3}{x_2-x_3}\right)
\left(\frac{x_1-x_3}{x_2-x_3}\right)^{\hbar t_{13}},\\
W_5(x_1,x_2,x_3)&=(x_2-x_3)^{\hbar(t_{12}+t_{23}+t_{13})}
H_5\left(\frac{x_2-x_1}{x_2-x_3}\right)
\left(\frac{x_2-x_1}{x_2-x_3}\right)^{\hbar t_{12}}.
\end{align*}
We require the functions $H_i$ to be analytic
on the unit disc with value $1$ at $0$.
The functions $W_i$ extend uniquely to solutions of
KZ$_3$ on $\Gamma$. By definition of $\Phi$ we immediately have
\begin{equation} \label{eW1}
W_0=W_1\Phi,\ \ W_2=W_3\Phi_{132},\ \ W_4=W_5\Phi_{312}.
\end{equation}

We next compare $W_2$ and $W_4$. Consider the set
$$
\Omega_2=\{(z_1,z_2,z_3)\in\Gamma: |z_3-z_1|<|z_2-z_1|\}.
$$
It has two connected components, $\Omega_2^+$ and $\Omega_2^-$,
corresponding to the two possible orientations of the pair of
vectors $(z_2-z_1,z_3-z_1)$ (if the vectors are colinear, first perturb
$(z_1,z_2,z_3)$ in $\Gamma$). The initial real domain of
definition of $W_2$ is contained in $\Omega_2^+$, so
\begin{equation} \label{eW21}
W_2(z_1,z_2,z_3)=(z_2-z_1)^{\hbar(t_{12}+t_{23}+t_{13})}
H_2\left(\frac{z_3-z_1}{z_2-z_1}\right)
\left(\frac{z_3-z_1}{z_2-z_1}\right)^{\hbar t_{13}}\hbox{ for }
(z_1,z_2,z_3)\in\Omega_2^+.
\end{equation}
Similarly the set $\Omega_4=\{(z_1,z_2,z_3)\in\Gamma:
|z_3-z_1|<|z_2-z_3|\}$ has two connected components $\Omega^+_4$ and
$\Omega^-_4$, with $\Omega^+_4$ containing the initial domain of
definition of $W_4$, and
\begin{equation}\label{eW41}
W_4(z_1,z_2,z_3)=(z_2-z_3)^{\hbar(t_{12}+t_{23}+t_{13})}
H_4\left(\frac{z_1-z_3}{z_2-z_3}\right)
\left(\frac{z_1-z_3}{z_2-z_3}\right)^{\hbar t_{13}}\hbox{ for }
(z_1,z_2,z_3)\in\Omega_4^+.
\end{equation}
In the latter expression $(z_1-z_3)^{\hbar t_{13}}$ means the
function on $\Gamma$ obtained by analytic continuation of
$(x_1-x_3)^{\hbar t_{13}}$ from the real domain $\{x_3<x_1<x_2\}$.
On the other hand, $(z_3-z_1)^{\hbar t_{13}}$ in (\ref{eW21}) is
obtained by analytic continuation from $\{x_1<x_3<x_2\}$. Going from
the first real domain to the second within $\Gamma$ changes the
argument of $x_1-x_3$ by $-\pi$, so that $(x_1-x_3)^{\hbar t_{13}}$
in the second domain is $(x_3-x_1)^{\hbar t_{13}}e^{-i\pi\hbar
t_{13}}$. In other words, we can rewrite (\ref{eW41}) as
\begin{equation}\label{eW42}
W_4(z_1,z_2,z_3)=(z_2-z_3)^{\hbar(t_{12}+t_{23}+t_{13})}
H_4\left(\frac{z_1-z_3}{z_2-z_3}\right)
\left(\frac{z_3-z_1}{z_2-z_3}\right)^{\hbar t_{13}}
e^{-\pi i\hbar t_{13}}
\end{equation}
for $(z_1,z_2,z_3)\in\Omega_4^+$,
and now all the power functions on the right hand sides of
(\ref{eW21}) and (\ref{eW42}) are obtained by analytic continuation
from the real domain $\{x_1<x_3<x_2\}$.

We are now in a position to compute the operator $\Theta$ such that
$W_2=W_4\Theta$. For a real point $(x_1,x_2,x_3)$ such that
$x_1<x_3<x_2$ and $x_3-x_1<x_2-x_3$, which belongs to
$\Omega_2^+\cap\Omega_4^+$, put
$$
x=\frac{x_3-x_1}{x_2-x_1}.
$$
Then by virtue of (\ref{eW21}) and (\ref{eW42}) the equality
$W_2=W_4\Theta$ implies
$$
H_2(x)x^{\hbar t_{13}}=(1-x)^{\hbar(t_{12}+t_{23}+t_{13})}
H_4\left(\frac{x}{x-1}\right)\left(\frac{x}{1-x}\right)^{\hbar
t_{13}} e^{-\pi i\hbar t_{13}}\Theta.
$$
Since $H_2$ and $H_4$ are analytic in a neighbourhood of zero and
$H_2(0)=H_4(0)=1$, we see that the function $x^{\hbar t_{13}}e^{-\pi
i\hbar t_{13}}\Theta x^{-\hbar t_{13}}$ extends to an analytic
function in a neighbourhood of zero with value~$1$ at $0$. By
Remark~\ref{lSpec} this is possible only when $e^{-\pi i\hbar
t_{13}}\Theta=1$.

Similar considerations apply to the pairs $(W_1,W_3)$ and
$(W_0,W_5)$, and we get
\begin{equation} \label{eW0}
W_0=W_5e^{\pi i\hbar (t_{13}+t_{23})},\ \ W_1=W_3e^{\pi i\hbar
t_{23}},\ \ W_2=W_4e^{\pi i\hbar t_{13}}.
\end{equation}
Identities (\ref{eW1}) and (\ref{eW0}) imply
$$
e^{-\pi i\hbar(t_{13}+t_{23})}\Phi_{312} e^{\pi i\hbar
t_{13}}\Phi^{-1}_{132}e^{\pi i\hbar t_{23}}\Phi=1.
$$
As $(\hat\Delta\otimes\iota)(t)=t_{13}+t_{23}$, this is exactly the
first identity in (\ref{equasi1}).
\end{proof}

\bigskip

\section{Theorem of Kazhdan and Lusztig} \label{s3}

For $q\in\C\setminus\{0\}$ not a root of unity consider the
quantized universal enveloping algebra $U_q\g$. To fix notation
recall that it is generated by elements $E_i$, $F_i$, $K_i$,
$K_i^{-1}$, $1\le i\le r$, satisfying the relations
$$
K_iK_i^{-1}=K_i^{-1}K_i=1,\ \ K_iK_j=K_jK_i,\ \
K_iE_jK_i^{-1}=q_i^{a_{ij}}E_j,\ \
K_iF_jK_i^{-1}=q_i^{-a_{ij}}F_j,
$$
$$
E_iF_j-F_jE_i=\delta_{ij}\frac{K_i-K_i^{-1}}{q_i-q_i^{-1}},
$$
$$
\sum^{1-a_{ij}}_{k=0}(-1)^k\begin{bmatrix}1-a_{ij}\\
k\end{bmatrix}_{q_i} E^k_iE_jE^{1-a_{ij}-k}_i=0,\ \
\sum^{1-a_{ij}}_{k=0}(-1)^k\begin{bmatrix}1-a_{ij}\\
k\end{bmatrix}_{q_i} F^k_iF_jF^{1-a_{ij}-k}_i=0,
$$
where $\displaystyle\begin{bmatrix}m\\
k\end{bmatrix}_{q_i}=\frac{[m]_{q_i}!}{[k]_{q_i}![m-k]_{q_i}!}$,
$[m]_{q_i}!=[m]_{q_i}[m-1]_{q_i}\dots [1]_{q_i}$,
$\displaystyle[n]_{q_i}=\frac{q_i^n-q_i^{-n}}{q_i-q_i^{-1}}$ and
$q_i=q^{d_i}$. This is a Hopf algebra with coproduct $\Dhat_q$ and
counit $\hat\eps_q$ defined by
$$
\Dhat_q(K_i)=K_i\otimes K_i,\ \
\Dhat_q(E_i)=E_i\otimes1+ K_i\otimes E_i,\ \
\Dhat_q(F_i)=F_i\otimes K_i^{-1}+1\otimes F_i,
$$
$$
\hat\eps_q(E_i)=\hat\eps_q(F_i)=0,\ \ \hat\eps_q(K_i)=1.
$$

If $V$ is a finite dimensional $U_q\g$-module and $\lambda\in
P\subset \h^*$ is an integral weight, denote by $V(\lambda)$ the
space of vectors $v\in V$ of weight $\lambda$, so that
$K_iv=q_i^{\lambda(h_i)}v$ for all $i$, where $h_i\in\h$ is such
that $\alpha_j(h_i)=a_{ij}$. Recall that $V$ is called admissible if
$V=\oplus_{\lambda\in P}V(\lambda)$. Consider the
tensor category of finite dimensional admissible $U_q\g$-modules. It
is a semisimple category with simple objects indexed by dominant
integral weights $\lambda\in P_+$. For each $\lambda\in P_+$ we fix
an irreducible $U_q\g$-module $V^q_\lambda$ with highest weight
$\lambda$. Denote by $\widehat{\C[G_q]}$ the discrete bialgebra
defined by our category, so $\widehat{\C[G_q]}\cong\oplus_{\lambda\in
P_+}\End(V^q_\lambda)$. Denote by $\U(G_q)$ the multiplier algebra
$M(\widehat{\C[G_q]})$.

The discrete bialgebra $\widehat{\C[G_q]}$ is quasitriangular. The
universal $R$-matrix $\RR_\hbar$ depends on the choice of $\hbar\in\C$
such that $q=e^{\pi i\hbar}$. From now on we write $q^x$ instead of
$e^{\pi i\hbar x}$, provided the choice of~$\hbar$ is clear from the
context. The $R$-matrix $\RR_\hbar$ can can be defined by an
explicit formula, see e.g. \cite[Theorem~8.3.9]{CP}, but for us it
will be enough to remember that it is characterized by
$\Dhat^{op}_q=\RR_\hbar\Dhat_q(\cdot)\RR^{-1}_\hbar$ and the
following property. Let $\lambda,\mu\in P_+$. Denote by
$\bar\lambda\in P_+$ the weight $-w_0\lambda$, where $w_0$ is the
longest element in the Weyl group. Then $-\lambda$ is the lowest
weight of $V^q_{\bar\lambda}$, so there exists a nonzero vector
$\zeta^q_{\bar\lambda}\in V^q_{\bar\lambda}(-\lambda)$ such that
$F_i\zeta^q_{\bar\lambda}=0$. Denote also by $\xi^q_\mu$ a highest
weight vector of $V^q_\mu$, so $E_i\xi^q_\mu=0$. Then
\begin{equation} \label{eRmat}
\RR_\hbar(\zeta^q_{\bar\lambda}\otimes\xi^q_\mu)
=q^{-(\lambda,\mu)}\zeta^q_{\bar\lambda}\otimes\xi^q_\mu.
\end{equation}
This indeed characterizes $\RR_\hbar$ since $\xi^q_\mu\otimes
\zeta^q_{\bar\lambda}$ is a cyclic vector in $V^q_\mu\otimes
V^q_{\bar\lambda}$. Notice that there exists $d\in\N$ such that
$d(\lambda,\mu)\in\Z$ for all $\lambda,\mu\in P$. Therefore for each
$q$ we get only finitely many different $R$-matrices $\RR_\hbar$.

Denote by $\CC(\g,\hbar)$ the strict braided monoidal category of
admissible finite dimensional $U_q\g$-modules with braiding defined
by $\RR_\hbar$.

Finally, if $q$ is real then $\widehat{\C[G_q]}$ is a discrete
$*$-bialgebra, with the $*$-operation defined on $U_q\g$ by for
example $K_i^*=K_i$, $E_i^*=F_iK_i$, $F_i^*=K_i^{-1}E_i$.
Furthermore, if $q>0$ then $q=e^{\pi i\hbar}$ for a unique $\hbar\in
i\R$. In this case $\RR_\hbar^*=(\RR_\hbar)_{21}$, so
$(\widehat{\C[G_q]},\Dhat_q,\hat\eps_q,\RR_\hbar)$ is a
quasitriangular discrete $*$-bialgebra.

\medskip

Since the irreducible $U_q\g$-modules and $\g$-modules are both
parameterized by dominant integral weights, we have a canonical
isomorphism between the centers of $\widehat{\C[G]}$ and
$\widehat{\C[G_q]}$. We can now formulate the main result.

\begin{theorem} \label{Dr}
Let $q>0$ and $\hbar\in i\R$ be such that $q=e^{\pi i\hbar}$. Then
there exists a unitary twist $\F\in\U(G\times G)$ such that the
quasitriangular discrete $*$-quasi-bialgebras $(\widehat{\C[G]},
\Dhat,\hat\eps,\Phi(\hbar t_{12},\hbar t_{23}),e^{\pi i\hbar t})_\F$
and $(\widehat{\C[G_q]},\Dhat_q,\hat\eps_q,1,\RR_\hbar)$ are
$*$-isomorphic, via an isomorphism extending the canonical
identification of the centers.
\end{theorem}

We call an element $\F$ in the above theorem a unitary Drinfeld twist.

\medskip

We shall say that a statement holds for generic $\hbar$ if it holds
for $\hbar$ outside a countable set.

\begin{lemma}
Assume a unitary Drinfeld twist exists for generic $\hbar\in i\R$.
Then a unitary Drinfeld twist exists for all $\hbar\in i\R$.
\end{lemma}

\bp It suffices to show that if $\hbar_n\to\hbar\in i\R^*$ and a
unitary Drinfeld twist exists for every $\hbar_n$ then it exists for
$\hbar$.

For each $n$ fix a $*$-isomorphism
$\varphi_n\colon\U(G_{q_n})\to\U(G)$, where $q_n=e^{\pi i\hbar_n} $,
and a unitary Drinfeld twist $\F_n$. By compactness of finite
dimensional unitary groups, passing to a subsequence we may assume
that $\{\F_n\}_n$ converges (in the strong operator topology) to a
unitary $\F\in W^*(G)\bar\otimes W^*(G)$.

Denote the generators of $U_{q_n}\g$ by $E_i(q_n)$, $F_i(q_n)$,
$K_i(q_n)$. Denote also by $\pi^{q_n}_\lambda\colon
U_{q_n}\g\to\End(V^{q_n}_\lambda)$, resp. $\pi_\lambda\colon
U\g\to\End(V_\lambda)$, an irreducible $*$-representation of
$U_{q_n}\g$, resp. $U\g$, with highest weight~$\lambda$. We claim
that the sequences $\{(\pi_\lambda\circ\varphi_n)(E_i(q_n))\}_n$ are
bounded for any $\lambda$. Indeed, since $\varphi_n$ extends the
canonical identification of the centers by assumption, the
representation $\pi_\lambda\circ\varphi_n$ is unitarily equivalent
to $\pi^{q_n}_{\lambda}$. Normalize the scalar product on
$V^{q_n}_{\lambda}$ by requiring that the highest weight vector
$\xi^{q_n}_{\lambda}$ has norm one. Then the scalar products
$$
(\pi^{q_n}_{\lambda}(F_{i_1}(q_n)\ldots
F_{i_k}(q_n))\xi^{q_n}_{\lambda}
,\pi^{q_n}_{\lambda}(F_{j_1}(q_n)\ldots
F_{j_l}(q_n))\xi^{q_n}_{\lambda})
$$
converge to similar scalar products for $q=e^{\pi i\hbar}$, which
can easily be checked by induction on $k+l$ using
$F_i^*=K_i^{-1}E_i$ and the quantum Serre relations. Choose a set of
multiindices $(i_1,\ldots,i_k)$ such that the vectors
$(\pi^{q}_{\lambda}(F_{i_1}(q)\ldots F_{i_k}(q))\xi^q_{\lambda}$
form a basis in $V^q_{\lambda}$. It then follows that the same
expressions for $q_n$ define a basis in $V^{q_n}_{\lambda}$ whenever
$n$ is sufficiently large. By applying the orthonormalization
procedure we obtain an orthonormal basis in $V^{q_n}_{\lambda}$. The
matrix coefficients of $\pi^{q_n}_{\lambda}(E_i(q_n))$ in this basis
are determined by the scalar products
$$
(\pi^{q_n}_{\lambda}(E_i(q_n))
\pi^{q_n}_{\lambda}(F_{i_1}(q_n)\ldots
F_{i_k}(q_n))\xi^{q_n}_{\lambda}
,\pi^{q_n}_{\lambda}(F_{j_1}(q_n)\ldots
F_{j_l}(q_n))\xi^{q_n}_{\lambda}).
$$
It follows that they converge to the corresponding matrix
coefficients of $\pi^q_{\lambda}(E_i(q))$. In particular, the
sequence $\{\pi^{q_n}_{\lambda}(E_i(q_n))\}_n$ is bounded, and hence
so is $\{(\pi_\lambda\circ\varphi_n)(E_i(q_n))\}_n$. Similar
arguments apply to the other generators of $U_{q_n}\g$.

By passing to a subsequence, we may therefore assume that the
operators $(\pi_\lambda\circ\varphi_n)(T(q_n))$, where $T(q_n)$ is
any of the generators $E_i(q_n)$, $F_i(q_n)$, $K_i(q_n)$ of
$U_{q_n}\g$, converge for every dominant integral weight $\lambda$.
For each $\lambda$ the  operators we get in the limit define a
$*$-representation $\tilde\pi_\lambda\colon U_q\g\to
\End(V_\lambda)$. It is a representation with highest weight
$\lambda$, so for dimension reasons it must be equivalent to the
irreducible representation with highest weight $\lambda$. The
representations $\tilde\pi_\lambda$ define a $*$-isomorphism
$\varphi\colon\U(G_q)\to\U(G)$. As
$\{(\pi_\lambda\circ\varphi_n)(T(q_n))\}_n$ converges to
$(\pi_\lambda\circ\varphi)(T(q))$ for each generator $T(q_n)$
of~$U_{q_n}\g$, the limit $\F$ of $\{\F_n\}_n$ is a unitary Drinfeld
twist with respect to $\varphi$ (e.g. the identity $\Phi(\hbar
t_{12},\hbar t_{23})_\F=1$ holds because $\Phi(\hbar_nt_{12},\hbar_n
t_{23})\to\Phi(\hbar t_{12},\hbar t_{23})$). \ep

Therefore it suffices to prove Theorem~\ref{Dr} for generic
$\hbar\in i\R$. Furthermore, by Proposition~\ref{Starreduction} it
is enough to show that $(\widehat{\C[G]}, \Dhat,\hat\eps,\Phi(\hbar
t_{12},\hbar t_{23}),e^{\pi i\hbar t})_\F$ and
$(\widehat{\C[G_q]},\Dhat_q,\hat\eps_q,1,\RR_\hbar)$ are isomorphic
for a (not necessarily unitary) twist $\F\in\U(G\times G)$. By
Proposition~\ref{QAvsME}(ii) the existence of such an isomorphism
can be reformulated in categorical terms as follows, where we now
consider complex parameters instead of only purely imaginary ones.

\begin{theorem} \label{KL}
For generic $\hbar\in\C$ and $q=e^{\pi i\hbar}$  there exists a
$\C$-linear braided monoidal equivalence
$F\colon\D(\g,\hbar)\to\CC(\g,\hbar)$ such that $F$ maps an
irreducible $\g$-module with highest weight $\lambda$ onto an
irreducible $U_q\g$-module with highest weight $\lambda$, and the
composition of $F$ with the forgetful functor
$\CC(\g,\hbar)\to\Vect$ is naturally isomorphic to the forgetful
functor $\D(\g,\hbar)\to \Vect$.
\end{theorem}

We will start proving this theorem in the next section. In the
remaining part of this section we want to make a few remarks that
will not be important later.

The result holds for all $\hbar\notin\Q^*$ by \cite{KL1,KL2,EK3}.
Recall that since $U_q\g$ is a Hopf algebra, the category
$\CC(\g,\hbar)$ is rigid, with a right dual to $V$ defined by
$V^\vee=V^*$, $af=f(\hat S_q(a)\,\cdot)$, where $\hat S_q$ is the
coinverse. It follows that $\D(\g,\hbar)$ is a rigid tensor category
as well. Let us show that rigidity for all $\hbar\notin\Q^*$ follows
already from Theorem~\ref{KL}\footnote{as well as from the original
result of Drinfeld in the formal deformation case}; in particular,
$(\widehat{\C[G]}, \Dhat,\hat\eps,\Phi(\hbar t_{12},\hbar
t_{23}),e^{\pi i\hbar t})$ is a discrete quasi-Hopf algebra for all
$\hbar\notin\Q^*$ by Proposition~\ref{RigVSCo}. As we have said,
this result will not be used later, but it is in fact the first step
in extending Theorem~\ref{KL} to all $\hbar\notin\Q^*$.

\smallskip

For an element $\beta=\sum_in_i\alpha_i$ of the root lattice put
$K_\beta=\prod_iK_i^{n_i}\in U_q\g$ and
$h_\beta=\sum_in_id_ih_i\in\h$, so that
$\lambda(h_\beta)=(\lambda,\beta)$. For a finite dimensional
$\g$-module $V$ denote by $d(V)$ the dimension of $V$ and by
$d_q(V)$ the quantity $\Tr(q^{h_{2\rho}})$, where $\rho$ is half the
sum of the positive roots. We use the same notation~$d_q(V)$ for the
quantum dimension $\Tr(K_{2\rho})$ of a module $V$ in $\CC(\g,\hbar)$.

Recall that we denote by $i_v\colon\C\to V\otimes V^*$ and
$e_v\colon V^*\otimes V\to\C$ the standard maps making $V^*$ a right
dual of $V$ in $\Vect$.

\begin{corollary}
Let $\hbar\notin\Q^*$, $q=e^{\pi i\hbar}$, and $V$ be an irreducible
$\g$-module. Then a right dual of~$V$ in $\D(\g,\hbar)$ can be
defined by $V^\vee=V^*$ with the usual $\g$-module structure given
by $Xf=-f(X\cdot)$ for $X\in \g$, and $i_V=i_v$, $\displaystyle
e_V=\frac{d_q(V)}{d(V)}e_v$.
\end{corollary}

\bp We shall only check that the composition
$$
V\xrightarrow{i_V\otimes\iota}(V\otimes V^*)\otimes V
\xrightarrow{\Phi}V\otimes (V^*\otimes V) \xrightarrow{\iota\otimes
e_V}V
$$
is the identity map. By continuity it suffices to prove this for
generic $\hbar$.

Assume $V$ is an irreducible module with highest weight $\lambda$.
The map $e_V$ coincides with the composition
$$
V^*\otimes V\xrightarrow{\Sigma q^t}V\otimes V^*
\xrightarrow{q^{(\lambda+2\rho,\lambda)}d_q(V)\ell_V}\C,
$$
where $\ell_V$ is the unique left inverse of $i_V$ in
$\D(\g,\hbar)$, that is, $\ell_V(v\otimes f)=d(V)^{-1}f(v)$. To see
this one just has to check how both maps act on the one-dimensional
submodule $i_V(\C)$ and then observe that $t$ acts on this submodule
as multiplication by $-(\lambda+2\rho,\lambda)$, which follows from
\eqref{eCas1}.

It follows that we equivalently have to show that the composition
$$
V\xrightarrow{i_V\otimes\iota}(V\otimes V^*)\otimes V
\xrightarrow{\Phi}V\otimes (V^*\otimes V)
\xrightarrow{\iota\otimes\Sigma q^t}V\otimes(V\otimes V^*)
\xrightarrow{\iota\otimes q^{(\lambda+2\rho,\lambda)}d_q(V)\ell_V}V
$$
is the identity map. This computation can be done in the equivalent
strict tensor category $\CC(\g,\hbar)$. In other words, we have to check
that for an irreducible module $V$ with highest weight $\lambda$ in
$\CC(\g,\hbar)$ the composition
\begin{equation} \label{erig}
V\xrightarrow{i'_V\otimes\iota}V\otimes V^*\otimes V
\xrightarrow{\iota\otimes\Sigma\RR_\hbar}V\otimes V\otimes V^*
\xrightarrow{\iota\otimes q^{(\lambda+2\rho,\lambda)}d_q(V)\ell'_V}V
\end{equation}
is the identity map, where $i_V'\colon\C\to V\otimes V^*$ is an
isomorphism onto the submodule with trivial $U_q\g$-action and
$\ell'_V$ is the unique left inverse of $i'_V$.
%This follows from
%\cite[Proposition 5.1]{Dr0}, but we will sketch a proof for the
%reader's convenience.
To show this, first of all notice that as $i'_V$ is unique up to a
scalar, the composition does not depend on the choice of $i'_V$.
Hence we may assume that $i'_V$ is given by the same formula as
$i_v$. Then the left inverse map $\ell'_V$ in $\CC(\g,\hbar)$ is given
by
$$
V\otimes V^*\to \C,\ \ v\otimes f\mapsto d_q(V)^{-1}f(K_{2\rho}v),
$$
as can be checked using that the coinverse $\hat S_q$ has the
property $\hat S^2_q(a)=K_{2\rho}aK_{2\rho}^{-1}$. Computing
composition~\eqref{erig} we are then left to check that $\hat
S_q((\RR_\hbar)_0)K_{2\rho}(\RR_\hbar)_1$ acts on $V$ as
multiplication by $q^{-(\lambda+2\rho,\lambda)}$. As $V$ is
irreducible, we know that $\hat
S_q((\RR_\hbar)_0)K_{2\rho}(\RR_\hbar)_1$ acts as a scalar, so it
suffices to check how it acts on a highest weight vector, which is
easy to compute using the explicit formula for the $R$-matrix. \ep

%Remark that the above proof implies that $\hat
%S_q(\RR_1)\RR_0=K_{2\rho}q^{-C_q}$.

\bigskip

\section{Representing the forgetful functor} \label{s4}

To prove Theorem~\ref{KL} we first of all have to introduce a tensor
structure on the forgetful functor $\D(\g,\hbar)\to\Vect$. The goal is to represent this functor by an object, then by
Lemma~\ref{lcomonoid} a weak tensor structure on the functor is
equivalent to a comonoid structure on the representing object.

It is clear that within $\D(\g,\hbar)$ we do not have a representing
object. If we however allow infinite dimensional modules then there
is an obvious choice, the universal enveloping algebra $U\g$.
Namely, for any $\g$-module $V$ we have a canonical isomorphism
$$
\Hom_\g(U\g,V)\to V,\ \ f\mapsto f(1).
$$
It is however more convenient to consider the Lie algebra $\tilde
\g=\g\oplus\h$. Viewing $\g$-modules as $\tilde \g$-modules (with
the second copy of $\h$ acting trivially), the forgetful functor is
clearly naturally isomorphic to~$\Hom_{\tilde\g}(U\tilde\g,\cdot)$.
Recall that $\tilde\g$ comes with a structure of a Manin triple.
Namely, denote by $\bb_+$ and~$\bb_-$ the Borel subalgebras of $\g$,
and by $\mathfrak n_\pm\subset\bb_\pm$ their nilpotent subalgebras.
Consider $\bb_+$ and $\bb_-$ as Lie subalgebras of $\tilde\g$ via
the embeddings $\eta_\pm\colon \bb_\pm\to \g\oplus\h$,
$\eta_\pm(x)=(x,\pm\bar x)$, where $x\mapsto\bar x$ is the
projection $\g=\mathfrak n_+\oplus \h \oplus\mathfrak n_-\to\h$.
Then  $(\tilde\g,\bb_+,\bb_-)$ is a Manin triple with the symmetric
form on~$\tilde \g$ given by
$\big((x_1,y_1),(x_2,y_2)\big)=(x_1,x_2)-(y_1,y_2)$. Denote by
$\tilde t$ the element of $\tilde\g\otimes\tilde\g$ defined by this
symmetric form.

Identifying $U\bb_+$ with $U\tilde\g\otimes_{U\bb_-}\C$, we consider
$U\bb_+$ as a $\tilde\g$-module, which we denote by $M_+$. Similarly
define $M_-$ as $U\tilde\g\otimes_{U\bb_+}\C$. Then $M=M_+\otimes
M_-$ is isomorphic to $U\tilde\g$ as a $\tilde\g$-module by the
Poincare-Birkhoff-Witt theorem, so $M$ represents the forgetful
functor. We now want to define a
comonoid structure on $M$.

Denote by $1_+$ the canonical cyclic vector of $M_+$. Then there
exists a unique $\tilde\g$-module map $\delta_+\colon M_+\to
M_+\otimes M_+$ such that $1_+\mapsto 1_+\otimes 1_+$. This is
nothing else than the comultiplication $\Dhat\colon U\bb_+\to
U\bb_+\otimes U\bb_+$. In particular, $\delta_+$ is coassociative.
Ignore for the moment that $M_+$ is infinite dimensional and observe
that $\delta_+$ is also coassociative with respect to $\tilde\Phi
=\Phi(\hbar\tilde t_{12},\hbar\tilde t_{23})$, that is,
$(\iota\otimes\delta_+)\delta_+
=\tilde\Phi(\delta_+\otimes\iota)\delta_+$. Indeed, formally it is
enough to check this on the vector $1_+$, and this follows
immediately as $\tilde\Phi$ acts trivially on the vector $1_+\otimes
1_+\otimes 1_+$ since the vector is annihilated by $\tilde t_{12}$
and $\tilde t_{23}$. We thus see that $M_+$ is a comonoid. For
similar reasons $M_-$ is a comonoid. Now we want to define a
comonoid structure on $M=M_+\otimes M_-$,  and there is basically
one way to define a morphism $\delta\colon M\to M\otimes M$ using
$\delta_+$ and $\delta_-$, namely, as the composition
\begin{equation} \label{eEKco}
\begin{split}
M_+&\otimes M_-\xrightarrow{\delta_+\otimes\delta_-}
(M_+\otimes M_+)\otimes (M_-\otimes M_-)
\xrightarrow{(\tilde\Phi\otimes \iota)\tilde\Phi^{-1}_{12,3,4}}
(M_+\otimes (M_+\otimes M_-))\otimes M_-\\
&
\xrightarrow{\iota\otimes\Sigma e^{\pi i\hbar\tilde t}\otimes\iota}
(M_+\otimes (M_-\otimes M_+))\otimes M_-
\xrightarrow{\tilde\Phi_{12,3,4}(\tilde\Phi^{-1}\otimes \iota)}
(M_+\otimes M_-)\otimes (M_+\otimes M_-).
\end{split}
\end{equation}
As $M_+$ and $M_-$ are infinite dimensional, it is not obvious how
to make sense of this construction. So our first goal is to find a
representing module which is approximated by finite dimensional
ones.

\medskip

For every dominant integral weight $\mu$ fix an irreducible
$\g$-module $V_\mu$ with highest weight $\mu$. Fix also a highest
weight vector $\xi_\mu\in V_\mu$. We assume that $V_0=\C$ and
$\xi_0=1$. The construction of the representing object is based on
the following standard representation theoretic fact, see
e.g.~\cite{Z}: if $V$ is a finite dimensional $\g$-module and
$\lambda$ an integral weight then the map
\begin{equation} \label{emult}
\Hom_\g(V_{\bar\mu}\otimes V_{\lambda+\mu},V)\to V(\lambda),\ \
f\mapsto f(\zeta_{\bar\mu}\otimes\xi_{\lambda+\mu}),
\end{equation}
is an isomorphism for sufficiently large dominant integral weights
$\mu$, where $V(\lambda)\subset V$ is the subspace of vectors of
weight $\lambda$ and $\zeta_{\bar\mu}$ is a lowest weight vector in
$V_{\bar\mu}$. Remark that the above map is always injective as the
vector $\zeta_{\bar\mu}\otimes\xi_{\lambda+\mu}$ is cyclic.

We need to make a consistent choice of lowest weight vectors. For
this recall that if we fix Chevalley generators $e_i$, $f_i$, $h_i$,
$1\le i\le r$, of $\g$ then for any $\g$-module~$V$ there is an
action of the braid group $B_\g$ associated to $\g$ on $V$, see
e.g.~\cite{L}. Consider the canonical section $W_\g\to B_\g$ and
denote by $\theta\in B_\g$ the transformation corresponding to the
longest element $w_0$ in the Weyl group~$W_\g$. Then $\theta\colon
V\to V$ is a natural isomorphism having the following properties.
If $V$ and $W$ are $\g$-modules then the action of $\theta$ on
$V\otimes W$ coincides with $\theta\otimes\theta$. Next, $\theta$ maps
$V(\lambda)$ onto $V(w_0\lambda)$. In particular, $\theta\xi_\mu$ is
a lowest weight vector in~$V_\mu$, which we denote by $\zeta_\mu$.
Finally, for all $1\le i\le r$ we have $\theta f_i=-e_{\bar
i}\theta$, where $\bar i$ is such that $\alpha_{\bar
i}=\bar\alpha_i=-w_0\alpha_i$.

For an integral weight $\lambda$ and dominant integral weights $\mu$
and $\eta$ such that $\lambda+\mu$ is dominant consider the
composition of morphisms
\begin{equation} \label{etr0}
\tr^{\eta}_{\mu,\lambda+\mu}\colon V_{\bar\mu+\bar\eta}\otimes
V_{\lambda+\mu+\eta} \xrightarrow{T_{\bar\mu,\bar\eta}\otimes
T_{\eta,\lambda+\mu}} V_{\bar\mu}\otimes V_{\bar\eta}\otimes
V_\eta\otimes V_{\lambda+\mu} \xrightarrow{\iota\otimes
S_\eta\otimes\iota}V_{\bar\mu}\otimes V_{\lambda+\mu},
\end{equation}
where the morphisms $T$ and $S$ are uniquely determined by
$$
T_{\mu,\eta}\colon V_{\mu+\eta}\to V_\mu\otimes V_\eta, \ \
\xi_{\mu+\eta}\mapsto\xi_\mu\otimes \xi_\eta,
$$
and
$$
S_\eta\colon V_{\bar\eta}\otimes V_\eta\to\C,\ \
\zeta_{\bar\eta}\otimes\xi_\eta\mapsto1.
$$
Notice that $T_{\mu,\eta}\zeta_{\mu+\eta}=
\zeta_\mu\otimes\zeta_\eta$ by the properties of $\theta$.
It follows that
$$
\tr^{\eta}_{\mu,\lambda+\mu}(\zeta_{\bar\mu+\bar\eta}
\otimes\xi_{\lambda+\mu+\eta})
=\zeta_{\bar\mu}\otimes\xi_{\lambda+\mu},
$$
and this completely determines $\tr^{\eta}_{\mu,\lambda+\mu}$.
Using these morphisms define the inverse limit $\g$-module
$$
M_\lambda=\lim_{\xleftarrow[\mu]{}}V_{\bar\mu}\otimes
V_{\lambda+\mu}.
$$
We consider $M_\lambda$ as a topological $\g$-module with a base of
neighborhoods of zero formed by the kernels of the canonical
morphisms $M_\lambda\to V_{\bar\mu}\otimes V_{\lambda+\mu}$. Observe
that $\tr^{\eta}_{\mu,\lambda+\mu}$ is surjective since its image
contains the cyclic vector
$\zeta_{\bar\mu}\otimes\xi_{\lambda+\mu}$. It follows that the
morphisms $M_\lambda\to V_{\bar\mu}\otimes V_{\lambda+\mu}$ are
surjective. Hence, if $V$ is a $\g$-module with discrete topology,
then any continuous morphism $M_\lambda\to V$ factors through
$V_{\bar\mu}\otimes V_{\lambda+\mu}$ for some~$\mu$, so that the
space $\Hom_\g(M_\lambda,V)$ of such morphisms is the inductive
limit of $\Hom_\g(V_{\bar\mu}\otimes V_{\lambda+\mu},V)$\footnote{
Alternatively one can consider $M_\lambda$ as an object in the
category $\operatorname{pro}$-$\CC(\g)$ obtained by free completion
of $\CC(\g)$ under inverse limits. Then by definition
$\Hom(M_\lambda,V)$ is the inductive limit of
$\Hom_\g(V_{\bar\mu}\otimes V_{\lambda+\mu},V)$.}. In particular,
for any finite dimensional $\g$-module $V$ the maps \eqref{emult}
induce a linear isomorphism
$$
\Hom_\g(M_\lambda,V)\to V(\lambda).
$$
Therefore the topological $\g$-module $M=\oplus_{\lambda\in P} M_\lambda$, where
$P$ is the lattice of integral weights, represents the forgetful
functor.

\smallskip

There is an obvious deficiency in the construction of the module
$M_\lambda$: we did not take into account the associativity
morphisms in the composition \eqref{etr0}. So a more natural
morphism in $\D(\g,\hbar)$ is the composition
$$
V_{\bar\mu+\bar\eta}\otimes V_{\lambda+\mu+\eta}
\to(V_{\bar\mu}\otimes V_{\bar\eta})\otimes (V_\eta\otimes
V_{\lambda+\mu}) \xrightarrow{(\Phi\otimes \iota)\Phi^{-1}_{12,3,4}}
(V_{\bar\mu}\otimes (V_{\bar\eta}\otimes V_\eta))\otimes
V_{\lambda+\mu} \to V_{\bar\mu}\otimes V_{\lambda+\mu}.
$$
Remark that we could instead use
$(\iota\otimes\Phi^{-1})\Phi_{1,2,34}$ as the middle morphism, but
by the coherence theorem of Mac Lane we would get the same composition.

The problem now is that we do not get a coherent system of morphisms
$V_{\bar\mu+\bar\eta}\otimes V_{\lambda+\mu+\eta}\to
V_{\bar\mu}\otimes V_{\lambda+\mu}$. It turns out that this can be
rectified by rescaling. First we need a lemma.

\begin{lemma}
Denote by $g^\hbar_{\mu,\eta}$ the image of
$\zeta_{\bar\mu+\bar\eta}\otimes \xi_{\mu+\eta}$ under the
composition
$$
V_{\bar\mu+\bar\eta}\otimes V_{\mu+\eta}
\xrightarrow{T_{\bar\mu,\bar\eta}\otimes T_{\eta,\mu}}
V_{\bar\mu}\otimes V_{\bar\eta}\otimes V_\eta\otimes V_{\mu}
\xrightarrow{(\iota\otimes S_\eta\otimes\iota)B} V_{\bar\mu}\otimes
V_{\mu}\xrightarrow{S_\mu}\C,
$$
where $B=(\Phi\otimes \iota)\Phi^{-1}_{12,3,4}$. Then for generic $\hbar$ the
map $(\mu,\eta)\mapsto g^\hbar_{\mu,\eta}$ is a $\C^*$-valued
symmetric normalized $2$-cocycle on the semigroup $P_+$ of dominant
integral weights, that is,
$$
g^\hbar_{\mu,\eta}=g^\hbar_{\eta,\mu},\ \
g^\hbar_{0,\eta}=g^\hbar_{\mu,0}=1, \ \
g^\hbar_{\lambda+\mu,\eta}g^\hbar_{\lambda,\mu}
=g^\hbar_{\lambda,\mu+\eta}g^\hbar_{\mu,\eta}.
$$
\end{lemma}

In fact using that $\D(\g,\hbar)$ is rigid one can show that
$g^\hbar_{\mu,\eta}\ne0$ for all $\hbar\notin\Q^*$.

\bp It is easy to see that $g^0_{\mu,\eta}=1$. As
$g^\hbar_{\mu,\eta}$ is analytic in $\hbar$ outside a discrete set,
we conclude that $g^\hbar_{\mu,\eta}\ne0$ for generic $\hbar$.

\smallskip

That $g^\hbar_{0,\eta}=g^\hbar_{\mu,0}=1$ is immediate as the
associator is equal to $1$ as long as one of the modules is trivial.

\smallskip

To show that $g^\hbar_{\mu,\eta}$ is a cocycle first observe that
the compositions
$$
V_{\lambda+\mu+\eta}\xrightarrow{T_{\lambda+\mu,\eta}}
V_{\lambda+\mu}\otimes V_\eta\xrightarrow{T_{\lambda,\mu}\otimes\iota}
V_\lambda\otimes V_\mu\otimes V_\eta\xrightarrow{\Phi}
V_\lambda\otimes V_\mu\otimes V_\eta
$$
and
$$
V_{\lambda+\mu+\eta}\xrightarrow{T_{\lambda,\mu+\eta}}
V_\lambda\otimes V_{\mu+\eta}\xrightarrow{\iota\otimes T_{\mu,\eta}}
V_\lambda\otimes V_\mu\otimes V_\eta
$$
coincide. To see this we just have to check how these morphisms act
on the highest weight vector and then observe that $\Phi$ acts
trivially on $\xi_\lambda\otimes\xi_\mu\otimes\xi_\eta$, since both
$t_{12}$ and $t_{23}$ preserve the one-dimensional space spanned by
this vector and in particular commute on this space. Next observe
that the composition in the formulation of the lemma coincides with
$g^\hbar_{\mu,\eta}S_{\mu+\eta}$ by definition. It turns out that
these two properties are enough to establish the cocycle property
$g^\hbar_{\lambda+\mu,\eta}g^\hbar_{\lambda,\mu}
=g^\hbar_{\lambda,\mu+\eta}g^\hbar_{\mu,\eta}$. To show this we can
and shall strictify the category $\D(\g,\hbar)$ and thus omit $\Phi$
in all computations. For example the equality of the above two
compositions now reads as
\begin{equation} \label{eTT}
(T_{\lambda,\mu}\otimes\iota)T_{\lambda+\mu,\eta}
=(\iota\otimes T_{\mu,\eta})T_{\lambda,\mu+\eta}.
\end{equation}
Then the morphisms
$$
V_{\bar\lambda+\bar\mu+\bar\eta}\otimes V_{\lambda+\mu+\eta}\to \C
$$
given by
$$
S_\lambda(\iota\otimes S_\mu\otimes\iota) (\iota\otimes\iota\otimes
S_\eta\otimes\iota\otimes\iota) (\iota\otimes
T_{\bar\mu,\bar\eta}\otimes
T_{\eta,\mu}\otimes\iota)(T_{\bar\lambda,\bar\mu+\bar\eta}\otimes
T_{\mu+\eta,\lambda})
$$
and
$$
S_\lambda(\iota\otimes S_\mu\otimes\iota)
(T_{\bar\lambda,\bar\mu}\otimes T_{\mu,\lambda}) (\iota\otimes
S_\eta\otimes\iota) (T_{\bar\lambda+\bar\mu,\bar\eta}\otimes
T_{\eta,\lambda+\mu})
$$
coincide. On the other hand, the first  morphism is equal to
$$
S_\lambda(\iota\otimes g^\hbar_{\mu,\eta}S_{\mu+\eta}\otimes\iota)
(T_{\bar\lambda,\bar\mu+\bar\eta}\otimes T_{\mu+\eta,\lambda})
=g^\hbar_{\mu,\eta}g^\hbar_{\lambda,\mu+\eta}S_{\lambda+\mu+\eta},
$$
whereas the second morphism equals
$$
g^\hbar_{\lambda,\mu}S_{\lambda+\mu} (\iota\otimes
S_\eta\otimes\iota) (T_{\bar\lambda+\bar\mu,\bar\eta}\otimes
T_{\eta,\lambda+\mu})
=g^\hbar_{\lambda,\mu}g^\hbar_{\lambda+\mu,\eta}S_{\lambda+\mu+\eta}.
$$
Since $S_{\lambda+\mu+\eta}\ne0$ we get
$g^\hbar_{\mu,\eta}g^\hbar_{\lambda,\mu+\eta}
=g^\hbar_{\lambda,\mu}g^\hbar_{\lambda+\mu,\eta}$.

\smallskip

It remains to check that the cocycle is symmetric. First observe
that $T_{\mu,\eta}$ coincides with the composition
$$
V_{\mu+\eta}\xrightarrow{T_{\eta,\mu}}V_\eta\otimes V_\mu
\xrightarrow{q^{-(\mu,\eta)}\Sigma q^t}V_\mu\otimes V_\eta,
$$
where $q=e^{\pi i\hbar}$. To see this we again look at the action on
the highest weight vector. Then the claim follows from
\begin{equation} \label{ethigh}
t(\xi_\eta\otimes\xi_\mu)=(\mu,\eta)\xi_\eta\otimes\xi_\mu, 
\end{equation}
which is a consequence of \eqref{eCas1} and the fact that $C$ acts on $V_\lambda$
as multiplication by $(\lambda,\lambda+2\rho)$.
We now
strictify $\D(\g,\hbar)$ and do all computations omitting~$\Phi$.
Denote by $\sigma$ the braiding in our new strict category.  By
definition we have
$$
S_\mu(\iota\otimes S_\eta\otimes\iota)(T_{\bar\mu,\bar\eta}\otimes
T_{\eta,\mu})=g^\hbar_{\mu,\eta}S_{\mu+\eta}.
$$
As
\begin{equation} \label{eTC}
T_{\mu,\eta}=q^{(\mu,\eta)}\sigma^{-1}T_{\eta,\mu},
\end{equation}
we can rewrite this as
$$
g^\hbar_{\mu,\eta}S_{\mu+\eta}
=S_\mu(\iota\otimes S_\eta\otimes\iota)(\sigma^{-1}\otimes\sigma)
(T_{\bar\eta,\bar\mu}\otimes T_{\mu,\eta}).
$$
By the hexagon identities
$\sigma_{12,3}=(\sigma\otimes\iota)(\iota\otimes\sigma)$ and
$\sigma^{-1}_{1,23}=(\sigma^{-1}\otimes\iota)(\iota\otimes\sigma^{-1})$
we have
$$
\sigma^{-1}\otimes\sigma
=(\sigma^{-1}_{1,23}\otimes\iota)(\iota\otimes\sigma_{12,3}).
$$
Therefore by naturality of $\sigma$ we get
\begin{align*}
g^\hbar_{\mu,\eta}S_{\mu+\eta}
&=S_\mu(\iota\otimes S_\eta\otimes\iota)
(\sigma^{-1}_{1,23}\otimes\iota)(\iota\otimes\sigma_{12,3})
(T_{\bar\eta,\bar\mu}\otimes T_{\mu,\eta})\\
&=S_\mu(S_\eta\otimes\iota\otimes\iota)
(\iota\otimes\sigma_{12,3})
(T_{\bar\eta,\bar\mu}\otimes T_{\mu,\eta})\\
&=S_\eta(\iota\otimes\iota\otimes S_\mu)
(\iota\otimes\sigma_{12,3})
(T_{\bar\eta,\bar\mu}\otimes T_{\mu,\eta})\\
&=S_\eta(\iota\otimes S_\mu\otimes\iota)
(T_{\bar\eta,\bar\mu}\otimes T_{\mu,\eta})\\
&=g^\hbar_{\eta,\mu}S_{\mu+\eta}.
\end{align*}
Hence $g^\hbar_{\mu,\eta}=g^\hbar_{\eta,\mu}$.
\ep

It is well-known that a symmetric cocycle must be a coboundary. We
formulate this in the following a bit more precise form.

\begin{lemma}
Let $(\mu,\eta)\mapsto c_{\mu,\eta}$ be a $\C^*$-valued symmetric
normalized $2$-cocycle on $P_+$. Then for any nonzero complex
numbers $b_1,\dots,b_r$ there exists a unique map $P_+\ni\mu\mapsto
b_\mu\in\C^*$ such that
$$
c_{\mu,\eta}=b_{\mu+\eta}b^{-1}_\mu b^{-1}_\eta,\ \
b_0=1,\ \ b_{\omega_i}=b_i\ \ \hbox{for} \ \ i=1,\ldots,r.
$$
\end{lemma}

Here $\omega_1,\dots,\omega_r$ are the fundamental weights.

\bp It is clear that the map $b$ is unique if it exists. To
show existence, for a weight $\mu\in P_+$,
$\mu=k_1\omega_1+\ldots+k_n\omega_r$, put $|\mu|=k_1+\ldots+k_r$.
Define $b_\mu$ by induction on $|\mu|$ as follows. If $\mu-\omega_i$
is dominant for some $i$ then put
$b_\mu=c_{\mu-\omega_i,\omega_i}b_{\mu-\omega_i}b_{\omega_i}$. We
have to check that $b_\mu$ is well-defined. In other words, if
$\mu=\nu+\omega_i+\omega_j$ then we must show that
$$
c_{\nu+\omega_j,\omega_i}b_{\nu+\omega_j}b_{\omega_i}
=c_{\nu+\omega_i,\omega_j}b_{\nu+\omega_i}b_{\omega_j}.
$$
Using the cocycle identities
$$
c_{\nu+\omega_j,\omega_i}c_{\nu,\omega_j}
=c_{\nu,\omega_i+\omega_j}c_{\omega_i,\omega_j}
\ \ \hbox{and}\ \
c_{\nu+\omega_i,\omega_j}c_{\nu,\omega_i}
=c_{\nu,\omega_j+\omega_i}c_{\omega_j,\omega_i}
$$
and that $c_{\omega_i,\omega_j}=c_{\omega_j,\omega_i}$, we
equivalently have to check that
$$
c_{\nu,\omega_i}b_{\nu+\omega_j}b_{\omega_i}
=c_{\nu,\omega_j}b_{\nu+\omega_i}b_{\omega_j}.
$$
Since $c_{\nu,\omega_i}=b_{\nu+\omega_i}b_\nu^{-1}b_{\omega_i}^{-1}$
and $c_{\nu,\omega_j}=b_{\nu+\omega_j}b_\nu^{-1}b_{\omega_j}^{-1}$
by the inductive assumption, the identity indeed holds.

Therefore we have constructed a map $b$ such that $b_0=1$,
$b_{\omega_i}=b_i$ and
$c_{\mu,\omega_i}=b_{\mu+\omega_i}b_\mu^{-1}b_{\omega_i}^{-1}$ for
$i=1,\ldots,r$ and $\mu\in P_+$. By induction on $|\eta|$ one can
easily check that the identity $c_{\mu,\eta}=b_{\mu+\eta}b^{-1}_\mu
b^{-1}_\eta$ holds for all $\mu,\eta\in P_+$. \ep

For generic $\hbar$ fix a map $P_+\ni\mu\mapsto
g^\hbar_\mu\in\C^*$ such that
$$
g^\hbar_{\mu,\eta}g^\hbar_\mu g^\hbar_\eta=g^\hbar_{\mu+\eta}.
$$
In Section~\ref{sequivalence} we shall require an additional
property of this map, which determines the cochain~$g^\hbar_\mu$ up
to a character of the quotient $P/Q$ of the weight lattice by the
root lattice, but in this section as well as in the next one any
$g^\hbar_\mu$ will do.

Define $S^\hbar_\mu=g^\hbar_\mu S_\mu \colon V_{\bar\mu}\otimes
V_\mu\to\C$. We modify \eqref{etr0} by introducing the maps
\begin{equation} \label{etr1}
\tr^{\eta,\hbar}_{\mu,\lambda+\mu}\colon V_{\bar\mu+\bar\eta}\otimes
V_{\lambda+\mu+\eta} \xrightarrow{T_{\bar\mu,\bar\eta}\otimes
T_{\eta,\lambda+\mu}} V_{\bar\mu}\otimes V_{\bar\eta}\otimes
V_\eta\otimes V_{\lambda+\mu} \xrightarrow{(\iota\otimes
S^\hbar_\eta\otimes\iota)B}V_{\bar\mu}\otimes V_{\lambda+\mu},
\end{equation}
where $B=(\Phi\otimes\iota)\Phi^{-1}_{12,3,4}$.

\begin{lemma}
The morphisms \eqref{etr1} are coherent, that is, the composition
$$
V_{\bar\mu+\bar\eta+\bar\nu}\otimes V_{\lambda+\mu+\eta+\nu}
\xrightarrow{\tr^{\nu,\hbar}_{\mu+\eta,\lambda+\mu+\eta}}
V_{\bar\mu+\bar\eta}\otimes V_{\lambda+\mu+\eta}
\xrightarrow{\tr^{\eta,\hbar}_{\mu,\lambda+\mu}} V_{\bar\mu}\otimes
V_{\lambda+\mu}
$$
coincides with $\tr^{\eta+\nu,\hbar}_{\mu,\lambda+\mu}$.
\end{lemma}

\bp We strictify $\D(\g,\hbar)$ once again. Then the composition
$\tr^{\eta,\hbar}_{\mu,\lambda+\mu}
\tr^{\nu,\hbar}_{\mu+\eta,\lambda+\mu+\eta}$ equals
\begin{multline*}
(\iota\otimes S^\hbar_\eta\otimes\iota)
(T_{\bar\mu,\bar\eta}\otimes T_{\eta,\lambda+\mu})
(\iota\otimes S^\hbar_\nu\otimes\iota)
(T_{\bar\mu+\bar\eta,\bar\nu}\otimes T_{\nu,\lambda+\mu+\eta})\\
=(\iota\otimes S^\hbar_\eta\otimes\iota) (\iota\otimes\iota\otimes
S^\hbar_\nu\otimes\iota\otimes\iota)
(T_{\bar\mu,\bar\eta}\otimes\iota\otimes\iota\otimes
T_{\eta,\lambda+\mu}) (T_{\bar\mu+\bar\eta,\bar\nu}\otimes
T_{\nu,\lambda+\mu+\eta}).
\end{multline*}
Using \eqref{eTT} this can be written
$$
(\iota\otimes S^\hbar_\eta\otimes\iota)
(\iota\otimes\iota\otimes S^\hbar_\nu\otimes\iota\otimes\iota)
(\iota\otimes T_{\bar\eta,\bar\nu}\otimes T_{\nu,\eta}\otimes\iota)
(T_{\bar\mu,\bar\eta+\bar\nu}\otimes T_{\eta+\nu,\lambda+\mu}).
$$
By definition of $g^\hbar_{\eta,\nu}$ and using that
$g^\hbar_{\eta,\nu}g^\hbar_\eta g^\hbar_\nu=g^\hbar_{\eta+\nu}$ we
have
$$
S^\hbar_\eta(\iota\otimes S^\hbar_\nu\otimes\iota)
(T_{\bar\eta,\bar\nu}\otimes T_{\nu,\eta})
=S^\hbar_{\eta+\nu}.
$$
Therefore the above expression equals
$$
(\iota\otimes S^\hbar_{\eta+\nu}\otimes\iota)
(T_{\bar\mu,\bar\eta+\bar\nu}\otimes T_{\eta+\nu,\lambda+\mu})
=\tr^{\eta+\nu,\hbar}_{\mu,\lambda+\mu}.
$$
\ep

Using the morphisms $\tr^{\eta,\hbar}_{\mu,\lambda+\mu}$ we can
therefore define a $\g$-module
$$
M^\hbar_\lambda=\lim_{\xleftarrow[\mu]{}}V_{\bar\mu}\otimes
V_{\lambda+\mu}.
$$
Again we consider $M^\hbar_\lambda$ as a topological $\g$-module
with a base of neighbourhoods of zero given by the kernels of the
maps $M^\hbar_\lambda\to V_{\bar\mu}\otimes V_{\lambda+\mu}$, while
any module in $\D(\g,\hbar)$ is considered with discrete topology.

\begin{proposition} \label{represento}
For $\lambda\in P$ and generic $\hbar\in\C$ the topological module
$M^\hbar_\lambda$ is isomorphic to $M_\lambda$. In particular, for
any such $\hbar$ the functor $\D(\g,\hbar)\to\Vect$, $V\mapsto
V(\lambda)$, is naturally isomorphic to
$\Hom_\g(M^\hbar_\lambda,\cdot)$.
\end{proposition}

\bp Fix a regular dominant integral weight $\mu$ (that is,
$\mu$ lies in the interior of the Weyl chamber). Then $n\mu$
dominates any other weight for sufficiently large $n$. Choose
$n_0\in\N$ such that $n_0\mu+\lambda\ge0$. Then $M_\lambda$ is
isomorphic to the inverse limit of
$$
V_{n_0\bar\mu}\otimes V_{\lambda+n_0\mu}
\xleftarrow{\tr^{\mu}_{n_0\mu,\lambda+n_0\mu}}
V_{(n_0+1)\bar\mu}\otimes V_{\lambda+(n_0+1)\mu}
\xleftarrow{\tr^{\mu}_{(n_0+1)\mu,\lambda+(n_0+1)\mu}}
V_{(n_0+2)\bar\mu}\otimes V_{\lambda+(n_0+2)\mu}\longleftarrow\dots,
$$
and $M^\hbar_\lambda$ is the inverse limit of
$$
V_{n_0\bar\mu}\otimes V_{\lambda+n_0\mu}
\xleftarrow{\tr^{\mu,\hbar}_{n_0\mu,\lambda+n_0\mu}}
V_{(n_0+1)\bar\mu}\otimes V_{\lambda+(n_0+1)\mu}
\xleftarrow{\tr^{\mu,\hbar}_{(n_0+1)\mu,\lambda+(n_0+1)\mu}}
V_{(n_0+2)\bar\mu}\otimes V_{\lambda+(n_0+2)\mu}\longleftarrow\dots.
$$
It is therefore enough to find isomorphisms $f_n$ of
$V_{n\bar\mu}\otimes V_{\lambda+n\mu}$ onto itself such that for all
$n\ge n_0$ we have
$$
f_n\tr^{\mu}_{n\mu,\lambda+n\mu}=
\tr^{\mu,\hbar}_{n\mu,\lambda+n\mu}f_{n+1}.
$$
We construct $f_n$ by induction on $n$. Take $f_{n_0}$ to be the
identity map. Assuming that $f_n$ is constructed, observe that
$\tr^{\mu}_{n\mu,\lambda+n\mu}$ is surjective since it maps the
vector $\zeta_{(n+1)\bar\mu}\otimes\xi_{\lambda+(n+1)\mu}$ onto the
cyclic vector $\zeta_{n\bar\mu}\otimes\xi_{\lambda+n\mu}$. It
follows that for generic $\hbar$ the map
$\tr^{\mu,\hbar}_{n\mu,\lambda+n\mu}$ is surjective as well.
Therefore both maps $f_n\tr^{\mu}_{n\mu,\lambda+n\mu}$ and
$\tr^{\mu,\hbar}_{n\mu,\lambda+n\mu}$ are surjective. This is enough
to conclude that $f_{n+1}$ exists. Indeed, the claim is that if $g_1$
and $g_2$ are surjective morphisms $V\to W$ of finite dimensional
$\g$-modules then there exists an isomorphism $f$ of $V$ onto itself
such that $g_1f=g_2$. To see this we can reduce to the situation
when $V=U\otimes\C^n$ and $W=U\otimes\C^m$ for some irreducible
$\g$-module $U$. Then $g_i=\iota\otimes h_i$, where
$h_i\colon\C^n\to \C^m$ is a linear surjective map. Clearly we can
find an invertible linear map $h\colon\C^n\to \C^n$ such that
$h_1h=h_2$, and then put $f=\iota\otimes h$. \ep

\bigskip

\section{A comonoid structure on the representing object} \label{s5}

In the previous section we showed that for generic $\hbar\in\C$ the
topological $\g$-module
$$
M^\hbar=\oplus_{\lambda\in P}M^\hbar_\lambda
$$
represents the forgetful functor $\D(\g,\hbar)\to\Vect$. In this
section we shall turn the functor $\Hom_\g(M^\hbar,\cdot)$ into a
tensor functor. To do this we introduce a comonoid structure on
$M^\hbar$.

Define
$$
M^\hbar_{\lambda_1}\hat\otimes M^\hbar_{\lambda_2}
=\lim_{\xleftarrow[\mu_1,\mu_2]{}}(V_{\bar\mu_1}\otimes
V_{\lambda_1+\mu_1})\otimes (V_{\bar\mu_2}\otimes V_{\lambda_2+\mu_2})
$$
and then
$$
M^\hbar\hat\otimes M^\hbar=\prod_{\lambda_1,\lambda_2\in P}
M^\hbar_{\lambda_1}\hat\otimes M^\hbar_{\lambda_2}.
$$
Higher tensor powers of $M^\hbar$ are defined similarly. We want to
define
$$
\delta^\hbar\colon M^\hbar \to M^\hbar\hat\otimes M^\hbar.
$$
The restriction of $\delta^\hbar$ to $M^\hbar_\lambda$ composed with
the projection $M^\hbar\hat\otimes M^\hbar\to
M^\hbar_{\lambda_1}\hat\otimes M^\hbar_{\lambda_2}$ will be nonzero
only if $\lambda=\lambda_1+\lambda_2$, so $\delta^\hbar$ is
determined by maps
$$
\delta^\hbar_{\lambda_1,\lambda_2}\colon
M^\hbar_{\lambda_1+\lambda_2} \to M^\hbar_{\lambda_1}\hat\otimes
M^\hbar_{\lambda_2}.
$$
Motivated by \eqref{eEKco} we define these morphisms using the
compositions
\begin{equation} \label{em}
\begin{split}
m^\hbar_{\mu,\eta,\lambda_1,\lambda_2}\colon
V_{\bar\mu+\bar\eta}\otimes V_{\lambda_1+\lambda_2+\mu+\eta}
&\xrightarrow{T_{\bar\mu,\bar\eta}\otimes
T_{\lambda_1+\mu,\lambda_2+\eta}}(V_{\bar\mu}\otimes
V_{\bar\eta})\otimes (V_{\lambda_1+\mu}
\otimes V_{\lambda_2+\eta})\\
&\xrightarrow{q^{(\lambda_1+\mu,\eta)}B^{-1}(\iota\otimes \Sigma
q^t\otimes\iota)B}(V_{\bar\mu}\otimes  V_{\lambda_1+\mu})\otimes
(V_{\bar\eta}\otimes V_{\lambda_2+\eta}),
\end{split}
\end{equation}
where $q=e^{\pi i\hbar}$ and $B=(\Phi\otimes\iota)\Phi^{-1}_{12,3,4}$.

\begin{lemma} \label{mtr}
The morphisms $m^\hbar$ are consistent with the morphisms
$\tr^{\cdot,\hbar}$ defining the inverse limits, so they define
morphisms $\delta^\hbar_{\lambda_1,\lambda_2}\colon
M^\hbar_{\lambda_1+\lambda_2} \to M^\hbar_{\lambda_1}\hat\otimes
M^\hbar_{\lambda_2}$.
\end{lemma}

\bp We have to check that
\begin{equation} \label{etrm}
(\tr^{\nu,\hbar}_{\mu,\lambda_1+\mu}\otimes
\tr^{\omega,\hbar}_{\eta,\lambda_2+\eta})
m^\hbar_{\mu+\nu,\eta+\omega,\lambda_1,\lambda_2}
=m^\hbar_{\mu,\eta,\lambda_1,\lambda_2}
\tr^{\nu+\omega,\hbar}_{\mu+\eta,\lambda_1+\lambda_2+\mu+\eta}.
\end{equation}
Since
$$
\tr^{\nu,\hbar}_{\mu,\lambda_1+\mu}\otimes
\tr^{\omega,\hbar}_{\eta,\lambda_2+\eta}=
(\tr^{\nu,\hbar}_{\mu,\lambda_1+\mu}\otimes\iota\otimes\iota)
(\iota\otimes\iota\otimes \tr^{\omega,\hbar}_{\eta,\lambda_2+\eta}),
$$
it suffices to check this assuming that either $\nu$ or $\omega$ is
zero. We shall only consider the case $\omega=0$. We therefore have
to check that
\begin{equation} \label{emtr}
(\tr^{\nu,\hbar}_{\mu,\lambda_1+\mu}\otimes\iota\otimes\iota)
m^\hbar_{\mu+\nu,\eta,\lambda_1,\lambda_2}
=m^\hbar_{\mu,\eta,\lambda_1,\lambda_2}
\tr^{\nu,\hbar}_{\mu+\eta,\lambda_1+\lambda_2+\mu+\eta}.
\end{equation}
We strictify $\D(\g,\hbar)$. Denote by $\sigma$ the braiding in the
strict tensor category. In the computation below we omit subindices
of the morphisms $T$ since they are completely determined by the
target modules. We will keep track of some of them to get the right
power of $q$.  Thus by definition of $\tr^{\cdot,\hbar}$ and
$m^\hbar$ the left hand side of \eqref{emtr} is equal to
\begin{multline*}
q^{(\lambda_1+\mu+\nu,\eta)}(\iota\otimes
S^\hbar_\nu\otimes\iota\otimes\iota\otimes\iota)
(T_{\bar\mu,\bar\nu}\otimes T\otimes\iota\otimes\iota)
(\iota\otimes\sigma\otimes\iota)
(T_{\bar\mu+\bar\nu,\bar\eta}\otimes T)\\
=q^{(\lambda_1+\mu+\nu,\eta)}(\iota\otimes
S^\hbar_\nu\otimes\iota\otimes\iota\otimes\iota)
(\iota\otimes\iota\otimes \sigma_{1,23}\otimes\iota)
(T_{\bar\mu,\bar\nu}\otimes\iota\otimes
T\otimes\iota)(T_{\bar\mu+\bar\nu,\bar\eta}\otimes T).
\end{multline*}
Using the identity $(T_{\bar\mu,\bar\nu}\otimes\iota)
T_{\bar\mu+\bar\nu,\bar\eta}=(\iota\otimes
T_{\bar\nu,\bar\eta})T_{\bar\mu,\bar\nu+\bar\eta}$, see \eqref{eTT},
the above expression can be written as
$$
q^{(\lambda_1+\mu+\nu,\eta)}(\iota\otimes
S^\hbar_\nu\otimes\iota\otimes\iota\otimes\iota)
(\iota\otimes\iota\otimes\sigma_{1,23}\otimes\iota) (\iota\otimes
T_{\bar\nu,\bar\eta}\otimes \iota\otimes
T)(T_{\bar\mu,\bar\nu+\bar\eta}\otimes T).
$$
By \eqref{eTC} we have $T_{\bar\nu,\bar\eta}
=q^{-(\bar\nu,\bar\eta)}\sigma
T_{\bar\eta,\bar\nu}=q^{-(\nu,\eta)}\sigma T_{\bar\eta,\bar\nu}$, so
we get
$$
q^{(\lambda_1+\mu,\eta)}(\iota\otimes
S^\hbar_\nu\otimes\iota\otimes\iota\otimes\iota)
(\iota\otimes\iota\otimes\sigma_{1,23}\otimes\iota) (\iota\otimes
\sigma T\otimes \iota\otimes T)(T\otimes T).
$$
By the hexagon identity we have $(\iota\otimes\sigma_{1,23})
(\sigma\otimes\iota\otimes\iota)=\sigma_{1,234}$, so the above
expression equals
$$
q^{(\lambda_1+\mu,\eta)}(\iota\otimes
S^\hbar_\nu\otimes\iota\otimes\iota\otimes\iota)
(\iota\otimes\sigma_{1,234}\otimes\iota) (\iota\otimes T\otimes
\iota\otimes T)(T\otimes T).
$$
Using again that $(T\otimes\iota)T=(\iota\otimes T)T$, we get
\begin{equation*}
\begin{split}
q^{(\lambda_1+\mu,\eta)}(\iota\otimes
S^\hbar_\nu&\otimes\iota\otimes\iota\otimes\iota)
(\iota\otimes\sigma_{1,234}\otimes\iota) (T\otimes
\iota\otimes\iota\otimes T)(T\otimes T)\\
&=q^{(\lambda_1+\mu,\eta)}
(\iota\otimes\sigma\otimes\iota)
(\iota\otimes\iota\otimes S^\hbar_\nu\otimes\iota\otimes\iota)
(T\otimes \iota\otimes\iota\otimes T)(T\otimes T)\\
&=q^{(\lambda_1+\mu,\eta)}
(\iota\otimes\sigma\otimes\iota)
(T\otimes T)(\iota\otimes S^\hbar_\nu\otimes\iota)(T\otimes T),
\end{split}
\end{equation*}
which is exactly the right hand side of \eqref{emtr}.
\ep

Using the morphisms $\delta^\hbar_{\lambda_1,\lambda_2}$ we can in
an obvious way define morphisms
$$
(\delta^\hbar\otimes\iota)\delta^\hbar,
(\iota\otimes\delta^\hbar)\delta^\hbar \colon M^\hbar\to
M^\hbar\hat\otimes M^\hbar \hat\otimes M^\hbar.
$$

\begin{lemma}
We have $\Phi(\delta^\hbar\otimes\iota)\delta^\hbar
=(\iota\otimes\delta^\hbar)\delta^\hbar$.
\end{lemma}

\bp For $\lambda_1,\lambda_2,\lambda_3\in P$ we have to check that
$$
\Phi(\delta^\hbar_{\lambda_1,\lambda_2}\otimes\iota)
\delta^\hbar_{\lambda_1+\lambda_2,\lambda_3}=
(\iota\otimes\delta^\hbar_{\lambda_2,\lambda_3})
\delta^\hbar_{\lambda_1,\lambda_2+\lambda_3}.
$$
This reduces to showing that
$$
\Phi_{12,34,56}
(m^\hbar_{\mu_1,\mu_2,\lambda_1,\lambda_2}\otimes\iota
\otimes\iota)
m^\hbar_{\mu_1+\mu_2,\mu_3,\lambda_1+\lambda_2,\lambda_3}
=(\iota\otimes\iota\otimes m^\hbar_{\mu_2,\mu_3,\lambda_2,\lambda_3})
m^\hbar_{\mu_1,\mu_2+\mu_3,\lambda_1,\lambda_2+\lambda_3}.
$$
Let us first check that the powers of $q$ in the definition of
$m^\hbar$ match. On the left hand side we have
$q^{(\lambda_1+\mu_1,\mu_2)
+(\lambda_1+\lambda_2+\mu_1+\mu_2,\mu_3)}$, whereas on the right
hand side we get $q^{(\lambda_2+\mu_2,\mu_3)
+(\lambda_1+\mu_1,\mu_2+\mu_3)}$, which obviously coincide.
Strictifying and omitting subindices in $T$, as we did in the proof
of the previous lemma, it remains to show that
$$
(\iota\otimes\sigma\otimes\iota\otimes\iota\otimes\iota)
(T\otimes T\otimes\iota\otimes\iota)(\iota\otimes\sigma\otimes\iota)
(T\otimes T)
=(\iota\otimes\iota\otimes\iota\otimes\sigma\otimes\iota)
(\iota\otimes\iota\otimes T\otimes T)(\iota\otimes\sigma\otimes\iota)
(T\otimes T).
$$
By naturality of the braiding, the left hand side equals
$$
(\iota\otimes\sigma\otimes\iota\otimes\iota\otimes\iota)
(\iota\otimes\iota\otimes\sigma_{1,23}\otimes\iota)
(T\otimes\iota\otimes T\otimes\iota)
(T\otimes T),
$$
whereas the right hand side equals
$$
(\iota\otimes\iota\otimes\iota\otimes\sigma\otimes\iota)
(\iota\otimes\sigma_{12,3}\otimes\iota\otimes\iota)
(\iota\otimes T\otimes\iota\otimes T)
(T\otimes T).
$$
As $(T\otimes\iota)T=(\iota\otimes T)T$, we thus only need to check
that
$$
(\sigma\otimes\iota\otimes\iota)
(\iota\otimes\sigma_{1,23})
=(\iota\otimes\iota\otimes\sigma)
(\sigma_{12,3}\otimes\iota),
$$
which is immediate from the hexagon identities
$\sigma_{1,23}=(\iota\otimes\sigma)(\sigma\otimes\iota)$
and $\sigma_{12,3}=(\sigma\otimes\iota)(\iota\otimes\sigma)$.
\ep

We next introduce a morphism $\eps^\hbar\colon M^\hbar\to\C$ by requiring it
to be nonzero only on $M^\hbar_0$,
where we set it to be the canonical morphism $M^\hbar_0\to
V_{\bar0}\otimes V_0=\C$, so that $\eps^\hbar\colon
M^\hbar_0\to\C$ is determined by the morphisms
$$
\tr^{\mu,\hbar}_{0,0}=S^\hbar_0\colon V_{\bar\mu}\otimes V_\mu\to\C.
$$

\begin{lemma}
We have $(\eps^\hbar\otimes\iota)\delta^\hbar=\iota
=(\iota\otimes\eps^\hbar)\delta^\hbar$.
\end{lemma}

\bp We have to check that on $M^\hbar_\lambda$ we have
$(\eps^\hbar\otimes\iota)\delta^\hbar_{0,\lambda}=\iota
=(\iota\otimes\eps^\hbar)\delta^\hbar_{\lambda,0}$. This 
follows from the fact that $m^\hbar_{0,\eta,0,\lambda}$
and $m^\hbar_{\mu,0,\lambda,0}$ are the identity maps. \ep

Therefore $M^\hbar$ is a comonoid, so $\Hom_\g(M^\hbar,\cdot)$
becomes a weak tensor functor $\D(\g,\hbar)\to \Vect$.

\begin{proposition} \label{ptensor}
For generic $\hbar\in\C$ the weak tensor functor
$\Hom_\g(M^\hbar,\cdot)\colon\D(\g,\hbar)\to\Vect$ is a tensor
functor.
\end{proposition}

\bp Let $V$ and $W$ be finite dimensional $\g$-modules.
We have to show that for generic $\hbar$ the map
$$
\Hom_\g(M^\hbar,V)\otimes\Hom_\g(M^\hbar,W)
\to \Hom_\g(M^\hbar,V\otimes W), \ \
f\otimes g\mapsto (f\otimes g)\delta^\hbar,
$$
is a linear isomorphism. As $\Hom_\g(M^\hbar,V) =\oplus_{\lambda\in
P}\Hom_\g(M^\hbar_\lambda,V)$ (notice that the direct sum is finite,
because $\Hom_\g(M^\hbar_\lambda,V)\ne0$ only if $V(\lambda)\ne0$),
we equivalently have to check that for any $\lambda\in P$ the map
$$
\bigoplus_{\lambda_1+\lambda_2=\lambda}
\Hom_\g(M^\hbar_{\lambda_1},V)\otimes\Hom_\g(M^\hbar_{\lambda_2},W)
\to \Hom_\g(M^\hbar_\lambda,V\otimes W),  \ \ f_{\lambda_1}\otimes
g_{\lambda_2} \mapsto (f_{\lambda_1}\otimes
g_{\lambda_2})\delta^\hbar_{\lambda_1,\lambda_2},
$$
is an isomorphism for generic $\hbar$. As
$\Hom_\g(M^\hbar_\lambda,V)$ is the inductive limit of
$\Hom_\g(V_{\bar\mu}\otimes V_{\lambda+\mu},V)$, it suffices to
check that for fixed $\lambda\in P$ and all sufficiently large
$\mu_1$ and $\mu_2$ the map
$$
\bigoplus_{\lambda_1+\lambda_2=\lambda} \Hom_\g(V_{\bar\mu_1}\otimes
V_{\lambda_1+\mu_1},V) \otimes\Hom_\g(V_{\bar\mu_2}\otimes
V_{\lambda_2+\mu_2},W) \to \Hom_\g(V_{\bar\mu_1+\bar\mu_2}\otimes
V_{\lambda_1+\lambda_2+\mu_1+\mu_2},V\otimes W),
$$
which maps $f_{\lambda_1}\otimes g_{\lambda_2}$ onto
$(f_{\lambda_1}\otimes g_{\lambda_2})
m^\hbar_{\mu_1,\mu_2,\lambda_1,\lambda_2}$, is an isomorphism for
generic $\hbar$. As the map is analytic in $\hbar$ outside a
discrete set, it suffices to check that it is an isomorphism for
$\hbar=0$. For sufficiently large $\mu_1$ we have isomorphisms
$\Hom_\g(V_{\bar\mu_1}\otimes V_{\lambda_1+\mu_1},V)\to
V(\lambda_1)$, $f\mapsto
f(\zeta_{\bar\mu_1}\otimes\xi_{\lambda_1+\mu_1})$, and similar
isomorphisms for $W$ and $V\otimes W$. It is then easy to verify
that under these isomorphisms the above map (for $\hbar=0$) becomes
$$
\bigoplus_{\lambda_1+\lambda_2=\lambda}
V(\lambda_1)\otimes W(\lambda_2)\to (V\otimes W)(\lambda),
\ \ v\otimes w\mapsto v\otimes w,
$$
which is clearly an isomorphism.
\ep

Recall that the construction of the comonoid $M^\hbar$ depends on the choice of a $1$-cochain $g^\hbar_\mu$ with coboundary
$g^\hbar_{\mu,\eta}$.

\begin{lemma} \label{lunicom}
Up to an isomorphism the comonoid $M^\hbar$ does not depend on the
choice of $g^\hbar_\mu$.
\end{lemma}

\bp Assume that $\tilde g^\hbar_{\mu}$ is another cochain. Denote by
$\tilde M^\hbar$ the new comonoid. The map $\chi\colon
P_+\to\C^*$ defined by $\chi(\mu)=\tilde g^\hbar_\mu
(g^\hbar_\mu)^{-1}$ is a homomorphism. Then it is straightforward
to check that the morphisms $V_{\bar\mu}\otimes V_{\lambda+\mu}\to
V_{\bar\mu}\otimes V_{\lambda+\mu}$ given by multiplication with
$\chi(\mu)$ induce an isomorphism of  $\tilde M^\hbar$ and $M^\hbar$
which respects their comonoid structures. \ep

So far we have constructed for generic $\hbar\in\C$ a tensor functor
$\D(\g,\hbar)\to\Vect$. Up to natural isomorphisms of tensor
functors the construction is canonical. Furthermore, disregarding
the tensor structure the functor is naturally isomorphic to the
forgetful functor. By the discussion after Proposition~\ref{QAvsME}
(or by combining Proposition~\ref{reconstruction} and
Proposition~\ref{QAvsME}(ii)) it already follows that for generic
$\hbar$ a twisting of $(\widehat{\C[G]}, \Dhat,\hat\eps,\Phi(\hbar
t_{12},\hbar t_{23}),e^{\pi i\hbar t})$ is isomorphic to a discrete
bialgebra, or equivalently, there exists a twist
$\F^\hbar\in\U(G\times G)$ such that $\Phi(\hbar t_{12},\hbar
t_{23})_{\F^\hbar}=1$. In the next section we will show that this
bialgebra is isomorphic to $\widehat{\C[G_q]}$ by turning the tensor
functor $\D(\g,\hbar)\to\Vect$ into an equivalence of the braided
monoidal categories $\D(\g,\hbar)$ and~$\CC(\g,\hbar)$.

\medskip

In the remaining part of the section we will summarize how one gets
a twist $\F^\hbar$ such that $\Phi(\hbar t_{12},\hbar
t_{23})_{\F^\hbar}=1$ in the form of an ``algorithm".

\smallskip

1) For $\mu,\eta\in P_+$ compute the image $g^\hbar_{\mu,\eta}$
of
$\zeta_{\bar\mu+\bar\eta}\otimes \xi_{\mu+\eta}$ under the
composition
$$
V_{\bar\mu+\bar\eta}\otimes V_{\mu+\eta}
\xrightarrow{T_{\bar\mu,\bar\eta}\otimes T_{\eta,\mu}}
V_{\bar\mu}\otimes V_{\bar\eta}\otimes V_\eta\otimes V_{\mu}
\xrightarrow{(\iota\otimes S_\eta\otimes\iota)B} V_{\bar\mu}\otimes
V_{\mu}\xrightarrow{S_\mu}\C,
$$
where $B$ is $(\Phi\otimes \iota)\Phi^{-1}_{12,3,4}$. Fix nonzero
numbers $z_1,\dots,z_r$. For $\mu=\omega_{i_1}+\dots+\omega_{i_k}$
put
$$
g^\hbar_\mu=z_{i_1}\prod^{k-1}_{l=1}g^\hbar_{\omega_{i_1}
+\dots+\omega_{i_l},\omega_{i_{l+1}}}z_{i_{l+1}}.
$$

2) Fix a regular dominant integral weight $\mu$. For each $\lambda
\in P$ choose $n_\lambda\in\N$ such that $n_\lambda\mu+\lambda\ge0$.
Then inductively choose isomorphisms $f^\lambda_n$, $n\ge
n_\lambda$, of $V_{n\bar\mu}\otimes V_{\lambda +n\mu}$ onto itself
such that $f^\lambda_{n_\lambda}$ is the identity map and for each
$n\ge n_\lambda$ the following diagram commutes:
$$
\begin{CD}
V_{n\bar\mu}\otimes V_{\lambda+n\mu}
@<{\tr^{\mu}_{n\mu,\lambda+n\mu}}<<
V_{(n+1)\bar\mu}\otimes V_{\lambda+(n+1)\mu}\\
@V{f^\lambda_n}VV     @VV{f^\lambda_{n+1}}V\\
V_{n\bar\mu}\otimes V_{\lambda+n\mu}
@<{\tr^{\mu,\hbar}_{n\mu,\lambda+n\mu}}<<
V_{(n+1)\bar\mu}\otimes V_{\lambda+(n+1)\mu},
\end{CD}
$$
where $\tr^{\mu,\hbar}_{n\mu,\lambda+n\mu}$ is the composition
$$
V_{(n+1)\bar\mu}\otimes V_{\lambda+(n+1)\mu}
\xrightarrow{T_{n\bar\mu,\bar\mu}\otimes T_{\mu,\lambda+n\mu}}
V_{n\bar\mu}\otimes V_{\bar\mu}\otimes V_\mu\otimes V_{\lambda+n\mu}
\xrightarrow{(\iota\otimes g^\hbar_\mu
S_\mu\otimes\iota)B}V_{n\bar\mu}\otimes V_{\lambda+n\mu}
$$
with $B=(\Phi\otimes\iota)\Phi^{-1}_{12,3,4}$, and
$\tr^{\mu}_{n\mu,\lambda+n\mu}$ is defined similarly with
$g^\hbar_\mu$ and $\Phi$ trivial.

3) Let $\eta,\nu\in P_+$. Then $\F^\hbar$ is defined
by requiring that it acts on the
space $V_\eta\otimes V_\nu$
by the operator~$\F^\hbar_{\eta,\nu}$
such that for weights $\lambda_1$ and $\lambda_2$ with
$V_\eta(\lambda_1)\ne0$, $V_\nu(\lambda_2)\ne0$ and
$\lambda=\lambda_1+\lambda_2$ the following diagram
commutes:
$$
\xymatrix{
\Hom_\g(V_{n\bar\mu}\otimes V_{\lambda_1+n\mu},V_\eta)\otimes
\Hom_\g(V_{m\bar\mu}\otimes V_{\lambda_2+m\mu},V_\nu)
\ar[r]\ar[d] & \Hom_\g(V_{(n+m)\bar\mu}\otimes V_{\lambda+(n+m)\mu},
V_\eta\otimes V_\nu)\ar[d]\\
V_\eta(\lambda_1)\otimes V_\nu(\lambda_2)
\ar[r]^{(\F^\hbar_{\eta,\nu})^{-1}} &(V_\eta\otimes V_\nu)(\lambda)
}
$$
where the left vertical arrow is the map
$$
f\otimes g\mapsto ff_n^{\lambda_1}(\zeta_{n\bar\mu}\otimes
\xi_{\lambda_1+n\mu})\otimes
gf_m^{\lambda_2}(\zeta_{m\bar\mu}\otimes
\xi_{\lambda_2+m\mu}),
$$
the right vertical arrow is the map
$$
f\mapsto ff_{n+m}^{\lambda_1+\lambda_2}(\zeta_{(n+m)\bar\mu}\otimes
\xi_{\lambda+(n+m)\mu})
$$
and finally the top horizontal arrow maps $f\otimes g$ onto the
composition
\begin{align*}
V_{(n+m)\bar\mu}\otimes V_{\lambda+(n+m)\mu}
&\xrightarrow{T_{n\bar\mu,m\bar\mu}\otimes
T_{\lambda_1+n\mu,\lambda_2+m\mu}}(V_{n\bar\mu}\otimes
V_{m\bar\mu})\otimes (V_{\lambda_1+n\mu}
\otimes V_{\lambda_2+m\mu})\\
&\xrightarrow{q^{(\lambda_1+n\mu,m\mu)}B^{-1}(\iota\otimes \Sigma
q^t\otimes\iota)B}(V_{n\bar\mu}\otimes  V_{\lambda_1+n\mu})\otimes
(V_{m\bar\mu}\otimes V_{\lambda_2+m\mu})\\
&\xrightarrow{f\otimes g}V_\eta\otimes V_\nu
\end{align*}
with $q=e^{\pi i\hbar}$
and $B=(\Phi\otimes\iota)\Phi^{-1}_{12,3,4}$. Here $n$ and $m$ can
be any natural numbers large enough so that $n\ge n_{\lambda_1}$,
$m\ge n_{\lambda_2}$, $n+m\ge n_{\lambda_1+\lambda_2}$ and the
vertical arrows are isomorphisms.

\bigskip

\section{Representing $U_q\g$ by endomorphisms of the functor}
\label{sequivalence}

In this section we will show that $U_q\g$, $q=e^{\pi i\hbar}$, can
be represented by endomorphisms of the functor
$\Hom_\g(M^\hbar,\cdot)$, so $\Hom_\g(M^\hbar,\cdot)$ can be
considered as a functor $\D(\g,\hbar)\to\CC(\g,\hbar)$. For this it
is natural to try to define an action of the opposite algebra
$(U_q\g)^{op}$ on~$M^\hbar$. We will show a bit less, namely, that
there is an action of a larger algebra $U_q\tilde\g$ such that the
corresponding action on the functor factors through $U_q\g$.

Denote by $U_q\tilde\g$ the universal algebra generated by
elements $E_i,F_i,K_i,K_i^{-1}$, $1\le i\le r$, such that
$$
K_iK_i^{-1}=K_i^{-1}K_i=1,\ \ K_iK_j=K_jK_i,\ \
K_iE_jK_i^{-1}=q_i^{-a_{ij}}E_j,\ \
K_iF_jK_i^{-1}=q_i^{a_{ij}}F_j,
$$
$$
E_iF_j-F_jE_i=-\delta_{ij}\frac{K_i-K_i^{-1}}{q_i-q_i^{-1}}.
$$
This is a Hopf algebra with coproduct $\Dhat_q$ defined by
$$
\Dhat_q(K_i)=K_i\otimes K_i,\ \
\Dhat_q(E_i)=E_i\otimes1+ K_i\otimes E_i,\ \
\Dhat_q(F_i)=F_i\otimes K_i^{-1}+1\otimes F_i.
$$
The action of $U_q\tilde\g$ on $M^\hbar=\oplus_{\lambda\in
P}M^\hbar_\lambda$ will be such that
$$
E_iM^\hbar_\lambda\subset M^\hbar_{\lambda-\alpha_i},\ \
F_iM^\hbar_\lambda\subset M^\hbar_{\lambda+\alpha_i}, \ \
K_i|_{M^\hbar_\lambda}=q_i^{\lambda(h_i)}.
$$
From now on we shall write $\lambda(i)$ instead of $\lambda(h_i)$
to simplify notation.
Therefore $\lambda(1),\dots,\lambda(r)$ are the coefficients of
$\lambda$ in the basis $\omega_1,\dots,\omega_r$.

Recalling that $M^\hbar_\lambda$ is the inverse limit of
$V_{\bar\mu}\otimes V_{\lambda+\mu}$, to define $F_i$ we need
consistent morphisms
$$
V_{\bar\mu+\bar\eta}\otimes V_{\lambda+\mu+\eta}\to
V_{\bar\mu}\otimes V_{\lambda+\alpha_i+\mu}.
$$
These will be defined using morphisms
$V_{\lambda+\mu+\eta}\to V_\eta\otimes V_{\lambda+\alpha_i+\mu}$,
or in other words, morphisms
$$
V_{\mu+\eta-\alpha_i}\to V_\mu\otimes V_\eta.
$$
Up to a scalar there exists only one such morphism. Indeed, if
$\mu(i),\eta(i)\ge1$, then the weight space $(V_\mu\otimes
V_\eta)(\mu+\eta-\alpha_i)$ is spanned by the vectors
$f_i\xi_\mu\otimes\xi_\eta$ and $\xi_\mu\otimes f_i\xi_\eta$. The
vector
$$
\mu(i)\xi_\mu\otimes f_i\xi_\eta-\eta(i)f_i\xi_\mu\otimes \xi_\eta
$$
is the only vector in this space, up to a scalar, which is killed by
$e_i$. The corresponding morphism is defined by
\begin{equation} \label{etau}
\tau_{i;\mu,\eta}\colon V_{\mu+\eta-\alpha_i}\to V_\mu\otimes
V_\eta, \ \ \xi_{\mu+\eta-\alpha_i}\mapsto \mu(i)\xi_\mu\otimes
f_i\xi_\eta-\eta(i)f_i\xi_\mu\otimes \xi_\eta.
\end{equation}
Remark that we also have
\begin{equation} \label{etaul}
\tau_{i;\mu,\eta}(\zeta_{\mu+\eta-\alpha_i})
=-\mu(i)\zeta_\mu\otimes e_{\bar i}\zeta_\eta+\eta(i)e_{\bar
i}\zeta_\mu\otimes \zeta_\eta,
\end{equation}
as can be easily checked using the properties of $\theta$ discussed in
Section~\ref{s4}.

Up to a scalar the morphism $V_{\bar\mu+\bar\eta}\otimes
V_{\lambda+\mu+\eta}\to V_{\bar\mu}\otimes V_{\lambda+\alpha_i+\mu}$
will be defined as the composition
$$
V_{\bar\mu+\bar\eta}\otimes V_{\lambda+\mu+\eta}
\xrightarrow{T_{\bar\mu,\bar\eta}\otimes
\tau_{i;\eta,\lambda+\alpha_i+\mu}} V_{\bar\mu}\otimes
V_{\bar\eta}\otimes V_\eta\otimes V_{\lambda+\alpha_i+\mu}
\xrightarrow{(\iota\otimes S_\eta^\hbar\otimes\iota)B}
V_{\bar\mu}\otimes V_{\lambda+\alpha_i+\mu},
$$
where $B=(\Phi\otimes\iota)\Phi_{12,3,4}^{-1}$. To find the right
normalization we want these maps to define the usual action of $\g$
on the forgetful functor for $\hbar=0$. It is not difficult to check
that for $\hbar=0$ we have to divide the above map by
$\eta(i)$. More importantly, we want the above maps to be consistent
with $\tr^{\cdot,\hbar}$ for all $\hbar$. We are then forced to find
out how the associator $\Phi$ composes with morphisms
$V_{\mu+\eta+\nu-\alpha_i}\to V_\mu\otimes V_\eta\otimes V_\nu$
obtained by combining the maps $\tau$ and $T$. The space of all
possible morphisms is isomorphic to the two-dimensional subspace of
$V_\mu\otimes V_\eta\otimes V_\nu$ of vectors of weight
$\mu+\eta+\nu-\alpha_i$ killed by $e_i$. Therefore we have to
compute the operator $\Phi(A,B)$ for two-by-two matrices $A$ and
$B$. 

\begin{lemma} \label{ltauT}
Let $A=\begin{pmatrix}a+b & 0\\ c & 0\end{pmatrix}$ and
$B=\begin{pmatrix}-b-c & a\\ 0 & 0\end{pmatrix}$ be such that
the numbers $a,b,c,a+b,a+c,b+c,a+b+c$ are non-integral.
Consider the eigenvectors
$$
e_1=\begin{pmatrix}a+b \\ c\end{pmatrix}, \ \ e_2=\begin{pmatrix}0 \\
b\end{pmatrix}
$$
of $A$ and the eigenvectors
$$
f_1=\begin{pmatrix}a \\ b+c\end{pmatrix}, \ \ f_2=\begin{pmatrix}b \\
0\end{pmatrix}
$$
of $B$. Then
$$
\begin{pmatrix}b & -c\\ 0 & b+c\end{pmatrix}
\begin{pmatrix}e_1\\ e_2\end{pmatrix}
=\begin{pmatrix}0 & a+b\\ b & -a\end{pmatrix}
\begin{pmatrix}f_1\\ f_2\end{pmatrix}
$$
and
\begin{align*}
\Phi(A,B)&\begin{pmatrix}\sin\pi b & -\sin\pi c\\ 0 & \sin
\pi(b+c)\end{pmatrix}
\begin{pmatrix}\frac{1}{\Gamma(1+a+b)\Gamma(1+c)
\Gamma(1-(a+b+c))} & 0\\ 0 & \frac{1}{\Gamma(1+a)\Gamma(1+b)
\Gamma(1-(a+b))}\end{pmatrix}
\begin{pmatrix}e_1\\ e_2\end{pmatrix}\\
=&\begin{pmatrix}0 & \sin\pi(a+b)\\ \sin\pi b & -
\sin\pi a\end{pmatrix}
\begin{pmatrix}\frac{1}{\Gamma(1+a)\Gamma(1+b+c)
\Gamma(1-(a+b+c))} & 0\\ 0  & \frac{1}{\Gamma(1+b)\Gamma(1+c)
\Gamma(1-(b+c))}\end{pmatrix}
\begin{pmatrix}f_1\\ f_2\end{pmatrix}.
\end{align*}
\end{lemma}

\bp The equation $v'=\left(\frac{A}{x}+\frac{B}{x-1}\right)v$, where
$v\colon(0,1)\to\C^2$, has the form
$$
\begin{cases}
v_0'=\left(\frac{a+b}{x}-\frac{b+c}{x-1}\right)v_0+\frac{a}{x-1}v_1\\
v_1'=\frac{c}{x}v_0.
\end{cases}
$$
It follows that $u=v_1$ satisfies the Gauss differential equation
\begin{equation} \label{eGauss}
x(1-x)u''+(\gamma-(\alpha+\beta+1)x)u'-\alpha\beta u=0,
\end{equation}
where $\alpha=-a$, $\beta=c$, $\gamma=1-a-b$. Denote by $\Gamma$ the
space of solutions of this equation on~$(0,1)$. By our discussion on
page~\pageref{MonDef} the operator $\Phi(A,B)$ can be written as
$\pi^{-1}_1\pi_0$, where the linear isomorphisms
$\pi_0,\pi_1\colon\C^2\to\Gamma$ are defined as follows. If $\xi$ is
an eigenvector of $A$ with eigenvalue $\lambda$ then $\pi_0(\xi)$ is
the unique solution $u\in\Gamma$ such that the vector valued
function
$$
(0,1)\ni x\mapsto x^{-\lambda}\begin{pmatrix}\frac{x}{c}u'(x)\\
u(x)\end{pmatrix}
$$
extends to a holomorphic function on the unit disc with value
$\xi$ at $x=0$. Similarly, if $\xi$ is an eigenvector of $B$
with eigenvalue $\lambda$ then $\pi_1(\xi)$ is the unique
solution $u\in\Gamma$ such that the vector valued function
$$
(0,1)\ni x\mapsto x^{-\lambda}
\begin{pmatrix}\frac{1-x}{c}u'(1-x)\\ u(1-x)\end{pmatrix}
$$
extends to a holomorphic function on the unit disc with value $\xi$
at $x=0$.

Recall that the Euler hypergeometric function
$F(\alpha,\beta,\gamma;\cdot)$ is the unique solution $u$ of
\eqref{eGauss} which is analytic on the unit disc and such that
$u(0)=1$, $u'(0)=\alpha\beta/\gamma$. Consider the following four
solutions of \eqref{eGauss}:
\begin{align*}
u_1&=x^{1-\gamma}(1-x)^{\gamma-\alpha-\beta}
F(1-\alpha,1-\beta,2-\gamma;x), \\
u_2&=F(\alpha,\beta,\gamma;x),\\
u_3&=F(\alpha,\beta,1+\alpha+\beta-\gamma;1-x),\\
 u_4&=x^{1-\gamma}(1-x)^{\gamma-\alpha-\beta}
F(1-\alpha,1-\beta,1-\alpha-\beta+\gamma;1-x).
\end{align*}
Then it is immediate that the isomorphisms $\pi_0$ and $\pi_1$ are
given by
$$
\pi_0(e_1)=cu_1, \ \ \pi_0(e_2)=bu_2,\ \
\pi_1(f_1)=(b+c)u_3, \ \
\pi_1(f_2)=-\frac{bc}{1-b-c}u_4.
$$
We have the following identity, see e.g.~\cite{Bat}:
$$
\frac{\Gamma(\alpha)\Gamma(\beta)}
{\Gamma(\alpha+\beta-\gamma+1)}u_3
=\frac{\Gamma(\alpha)\Gamma(\beta)\Gamma(1-\gamma)}
{\Gamma(\alpha-\gamma+1)\Gamma(\beta-\gamma+1)}u_2
+\Gamma(\gamma-1)u_1.
$$
Substituting $x$ for $1-x$ and $\gamma$ for
$1+\alpha+\beta-\gamma$ we also get
$$
\frac{\Gamma(\alpha)\Gamma(\beta)}{\Gamma(\gamma)}u_2
=\frac{\Gamma(\alpha)\Gamma(\beta)\Gamma(\gamma-\alpha-\beta)}
{\Gamma(\gamma-\alpha)\Gamma(\gamma-\beta)}u_3
+\Gamma(\alpha+\beta-\gamma)u_4.
$$
A direct but tedious computation using these identities together
with the identities $\Gamma(1+x)=x\Gamma(x)$ and
$\Gamma(x)\Gamma(1-x)=\pi/\sin\pi x$ yields the result. \ep

Define morphisms
\begin{equation} \label{etauh}
\tau^\hbar_{i;\mu,\eta} =\frac{\tau_{i;\mu,\eta}}{\Gamma(1+\hbar
d_i\mu(i)) \Gamma(1+\hbar d_i\eta(i))\Gamma(1-\hbar
d_i(\mu(i)+\eta(i)))} \colon V_{\mu+\eta-\alpha_i}\to V_\mu\otimes
V_\eta,
\end{equation}
where $\tau_{i;\mu,\eta}$ is defined by \eqref{etau}.

The subspace of $V_\mu\otimes V_\eta\otimes V_\nu$ of vectors of
weight $\mu+\eta+\nu-\alpha_i$ killed by $e_i$ is spanned by the
vectors
$$
\eta(i)\xi_\mu\otimes\xi_\eta\otimes f_i\xi_\nu-
\nu(i)\xi_\mu\otimes f_i\xi_\eta\otimes\xi_\nu \ \ \hbox{and} \ \
\mu(i)\xi_\mu\otimes f_i\xi_\eta\otimes \xi_\nu-
\eta(i)f_i\xi_\mu\otimes \xi_\eta\otimes\xi_\nu.
$$
This space is invariant under the operators $t_{12}$ and $t_{23}$.
In the above basis these operators have the form
$$
t_{12}=\begin{pmatrix}
(\mu,\eta) & 0\\ d_i\nu(i) & (\mu,\eta)-d_i\mu(i)-d_i\eta(i)
\end{pmatrix}, \ \
t_{23}=\begin{pmatrix}
(\eta,\nu)-d_i\eta(i)-d_i\nu(i) & d_i\mu(i)\\ 0 & (\eta,\nu)
\end{pmatrix}.
$$
To see this first recall that $t(\xi_\mu\otimes\xi_\eta)=(\mu,\eta)
\xi_\mu\otimes\xi_\eta$ by \eqref{ethigh}. Using $\g$-invariance
of $t$ we therefore get
\begin{equation} \label{etaction1}
t(f_i\xi_\mu\otimes\xi_\eta+\xi_\mu\otimes f_i\xi_\eta)
=f_it(\xi_\mu\otimes\xi_\eta)
=(\mu,\eta)(f_i\xi_\mu\otimes\xi_\eta+\xi_\mu\otimes f_i\xi_\eta).
\end{equation}
Using $(\mu,\alpha_i)=d_i\mu(i)$ and $(\alpha_i,\rho)=d_i$
and arguing as for \eqref{ethigh} we get
\begin{equation} \label{etaction2}
t|_{\tau_{i;\mu,\eta}(V_{\mu+\eta-\alpha_i})}
=(\mu,\eta)-d_i\mu(i)-d_i\eta(i),
\end{equation}
whence
$$
t(\mu(i)\xi_\mu\otimes f_i\xi_\eta-
\eta(i)f_i\xi_\mu\otimes \xi_\eta)=((\mu,\eta)-d_i\mu(i)-d_i\eta(i))
(\mu(i)\xi_\mu\otimes f_i\xi_\eta-
\eta(i)f_i\xi_\mu\otimes \xi_\eta).
$$
By virtue of this identity and \eqref{etaction1} we conclude that
$$
t(f_i\xi_\mu\otimes\xi_\eta)
=((\mu,\eta)-d_i\eta(i))f_i\xi_\mu\otimes\xi_\eta
+d_i\mu(i)\xi_\mu\otimes f_i\xi_\eta.
$$
Applying the flip we also get
$$
t(\xi_\eta\otimes f_i\xi_\mu)
=((\mu,\eta)-d_i\eta(i))\xi_\eta\otimes f_i\xi_\mu
+d_i\mu(i)f_i\xi_\eta\otimes \xi_\mu.
$$
These two identities and \eqref{ethigh} imply the above matrix forms of $t_{12}$ and $t_{23}$.

Recall now that $\Phi(A,B)=\Phi(A-\alpha,B-\beta)$ for any scalars
$\alpha$ and $\beta$. So replacing $t_{12}$ by $t_{12}-
((\mu,\eta)-d_i\mu(i)-d_i\eta(i))1$ and $t_{23}$ by
$t_{23}-(\eta,\nu)1$ we are in a position to apply Lemma~\ref{ltauT}
with $a=\hbar d_i\mu(i)$, $b=\hbar d_i\eta(i)$ and $c=\hbar
d_i\nu(i)$. One checks that the vectors $e_1,e_2,f_1,f_2$ in the
lemma are exactly the images of $\hbar
d_i\eta(i)\xi_{\mu+\eta+\nu-\alpha_i}$ under the morphisms
$(T_{\mu,\eta}\otimes\iota)\tau_{i;\mu+\eta,\nu}$,
$(\tau_{i;\mu,\eta}\otimes\iota)T_{\mu+\eta-\alpha_i,\nu}$,
$(\iota\otimes T_{\eta,\nu})\tau_{i;\mu,\eta+\nu}$ and
$(\iota\otimes\tau_{i;\eta,\nu})T_{\mu,\eta+\nu-\alpha_i}$,
respectively. As
$$
\sin \pi\hbar d_i x=\frac{q_i^x-q_i^{-x}}{2\sqrt{-1}}
=\frac{q_i-q_i^{-1}}{2 \sqrt{-1}}[x]_{q_i},
$$
Lemma~\ref{ltauT} can therefore be reformulated as the following
identity between morphisms $V_{\mu+\eta+\nu-\alpha_i}\to
V_\mu\otimes V_\eta\otimes V_\nu$:
\begin{equation} \label{etauT}
\begin{split}
\Phi(\hbar t_{12},\hbar t_{23})
\begin{pmatrix}[\eta(i)]_{q_i} & -[\nu(i)]_{q_i}\\
0 & [\eta(i)+\nu(i)]_{q_i}\end{pmatrix}&
\begin{pmatrix}(T_{\mu,\eta}\otimes\iota)\tau^\hbar_{i;\mu+\eta,\nu}\\
(\tau^\hbar_{i;\mu,\eta}\otimes\iota)
T_{\mu+\eta-\alpha_i,\nu}\end{pmatrix}\\
&=\begin{pmatrix}0 & [\mu(i)+\eta(i)]_{q_i}\\ [\eta(i)]_{q_i} & -
[\mu(i)]_{q_i}\end{pmatrix}
\begin{pmatrix}(\iota\otimes T_{\eta,\nu})\tau^\hbar_{i;\mu,\eta+\nu}\\
(\iota\otimes\tau^\hbar_{i;\eta,\nu})
T_{\mu,\eta+\nu-\alpha_i}\end{pmatrix}.
\end{split}
\end{equation}
It is remarkable that the proof of this identity is the first and only place
where one uses nontrivial specific properties of $\Phi$ beyond being
an associator; the only special property which we used before
Lemma~\ref{ltauT} was that $\Phi$ acts trivially on the highest
weight space of $V_\mu\otimes V_\eta\otimes V_\nu$.

\begin{proposition} \label{Maction}
The morphisms
\begin{equation} \label{ePsi}
V_{\bar\mu+\bar\eta}\otimes V_{\lambda+\mu+\eta}
\xrightarrow{[\eta(i)]_{q_i}^{-1}T_{\bar\mu,\bar\eta}\otimes
\tau^\hbar_{i;\eta,\lambda+\alpha_i+\mu}} V_{\bar\mu}\otimes
V_{\bar\eta}\otimes V_\eta\otimes V_{\lambda+\alpha_i+\mu}
\xrightarrow{(\iota\otimes S_\eta^\hbar\otimes\iota)B}
V_{\bar\mu}\otimes V_{\lambda+\alpha_i+\mu},
\end{equation}
where $B=(\Phi\otimes\iota)\Phi_{12,3,4}^{-1}$, are consistent
with $\tr^{\cdot,\hbar}$ and
hence define a morphism $F_i\colon M^\hbar_\lambda\to
M^\hbar_{\lambda+\alpha_i}$. Similarly the morphisms
\begin{equation} \label{ePhi}
V_{\bar\mu+\bar\eta}\otimes V_{\lambda+\mu+\eta}
\xrightarrow{[\eta(i)]_{q_i}^{-1}\tau^\hbar_{\bar
i;\bar\mu+\bar\alpha_i,\bar\eta}\otimes  T_{\eta,\lambda+\mu}}
V_{\bar\mu+\bar\alpha_i}\otimes V_{\bar\eta}\otimes V_\eta\otimes
V_{\lambda+\mu} \xrightarrow{(\iota\otimes
S_\eta^\hbar\otimes\iota)B} V_{\bar\mu+\bar\alpha_i}\otimes
V_{\lambda+\mu}
\end{equation}
define a morphism $E_i\colon M^\hbar_\lambda\to
M^\hbar_{\lambda-\alpha_i}$.

Furthermore, for generic $\hbar$ we can choose the $1$-cochain
$g^\hbar_\mu$ such that for each $i$ the composition
\begin{equation} \label{erho}
V_{2\bar\omega_i-\bar\alpha_i}\otimes V_{2\omega_i-\alpha_i}
\xrightarrow{\tau^\hbar_{\bar i;\bar\omega_i,\bar\omega_i}\otimes
\tau^\hbar_{i;\omega_i,\omega_i}}V_{\bar\omega_i} \otimes
V_{\bar\omega_i}\otimes V_{\omega_i}\otimes V_{\omega_i}
\xrightarrow{(\iota\otimes S_{\omega_i}\otimes\iota)B}
V_{\bar\omega_i}\otimes V_{\omega_i}\xrightarrow{S_{\omega_i}}\C
\end{equation}
coincides with $\displaystyle -\frac{g^\hbar_{2\omega_i-\alpha_i}}{
(g^\hbar_{\omega_i})^2}[2]_{q_i}S_{2\omega_i-\alpha_i}$. If
$g^\hbar_\mu$ is chosen this way then the morphisms $E_i$ and $F_i$
together with the morphism $K_i\colon M^\hbar\to M^\hbar$ acting on
$M^\hbar_\lambda$ as multiplication by $q^{\lambda(i)}_i$, define an
action of the algebra $U_q\tilde\g$ on $M^\hbar$. This action
respects the comonoid structure of $M^\hbar$ in the sense that
$\delta^\hbar(\omega x)=\Dhat_q(\omega)\delta^\hbar(x)$ and
$\eps^\hbar(\omega x)=\hat\eps_q(\omega)\eps^\hbar(x)$ for all
$\omega\in U_q\tilde\g$ and $x\in M^\hbar$.
\end{proposition}

\bp Denote the morphism \eqref{ePsi} by
$\Psi^{\eta,\hbar}_{i;\mu,\lambda+\alpha_i+\mu}$. For consistency we
have to check that
$$
\tr^{\eta,\hbar}_{\mu,\lambda+\alpha_i+\mu}
\Psi^{\nu,\hbar}_{i;\mu+\eta,\lambda+\alpha_i+\mu+\eta}
=\Psi^{\eta+\nu,\hbar}_{i;\mu,\lambda+\alpha_i+\mu}\ \
\hbox{and}\ \
\Psi^{\eta,\hbar}_{i;\mu,\lambda+\alpha_i+\mu}
\tr^{\nu,\hbar}_{\mu+\eta,\lambda+\alpha_i+\mu+\eta}
=\Psi^{\eta+\nu,\hbar}_{i;\mu,\lambda+\alpha_i+\mu}.
$$
We shall only check the first identity. Once again we strictify
$\D(\g,\hbar)$. As we have done before, we shall often skip the
lower indices of maps when they are determined by the target
modules. By definition we have
$$
\Psi^{\eta+\nu,\hbar}_{i;\mu,\lambda+\alpha_i+\mu}
=\frac{1}{[\eta(i)+\nu(i)]_{q_i}}(\iota\otimes
S^\hbar_{\eta+\nu}\otimes\iota)(T\otimes\tau^\hbar_i).
$$
Using that $S^\hbar_{\eta+\nu}=S^\hbar_\eta(\iota\otimes
S^\hbar_{\nu}\otimes\iota)(T\otimes T_{\nu,\eta})$ by definition, we
can rewrite the right hand side as
\begin{equation} \label{ePsitr0}
\frac{1}{[\eta(i)+\nu(i)]_{q_i}}(\iota\otimes
S^\hbar_\eta\otimes\iota)(\iota\otimes\iota\otimes
S^\hbar_{\nu}\otimes\iota\otimes\iota)(\iota\otimes T\otimes
T_{\nu,\eta}\otimes\iota)(T\otimes\tau^\hbar_i).
\end{equation}
On the other hand,
\begin{align*}
\tr^{\eta,\hbar}_{\mu,\lambda+\alpha_i+\mu}
\Psi^{\nu,\hbar}_{i;\mu+\eta,\lambda+\alpha_i+\mu+\eta}
&=\frac{1}{[\nu(i)]_{q_i}}(\iota\otimes S^\hbar_\eta\otimes\iota)
(T\otimes T)(\iota\otimes
S^\hbar_{\nu}\otimes\iota)(T\otimes\tau^\hbar_i)\\
&=\frac{1}{[\nu(i)]_{q_i}}(\iota\otimes S^\hbar_\eta\otimes\iota)
(\iota\otimes\iota\otimes S^\hbar_{\nu}\otimes\iota\otimes\iota)
(T\otimes\iota\otimes \iota\otimes T)(T\otimes\tau^\hbar_i)
\end{align*}
Using that $(T\otimes\iota)T=(\iota\otimes T)T$ by \eqref{eTT} we
can rewrite this as
\begin{equation} \label{ePsitr1}
\frac{1}{[\nu(i)]_{q_i}}(\iota\otimes S^\hbar_\eta\otimes\iota)
(\iota\otimes\iota\otimes S^\hbar_{\nu}\otimes\iota\otimes\iota)
(\iota\otimes T\otimes \iota\otimes T)(T\otimes\tau^\hbar_i).
\end{equation}
It follows from \eqref{etauT} that up to a scalar factor the
difference
$$
\frac{1}{[\eta(i)+\nu(i)]_{q_i}}
(T_{\nu,\eta}\otimes\iota)\tau^\hbar_i
-\frac{1}{[\nu(i)]_{q_i}}(\iota\otimes T)\tau^\hbar_i
$$
(in our strictified category) is equal to
$(\tau^\hbar_i\otimes\iota)T$. Therefore to show that
\eqref{ePsitr0} and \eqref{ePsitr1} are equal we have to check that
the morphism
$$
(\iota\otimes S^\hbar_\eta\otimes\iota)
(\iota\otimes\iota\otimes S^\hbar_{\nu}\otimes\iota\otimes\iota)
(\iota\otimes T\otimes\tau^\hbar_i\otimes\iota)(T\otimes T)
\colon V_{\bar\mu+\bar\eta+\bar\nu}
\otimes V_{\lambda+\mu+\eta+\nu}
\to V_{\bar\mu}\otimes V_{\lambda+\alpha_i+\mu}
$$
is zero. In fact already $S^\hbar_\eta (\iota\otimes
S^\hbar_{\nu}\otimes\iota) (T\otimes\tau^\hbar_i)=0$ since zero is
the only morphism from $V_{\bar\eta+\bar\nu}\otimes
V_{\eta+\nu-\alpha_i}$ to $\C$.

Thus $F_i$ is well-defined. Similarly one proves that $E_i$ is
well-defined.

\smallskip

Next we have to check that under a specific choice of $g^\hbar_\mu$
the morphisms $E_i,F_i,K_i$ satisfy the defining relations of
$U_q\tilde\g$. The only nontrivial relation is
\begin{equation} \label{ecomm}
E_iF_j-F_jE_i=-\delta_{ij}\frac{K_i-K_i^{-1}}{q_i-q_i^{-1}}.
\end{equation}
The rest clearly holds without any assumptions on the cochain by
using that $\alpha_j(i)=a_{ij}$.

The composition \eqref{erho} coincides with $S_{2\omega_i-\alpha_i}$
up to a scalar factor since the space of morphisms
$V_{2\bar\omega_i-\bar\alpha_i}\otimes V_{2\omega_i-\alpha_i}\to\C$
is one-dimensional. This factor is nonzero for generic $\hbar$ since
it is equal to~$-2$ for $\hbar=0$. Indeed, by virtue of
\eqref{etau}-\eqref{etaul} we have to show that the image of
$$
(-\zeta_{\bar\omega_i}\otimes e_i\zeta_{\bar\omega_i}
+e_i\zeta_{\bar\omega_i}\otimes\zeta_{\bar\omega_i})
\otimes (\xi_{\omega_i}\otimes f_i\xi_{\omega_i}
-f_i\xi_{\omega_i}\otimes\xi_{\omega_i})
$$
under the map $S_{\omega_i}(\iota\otimes S_{\omega_i}\otimes\iota)$
equals $-2$. This follows from
$S_{\omega_i}(e_i\zeta_{\bar\omega_i}\otimes f_i\xi_{\omega_i})=-1$
used twice, which in turn follows from the identities
$$
e_i\zeta_{\bar\omega_i}\otimes f_i\xi_{\omega_i}
=e_i(\zeta_{\bar\omega_i}\otimes f_i\xi_{\omega_i})-
\zeta_{\bar\omega_i}\otimes e_if_i\xi_{\omega_i}
=e_i(\zeta_{\bar\omega_i}\otimes f_i\xi_{\omega_i})-
\zeta_{\bar\omega_i}\otimes \xi_{\omega_i}.
$$
So to show that we can make the specific choice of the cochain
$g^\hbar_\mu$ stated in the formulation we just have to check that
the ratios $g^\hbar_{2\omega_i-\alpha_i}/ (g^\hbar_{\omega_i})^2$
can take arbitrary values. As we already remarked in the proof of
Lemma~\ref{lunicom}, the cochain $g^\hbar_\mu$ is defined up to
multiplication by a homomorphism $\chi\colon P\to\C^*$. If we
replace $g^\hbar_\mu$ by $g^\hbar_\mu\chi(\mu)$ then
$g^\hbar_{2\omega_i-\alpha_i}/ (g^\hbar_{\omega_i})^2$ changes by
the factor $\chi(\alpha_i)^{-1}$. Therefore the claim follows from
the fact that any homomorphism from the root lattice $Q$ into $\C^*$
can be extended to the weight lattice~$P$. This is well-known and
easy to see using infinite divisibility of $\C^*$.

Assuming now that the cochain $g^\hbar_\mu$ is chosen as stated we
want to check \eqref{ecomm}. Denoting the composition \eqref{ePhi}
by $\Phi^{\eta,\hbar}_{i;\mu+\alpha_i,\lambda+\mu}$, to prove
\eqref{ecomm} for $i=j$ it suffices to show that
\begin{equation}  \label{ecomm1}
\Phi^{\omega_i,\hbar}_{i;\mu+\alpha_i,\lambda+\alpha_i+\mu}
\Psi^{\omega_i,\hbar}_{i;\mu+\omega_i,\lambda+\alpha_i+\mu+\omega_i}-
\Psi^{\omega_i,\hbar}_{i;\mu+\alpha_i,\lambda+\alpha_i+\mu}
\Phi^{\omega_i,\hbar}_{i;\mu+\alpha_i+\omega_i,\lambda+\mu+\omega_i}
=-[\lambda(i)]_{q_i}
\tr^{2\omega_i-\alpha_i,\hbar}_{\mu+\alpha_i,\lambda+\mu+\alpha_i}.
\end{equation}
The first term on the left hand side in our strictified category is
$$
(\iota\otimes S^\hbar_{\omega_i}\otimes\iota) (\tau^\hbar_{\bar
i}\otimes T) (\iota\otimes S^\hbar_{\omega_i}\otimes\iota) (T\otimes
\tau^\hbar_i) =(\iota\otimes S^\hbar_{\omega_i}\otimes\iota)
(\iota\otimes\iota\otimes
S^\hbar_{\omega_i}\otimes\iota\otimes\iota) (\tau^\hbar_{\bar
i}\otimes\iota\otimes\iota\otimes T)(T\otimes \tau^\hbar_i).
$$
Expressing similarly the second term we get that the left hand side
of \eqref{ecomm1} equals
$$
(\iota\otimes S^\hbar_{\omega_i}\otimes\iota)
(\iota\otimes\iota\otimes
S^\hbar_{\omega_i}\otimes\iota\otimes\iota) \big((\tau^\hbar_{\bar
i}\otimes\iota\otimes\iota\otimes T)(T\otimes \tau^\hbar_i)
-(T\otimes\iota\otimes\iota\otimes\tau^\hbar_i)(\tau^\hbar_{\bar
i}\otimes T)\big).
$$
Next we use identities \eqref{etauT} to express $(\tau^\hbar_{\bar
i}\otimes\iota\otimes\iota\otimes T)(T\otimes \tau^\hbar_i)$ and
$(T\otimes\iota\otimes\iota\otimes\tau^\hbar_i)(\tau^\hbar_{\bar
i}\otimes T)$ in the form $(\iota\otimes *\otimes
*\otimes\iota)(*\otimes*)$. A tedious but straightforward
computation keeping track of subindices shows that the terms
$(\iota\otimes T\otimes T\otimes\iota)(\tau^\hbar_{\bar i}\otimes
\tau^\hbar_i)$ cancel, and what is left is the term
$\frac{[\lambda(i)]_{q_i}}{[2]_{q_i}}(\iota\otimes\tau^\hbar_{\bar
i}\otimes \tau^\hbar_i \otimes\iota)(T\otimes T)$ and scalar
multiples of $(\iota\otimes T\otimes
\tau^\hbar_i\otimes\iota)(\tau^\hbar_{\bar i}\otimes T)$ and
$(\iota\otimes \tau^\hbar_{\bar i}\otimes T \otimes\iota)(T\otimes
\tau^\hbar_i)$. The last two terms vanish when composed with
$(\iota\otimes S^\hbar_{\omega_i}\otimes\iota)
(\iota\otimes\iota\otimes
S^\hbar_{\omega_i}\otimes\iota\otimes\iota)$ for the same reason as
in the proof of consistency of $\Psi^{\cdot,\hbar}$. Therefore the
left hand side of \eqref{ecomm1} equals
$$
\frac{[\lambda(i)]_{q_i}}{[2]_{q_i}}(\iota\otimes
S^\hbar_{\omega_i}\otimes\iota) (\iota\otimes\iota\otimes
S^\hbar_{\omega_i}\otimes\iota\otimes\iota)
(\iota\otimes\tau^\hbar_{\bar i}\otimes \tau^\hbar_i
\otimes\iota)(T\otimes T).
$$
We have $S^\hbar_{\omega_i}(\iota\otimes
S^\hbar_{\omega_i}\otimes\iota) (\tau^\hbar_{\bar i}\otimes
\tau^\hbar_i)=-[2]_{q_i}S^\hbar_{2\omega_i-\alpha_i}$ by our choice
of the cochain $g^\hbar_\mu$, so the above expression is equal to
$$
-[\lambda(i)]_{q_i}(\iota\otimes S^\hbar_{2\omega_i-\alpha_i}
\otimes\iota)(T\otimes T),
$$
which by definition is the right hand side of \eqref{ecomm1}.

The relation $E_iF_j-F_jE_i=0$ for $i\ne j$ is proved similarly by
showing that
$$
\Phi^{\eta,\hbar}_{i;\mu+\alpha_i,\lambda+\alpha_j+\mu}
\Psi^{\nu,\hbar}_{j;\mu+\eta,\lambda+\alpha_j+\mu+\eta}-
\Psi^{\eta,\hbar}_{j;\mu+\alpha_i,\lambda+\alpha_j+\mu}
\Phi^{\nu,\hbar}_{i;\mu+\alpha_i+\eta,\lambda+\mu+\eta}=0.
$$
We only remark that in this case the morphism
$S^\hbar_\eta(\iota\otimes S^\hbar_\nu\otimes\iota)(\tau^\hbar_{\bar
i}\otimes \tau^\hbar_j)$ vanishes as there are no nonzero morphisms
$V_{\bar\eta+\bar\nu-\bar\alpha_i}\otimes
V_{\eta+\nu-\alpha_j}\to\C$.

\smallskip

It remains to show that the action of $U_q\tilde\g$ respects the
comonoid structure of $M^\hbar$. We shall only check that
$\delta^\hbar(F_ix)=\Dhat_q(F_i)\delta^\hbar(x)$, that is,
$$
\delta^\hbar_{\lambda_1,\lambda_2}F_i
=q_i^{-\lambda_2(i)}(F_i\otimes\iota)
\delta^\hbar_{\lambda_1-\alpha_i,\lambda_2}+(\iota\otimes
F_i)\delta^\hbar_{\lambda_1,\lambda_2-\alpha_i}.
$$
The morphisms $\delta^\hbar$ are induced by the morphisms $m^\hbar$
defined by \eqref{em}. Therefore it suffices to check that
\begin{equation}\label{emPsi}
\begin{split}
m^\hbar_{\mu,\eta,\lambda_1,\lambda_2}&
\Psi^{\nu,\hbar}_{i;\mu+\eta,\lambda_1+\lambda_2+\mu+\eta}\\
&=q_i^{-\lambda_2(i)}
(\Psi^{\nu,\hbar}_{i;\mu,\lambda_1+\mu}\otimes\iota\otimes\iota)
m^\hbar_{\mu+\nu,\eta,\lambda_1-\alpha_i,\lambda_2}
+(\iota\otimes\iota\otimes\Psi^{\nu,\hbar}_{i;\eta,\lambda_2+\eta})
m^\hbar_{\mu,\eta+\nu,\lambda_1,\lambda_2-\alpha_i}.
\end{split}
\end{equation}
The left hand side multiplied by $[\nu(i)]_{q_i}$ in our strictified
category with braiding $\sigma$ is
\begin{multline*}
q^{(\lambda_1+\mu,\eta)}(\iota\otimes\sigma\otimes\iota)
(T_{\bar\mu,\bar\eta}\otimes T_{\lambda_1+\mu,\lambda_2+\eta})
(\iota\otimes S^\hbar_\nu\otimes\iota)
(T_{\bar\mu+\bar\eta,\bar\nu}\otimes
\tau^\hbar_{i;\nu,\lambda_1+\lambda_2+\mu+\eta})\\
=q^{(\lambda_1+\mu,\eta)}(\iota\otimes\sigma\otimes\iota)
(\iota\otimes\iota\otimes S^\hbar_\nu\otimes\iota\otimes\iota)
(T_{\bar\mu,\bar\eta}\otimes \iota\otimes\iota\otimes
T_{\lambda_1+\mu,\lambda_2+\eta})
(T_{\bar\mu+\bar\eta,\bar\nu}\otimes
\tau^\hbar_{i;\nu,\lambda_1+\lambda_2+\mu+\eta}).
\end{multline*}
We claim that
\begin{equation} \label{etauC0}
\begin{split}
&(\iota\otimes T_{\lambda_1+\mu,\lambda_2+\eta})
\tau^\hbar_{i;\nu,\lambda_1+\lambda_2+\mu+\eta}\\
&=q_i^{-\lambda_2(i)-\eta(i)}
(\tau^\hbar_{i;\nu,\lambda_1+\mu}\otimes\iota)
T_{\lambda_1-\alpha_i+\mu+\nu,
\lambda_2+\eta}+q^{(\lambda_1+\mu,\nu)} (\sigma^{-1}\otimes\iota)
(\iota\otimes\tau^\hbar_{i;\nu,\lambda_2+\eta})
T_{\lambda_1+\mu,\lambda_2-\alpha_i+\eta+\nu}.
\end{split}
\end{equation}
We postpone the proof of this equality. Using it we see that the
left hand side of \eqref{emPsi} multiplied by $[\nu(i)]_{q_i}$ is
the sum of the term
$$
q^{(\lambda_1+\mu,\eta)}q_i^{-\lambda_2(i)-\eta(i)}
(\iota\otimes\sigma\otimes\iota) (\iota\otimes\iota\otimes
S^\hbar_\nu\otimes\iota\otimes\iota) (T_{\bar\mu,\bar\eta}\otimes
\iota\otimes\tau^\hbar_{i;\nu,\lambda_1+\mu}\otimes\iota)
(T_{\bar\mu+\bar\eta,\bar\nu}\otimes T_{\lambda_1-\alpha_i+\mu+\nu,
\lambda_2+\eta})
$$
\begin{align}
&=q^{(\lambda_1+\mu+\nu,\eta)}q_i^{-\lambda_2(i)-\eta(i)}
(\iota\otimes\sigma\otimes\iota)
(\iota\otimes\iota\otimes S^\hbar_\nu\otimes\iota\otimes\iota)
%\nonumber\\&\quad\quad\quad\quad
%(\iota\otimes\sigma^{-1}\otimes\iota\otimes\iota\otimes\iota)
(\iota\otimes
\sigma^{-1}T_{\bar\nu,\bar\eta}\otimes\tau^\hbar_i\otimes\iota)
(T_{\bar\mu,\bar\eta+\bar\nu}\otimes T)\nonumber\\
&=q^{(\lambda_1+\mu+\nu,\eta)}q_i^{-\lambda_2(i)-\eta(i)}
(\iota\otimes\sigma\otimes\iota)
(\iota\otimes\iota\otimes S^\hbar_\nu\otimes\iota\otimes\iota)
%\nonumber\\&\quad\quad\quad\quad
(\iota\otimes\sigma^{-1}\otimes\iota\otimes\iota\otimes\iota)
(T\otimes\iota\otimes\tau^\hbar_i\otimes\iota) \label{emPsi1}
(T\otimes T),
\end{align}
where we have used $(T\otimes \iota)T=(\iota\otimes T)T$ twice and
that $T_{\bar\eta,\bar\nu}=q^{(\eta,\nu)}
\sigma^{-1}T_{\bar\nu,\bar\eta}$ by \eqref{eTT} and \eqref{eTC}, and
the term
\begin{equation} \label{emPsi2}
\begin{split}
&q^{(\lambda_1+\mu,\eta+\nu)}(\iota\otimes\sigma\otimes\iota)
(\iota\otimes\iota\otimes S^\hbar_\nu\otimes\iota\otimes\iota)
(\iota\otimes\iota\otimes\iota\otimes\sigma^{-1}\otimes\iota)
(T\otimes \iota\otimes\iota\otimes\tau^\hbar_i)
(T\otimes T)\\
&=q^{(\lambda_1+\mu,\eta+\nu)}(\iota\otimes\sigma\otimes\iota)
(\iota\otimes\iota\otimes S^\hbar_\nu\otimes\iota\otimes\iota)
(\iota\otimes\iota\otimes\iota\otimes\sigma^{-1}\otimes\iota)
(\iota\otimes T\otimes\iota\otimes\tau^\hbar_i)
(T\otimes T).
\end{split}
\end{equation}
On the other hand, the first term on the right hand side of
\eqref{emPsi} multiplied by $[\nu(i)]_{q_i}$ equals
$$
q_i^{-\lambda_2(i)} q^{(\lambda_1-\alpha_i+\mu+\nu,\eta)}
(\iota\otimes S^\hbar_\nu\otimes\iota\otimes\iota\otimes\iota)
(T\otimes\tau^\hbar_i\otimes\iota\otimes\iota)
(\iota\otimes\sigma\otimes\iota)(T\otimes T).
$$
By naturality of $\sigma$ this expression can be written as
$$
q_i^{-\lambda_2(i)}q^{(\lambda_1-\alpha_i+\mu+\nu,\eta)}
(\iota\otimes S^\hbar_\nu\otimes\iota\otimes\iota\otimes\iota)
(\iota\otimes\iota\otimes\sigma_{1,23}\otimes\iota)
(T\otimes\iota\otimes \tau^\hbar_i\otimes\iota)(T\otimes T).
$$
As $(\alpha_i,\eta)=d_i\eta(i)$, to see that this is equal to
\eqref{emPsi1} we just have to check that
$$
\sigma(\iota\otimes S^\hbar_\nu\otimes\iota)
(\sigma^{-1}\otimes\iota\otimes\iota)
=(S^\hbar_\nu\otimes\iota\otimes\iota)(\iota\otimes\sigma_{1,23}).
$$
Writing $\sigma\colon U\otimes V\to V\otimes U$ as
$(\iota\otimes\sigma)(\sigma\otimes\iota)\colon U\otimes\C\otimes V
\to\C\otimes V\otimes U$ and using naturality of $\sigma$ we have
$$
\sigma(\iota\otimes S^\hbar_\nu\otimes\iota)
=(S^\hbar_\nu\otimes\iota\otimes\iota)
(\iota\otimes\iota\otimes\sigma) (\sigma_{1,23}\otimes\iota).
$$
As $(\iota\otimes\iota\otimes\sigma)(\sigma_{1,23}\otimes\iota)
=(\iota\otimes\sigma_{1,23})(\sigma\otimes\iota\otimes\iota)$
by the hexagon identities, we get the required equality.

Similarly it is proved that \eqref{emPsi2} coincides with the second
term on the right hand side of \eqref{emPsi} multiplied by
$[\nu(i)]_{q_i}$.

\smallskip

Therefore it remains to check identity \eqref{etauC0}. Replacing
$\lambda_1+\mu$ by $\mu$ and $\lambda_2+\eta$ by $\eta$, we have to
show that
$$
(\iota\otimes T_{\mu,\eta}) \tau^\hbar_{i;\nu,\mu+\eta}
=q_i^{-\eta(i)}(\tau^\hbar_{i;\nu,\mu}\otimes\iota)
T_{\mu+\nu-\alpha_i, \eta}+q^{(\mu,\nu)}(\sigma^{-1}\otimes\iota)
(\iota\otimes\tau^\hbar_{i;\nu,\eta})T_{\mu,\eta+\nu-\alpha_i}.
$$
It follows from identities \eqref{etauT} that
$$
[\mu(i)+\nu(i)]_{q_i}(\iota\otimes T_{\mu,\eta})
\tau^\hbar_{i;\nu,\mu+\eta}
=[\nu(i)]_{q_i}(T_{\nu,\mu}\otimes\iota)\tau^\hbar_{i;\mu+\nu,\eta}
+[\mu(i)+\eta(i)+\nu(i)]_{q_i}(\tau^\hbar_{i;\nu,\mu}\otimes\iota)
T_{\mu+\nu-\alpha_i,\eta}.
$$
Therefore we equivalently have to check that
\begin{multline*}
[\nu(i)]_{q_i}(\sigma T_{\nu,\mu}\otimes\iota)
\tau^\hbar_{i;\mu+\nu,\eta} +[\mu(i)+\eta(i)+\nu(i)]_{q_i}
(\sigma\tau^\hbar_{i;\nu,\mu}\otimes\iota)
T_{\mu+\nu-\alpha_i,\eta}\\
=q_i^{-\eta(i)}[\mu(i)+\nu(i)]_{q_i}
(\sigma\tau^\hbar_{i;\nu,\mu}\otimes\iota)T_{\mu+\nu-\alpha_i,
\eta}+q^{(\mu,\nu)}[\mu(i)+\nu(i)]_{q_i}
(\iota\otimes\tau^\hbar_{i;\nu,\eta})T_{\mu,\eta+\nu-\alpha_i}.
\end{multline*}
But up to the factor $q^{(\mu,\nu)}$ this is exactly the identity
$$
[\nu(i)]_{q_i}(T_{\mu,\nu}\otimes\iota)\tau^\hbar_{i;\mu+\nu,\eta}
-[\eta(i)]_{q_i}(\tau^\hbar_{i;\mu,\nu}\otimes\iota)
T_{\mu+\nu-\alpha_i,\eta} =[\mu(i)+\nu(i)]_{q_i}
(\iota\otimes\tau^\hbar_{i;\nu,\eta})T_{\mu,\eta+\nu-\alpha_i}
$$
from \eqref{etauT}, if we take into account that $\sigma
T_{\nu,\mu}=q^{(\mu,\nu)}T_{\mu,\nu}$ by \eqref{eTC} and
$$
%\begin{equation} \label{etauC}
\sigma\tau^\hbar_{i;\nu,\mu}=-q^{(\mu,\nu)}q_i^{-\mu(i)-\nu(i)}
\tau^\hbar_{i;\mu,\nu},
%\end{equation}
$$
which in turn follows from \eqref{etaction2}
and $\Sigma\tau_{i;\nu,\mu} =-\tau_{i;\mu,\nu}$. \ep

\begin{remark} \label{cochainchange}
If we replace the cochain $g^\hbar_\mu$ by the cochain
$g^\hbar_\mu\chi(\mu)$, where $\chi\colon P\to\C^*$ is a
homomorphism, then by Lemma~\ref{lunicom} the comonoid remains
unaltered up to an isomorphism. One can easily check that if use the
same formulas to define the morphisms $F_i$ and $E_i$ with the new
cochain then the morphism $F_i$ remains unchanged, while $E_i$
changes to $\chi(\alpha_i)E_i$.
\end{remark}

\begin{lemma} \label{uqgaction}
Let $V$ be a finite dimensional $\g$-module. Assume the cochain
$g^\hbar_\mu$ is chosen as in Proposition~\ref{Maction}. Then for
generic $\hbar$ the action of $U_q\tilde\g$ on $M^\hbar$ defines an
action of $U_q\g$ on $\Hom_\g(M^\hbar,V)$.
\end{lemma}

\bp The action of $U_q\tilde\g$ on $M^\hbar$ by $\g$-endomorphisms
defines an action of $(U_q\tilde\g)^{op}$ on $\Hom_\g(M^\hbar,V)$.
To show that this action defines an action of $U_q\g$ we just have
to check that the relations
\begin{equation} \label{eSerre}
\sum^{1-a_{ij}}_{k=0}(-1)^k\begin{bmatrix}1-a_{ij}\\
k\end{bmatrix}_{q_i} E^k_iE_jE^{1-a_{ij}-k}_i=0\ \ \hbox{and}\ \
\sum^{1-a_{ij}}_{k=0}(-1)^k\begin{bmatrix}1-a_{ij}\\
k\end{bmatrix}_{q_i} F^k_iF_jF^{1-a_{ij}-k}_i=0
\end{equation}
are satisfied for $i\ne j$.

We may assume that $V=V_\lambda$ for some $\lambda$. The morphisms
$\tr^{\mu,\hbar}_{0,\lambda}\colon V_{\bar\mu}\otimes
V_{\lambda+\mu}\to V_0\otimes V_\lambda=V_\lambda$ define a morphism
$\xi^\hbar_\lambda\colon M^\hbar_\lambda\to V_\lambda$, which we
consider as a vector in $\Hom_\g(M^\hbar,V_\lambda)$. We have
$E_i\xi^\hbar_\lambda=\xi^\hbar_\lambda\circ E_i=0$ as there are no
nonzero morphisms $M^\hbar_{\lambda+\alpha_i}\to V_\lambda$, so
$\xi^\hbar_\lambda$ is a highest weight vector in
$\Hom_\g(M^\hbar,V_\lambda)$. In particular, if we denote by
$G_{ij}\in (U_q\tilde\g)^{op}$ the left hand side of the first equation in \eqref{eSerre}
then $G_{ij}\xi^\hbar_\lambda=0$. Using the relations in $U_q\tilde\g$
it can be checked that $G_{ij}$ commutes with~$F_l$ for all $l$. Therefore to
prove that $G_{ij}=0$ it suffices to show that
$\Hom_\g(M^\hbar,V_\lambda)$ is spanned by $F_{i_1}\dots
F_{i_m}\xi^\hbar_\lambda=\xi^\hbar_\lambda\circ
F_{i_m}\circ\dots\circ F_{i_1}$. By Remark~\ref{cochainchange} the
latter property is independent of the choice of $g^\hbar_\mu$, so we
may assume that $g^\hbar_\mu$ is an analytic function in $\hbar$
with $g^0_\mu=1$, e.g. by choosing $g^\hbar_{\omega_k}=1$ for all
$k$.

Choose a finite set $I$ of multiindices $(i_1,\dots,i_m)$ such that
the vectors $f_{i_1}\dots f_{i_m}\xi_\lambda$ form a basis of
$V_\lambda$. Since $\dim\Hom_\g(M^\hbar,V_\lambda)\le\dim V_\lambda$
it then suffices to check that for generic $\hbar$ the vectors
$F_{i_1}\dots F_{i_m}\xi^\hbar_\lambda$, $(i_1,\dots,i_m)\in I$, are
linearly independent. The vectors
$$
F_{i_1}\dots F_{i_m}\xi^\hbar_\lambda\in
\Hom_\g(M^\hbar_{\lambda-\alpha_{i_1}-\dots-\alpha_{i_m}},V_\lambda)
$$ are defined by morphisms $V_{\bar\mu}\otimes
V_{\lambda-\alpha_{i_1}-\dots-\alpha_{i_m}+\mu}\to V_\lambda$.
Therefore it suffices to check that the latter morphisms are
linearly independent for generic $\hbar$. Since they depend
analytically on $\hbar$, it is enough to check linear independence
for $\hbar=0$. Under the injective maps $\Hom_\g(V_{\bar\mu}\otimes
V_{\lambda-\eta+\mu}, V_\lambda)\to V_\lambda(\lambda-\eta)$, $f\mapsto
f(\zeta_{\bar\mu}\otimes\xi_{\lambda-\eta+\mu})$, the morphisms are
mapped onto the vectors $f_{i_1}\dots f_{i_m}\xi_\lambda$, which are
linearly independent by assumption. To see that we indeed get the
vectors $f_{i_1}\dots f_{i_m}\xi_\lambda$ we just have to observe
that $\tr^\mu_{0,\lambda}\colon V_{\bar\mu}\otimes
V_{\lambda+\mu}\to V_\lambda$ is mapped onto $\xi_\lambda$ and that
the diagrams
$$
\xymatrix{
\Hom_\g(V_{\bar\mu}\otimes V_{\nu+\alpha_k+\mu},V_\lambda)
\ar[rr]^{\circ\Psi^{\eta,0}_{k;\mu,\nu+\alpha_k+\mu}} \ar[d]& &
\Hom_\g(V_{\bar\mu+\bar\eta}\otimes V_{\nu+\mu+\eta},V_\lambda)\ar[d]\\
V_\lambda(\nu+\alpha_k)\ar[rr]^{f_k} & & V_\lambda(\nu)}
$$
commute, where the top arrow is defined by the morphism
$$
\Psi^{\eta,0}_{k;\mu,\nu+\alpha_k+\mu}\colon
V_{\bar\mu+\bar\eta}\otimes V_{\nu+\mu+\eta}\to V_{\bar\mu}\otimes
V_{\nu+\alpha_k+\mu}
$$
given by~\eqref{ePsi} (with $\hbar=0$ and $g^0_\mu=1$).

Therefore we have proved the first relation in \eqref{eSerre}. The
second is proved similarly by considering the lowest weight vector
$\zeta^\hbar_\lambda\in\Hom_\g(M^\hbar_{-\bar\lambda},V_\lambda)$
defined by $\tr^{\mu-\bar\lambda,\hbar}_{\bar\lambda,0}\colon
V_{\bar\mu}\otimes V_{-\bar\lambda+\mu}\to V_\lambda\otimes
V_0=V_\lambda$. \ep

Thus for generic $\hbar$ we have a well-defined action of $U_q\g$ on
$\Hom_\g(M^\hbar,V)$, so $\Hom_\g(M^\hbar,\cdot)$ can be considered
as a functor $\D(\g,\hbar)\to\CC(\g,\hbar)$. By
Proposition~\ref{ptensor} and the last part of
Proposition~\ref{Maction} it is a tensor functor. Furthermore, by
Proposition~\ref{represento} for generic $\hbar$ the module
$M^\hbar$ is isomorphic to the module $M$ representing the forgetful
functor. Therefore the following theorem finishes the proof of
Theorem~\ref{KL} and thus also of Theorem~\ref{Dr}.

\begin{theorem}
If the cochain $g^\hbar_\mu$ is chosen as in
Proposition~\ref{Maction} then for generic $\hbar$ and $q=e^{\pi
i\hbar}$ the functor $\Hom_\g(M^\hbar,\cdot)$ is a $\C$-linear
braided monoidal equivalence of the categories $\D(\g,\hbar)$ and
$\CC(\g,\hbar)$. This functor maps an irreducible $\g$-module with
highest weight $\lambda$ onto an irreducible $U_q\g$-module with
highest weight $\lambda$.
\end{theorem}

\bp We have already proved that for generic $\hbar$ the functor
$F^\hbar=\Hom_\g(M^\hbar,\cdot)$ is a tensor functor. Furthermore,
by the proof of Lemma~\ref{uqgaction} for any $\lambda\in P$ the
$U_q\g$-module $F^\hbar(V_\lambda)$ has a highest weight
vector~$\xi^\hbar_\lambda$ of weight $\lambda$. Since the dimension
of this module is not bigger than that of $V_\lambda$, we conclude
that $F^\hbar(V_\lambda)$ must be an irreducible $U_q\g$-module with
highest weight $\lambda$. Therefore the image of the functor
contains all simple objects in $\CC(\g,\hbar)$ up to isomorphism. Since
the functor $F^\hbar$ respects direct sums, we conclude that it is
an equivalence of tensor categories.

It remains to check that the functor respects braiding, that is, the
diagram
$$
\xymatrix{
F^\hbar(U)\otimes F^\hbar(V)
\ar[r]^{\Sigma\RR_\hbar} \ar[d]_{F^\hbar_2} &
F^\hbar(V)\otimes F^\hbar(U) \ar[d]^{F^\hbar_2}\\
F^\hbar(U\otimes V) \ar[r]^{F^\hbar(\Sigma q^t)} &
F^\hbar(V\otimes U)
}
$$
commutes. It suffices to consider $U=V_{\bar\lambda}$ and $V=V_\mu$.
Consider the lowest weight vector $\zeta^\hbar_{\bar\lambda}\in
F^\hbar(V_{\bar\lambda})$ and the highest weight vector
$\xi^\hbar_\mu\in F^\hbar(V_\mu)$ defined in the proof of
Lemma~\ref{uqgaction}. It suffices to compute how the morphisms in
the above diagram act on $\zeta^\hbar_{\bar\lambda}\otimes
\xi^\hbar_\mu$. By \eqref{eRmat} we have
$\RR_\hbar(\zeta^\hbar_{\bar\lambda}\otimes
\xi^\hbar_\mu)=q^{-(\lambda,\mu)}\zeta^\hbar_{\bar\lambda}\otimes
\xi^\hbar_\mu$. Recalling that $F^\hbar_2$ is defined using
$\delta^\hbar\colon M^\hbar\to M^\hbar\hat\otimes M^\hbar$, we just
have to check that
$$
q^{-(\lambda,\mu)}(\xi^\hbar_\mu\otimes\zeta^\hbar_{\bar\lambda})
\delta^\hbar =\Sigma q^t(\zeta^\hbar_{\bar\lambda}\otimes
\xi^\hbar_\mu)\delta^\hbar
$$
as morphisms $M^\hbar\to V_\mu\otimes V_{\bar\lambda}$. Recall that
$\delta^\hbar$ is induced by the morphisms $m^\hbar$ defined by
\eqref{em}. Since $\xi^\hbar_\mu$ and $\zeta^\hbar_{\bar\lambda}$
are defined by the morphisms $\tr^{\eta,\hbar}_{0,\mu}$ and
$\tr^{\nu,\hbar}_{\lambda,0}$, respectively, by equality
\eqref{etrm} it suffices to show the following equality of
endomorphisms of $V_\mu\otimes V_{\bar\lambda}$:
$$
q^{-(\lambda,\mu)} m^\hbar_{0,\lambda,\mu,-\lambda} =\Sigma q^t
m^\hbar_{\lambda,0,-\lambda,\mu}.
$$
This is immediate by definition \eqref{em}, since the associator
$\Phi$ acts trivially on a tensor product of three modules if at
least one module is trivial. \ep

\end{document}